\documentclass[12pt]{article}

\usepackage{geometry}                
\usepackage{amsmath,amsfonts,amssymb,amsthm,amscd,enumerate,bm}
\usepackage[all]{xy}
\usepackage{graphicx}
\usepackage{epstopdf} 
\usepackage{enumitem} 
\usepackage{hyperref}

\usepackage[dvipsnames]{xcolor}

\usepackage{ulem}

\title{A Dolbeault--Dirac Spectral Triple for the $B_2$-Irreducible Quantum Flag Manifold}

\author{Fredy D\'iaz Garc\'ia$^{1,3}$ \and R\'eamonn \'O Buachalla$^2$ \and  
Elmar Wagner$^3$  
}

\date{ \normalsize
$^1$Centro de Ciencias Matem\'aticas, UNAM, Morelia, Mexico\\
E-mail address: {\it lenonndiaz@gmail.com}\\[2ex]

$^2$Mathematical Institute of Charles University, Prague, Czech Republic \\
E-mail address: {\it obuachalla@karlin.mff.cuni.cz} \\[2ex]

$^3$Instituto de F\'isica y Matem\'aticas, UMSNH, Morelia, Mexico\\
E-mail address: {\it elmar.wagner@umich.mx}
}           

\newtheorem{thm}{Theorem}
\newtheorem{prop}[thm]{Proposition}
\newtheorem{lem}[thm]{Lemma}
\newtheorem{cor}[thm]{Corollary}

\theoremstyle{definition}

\newtheorem{rem}[thm]{Remark}

\newcommand{\rnc}[2]{\renewcommand{#1}{#2}}
\rnc{\[}{\begin{equation}}
\rnc{\]}{\end{equation}}

\newcommand{\Z}{\mathbb{Z}}
\newcommand{\N}{\mathbb{N}}
\newcommand{\R}{\mathbb{R}}
\newcommand{\C}{\mathbb{C}}
\newcommand{\hsp}{{\hspace{-1pt}}}
\newcommand{\hs}{{\hspace{1pt}}}

\newcommand{\vare}{\varepsilon} 
\newcommand{\ha}{\alpha} 
\newcommand{\hb}{\beta} 
 
\newcommand{\ho}{\omega} 
\newcommand{\hu}{\upsilon} 
\newcommand{\hl}{\lambda} 
\newcommand{\hO}{\Omega}

\newcommand{\D}{\mathcal{D}}
\newcommand{\hH}{\mathcal{H}}
\newcommand{\cO}{\mathcal{O}}
\newcommand{\cU}{{\mathcal{U}}}
\newcommand{\K}{{\mathcal{K}}}

\newcommand{\A}{{\mathcal{A}}}

\newcommand{\cC}{{\mathcal{C}}}

\newcommand{\E}{{\mathcal{E}}}
\newcommand{\F}{{\mathcal{F}}}

\newcommand{\w}{w}

\newcommand{\db}{\bar{\partial}}

\newcommand{\lin}{{\mathrm{span}}}
\newcommand{\id}{{\mathrm{id}}}

\newcommand{\spec}{\mathrm{spec}}
 
\newcommand{\ran}{\mathrm{im}}

\newcommand{\inv}{{\mathrm{inv}}}
\newcommand{\pr}{{\mathrm{pr}}}

\newcommand{\rmB}{\mathrm{B}}

\newcommand{\half}{\mbox{$\frac{1}{2}$}}
\newcommand{\ot}{\otimes} 
\newcommand{\ut}{\,\underline{\otimes}\,}

\newcommand{\lra}{\longrightarrow}
\newcommand{\lla}{\longleftarrow}
\newcommand{\rtr}{\hsp\triangleleft\hsp}
\newcommand{\ltr}{\hsp\triangleright\hsp}
\newcommand{\btr}{\hs\raisebox{0.2ex}{\mbox{${\scriptstyle{\blacktriangleright}}$}}\hs}

\newcommand{\ip}[2]{\langle{#1},{#2}\rangle}

\newcommand{\msum}[2]{\underset{{#1}}{\overset{{#2}}{\mbox{$\sum$}}}}
\newcommand{\ebox}[1]{\mbox{$#1$}} 
\newcommand{\efrac}[2]{\mbox{$\frac{#1}{#2}$}}
\newcommand{\up}[1]{{}^{#1}} 
\newcommand{\csqrt}[1]{\sqrt{\phantom{c^k}}\hspace{-11pt}{ #1}}

\newcommand{\uhr}[1]{\!\! \upharpoonright_{\!#1}}
\newcommand{\Uqg}{\cU_q(\mathrm{so}(5))}
\newcommand{\OGq}{\cO(\mathrm{SO}_q(5))}
\newcommand{\Uql}{\cU_q(\mathfrak{l})}

\newcommand{\cc}{(\!(} 
\newcommand{\CC}{)\!)}

\begin{document}
\maketitle

\begin{abstract} 
The quantum version of the Bernstein--Gelfand--Gelfand resolution  
is used to construct a Dolbeault--Dirac operator 
on the anti-holomorphic forms of the Heckenberger--Kolb calculus 
for the $B_2$-irreducible quantum flag manifold. 
The spectrum and the multiplicities of the eigenvalues 
of the Dolbeault--Dirac operator are computed. 
It is shown that this construction yields an  
equivariant, even, $0^+$-summable spectral triple. 
\end{abstract}

\section{Introduction}

Soon after the simultaneous and independent emergence 
of quantum groups and noncommutative geometry in the 1980's, 
mathematicians have tried to relate these two theories to each other. 
In fact, the question of how to reconcile the theory of spectral triples 
with Drinfeld--Jimbo quantum groups is one of the major open problems 
in noncommutative geometry. The best-studied example is the spectral triple 
on the standard Podle\'s constructed by D\c{a}browski and Sitarz \cite{DS}. 
Later it was shown in \cite{SW} that this Dirac operator 
yields a commutator representation of Podle\'s' distinguished 2-dimensional covariant 
first order differential calculus \cite{P}. 
Furthermore, this spectral triple satisfies 
most of Connes' conditions for a noncommutative spin geometry, 
albeit in a modified manner \cite{RS, Wa}. Most notably, 
the eigenvalues of the Dirac operator grow exponentially 
so that the spectral triple has spectral dimension~0. 

Other attempts to reconcile Connes' noncommutative geometry with quantum groups 
led to equivariant isospectral Dirac operators on all Podle\'s spheres \cite{DDLW}, 
on quantum SU(2) \cite{DLSSV}, 
on the 4-dimensional orthogonal quantum sphere \cite{DDL}, 
and on all compact quantum groups \cite{NT}, 
but none of these equivariant spectral triples has been related to a 
finite-dimensional covariant differential calculus. Moreover, basic conditions 
on a noncommutative spin geometry, like the first order condition 
and the real structure, fail or can be satisfied only up to smoothing operators, 
which makes the standard Podle\'s sphere even more outstanding. 

What has become increasingly clear in recent years is that quantum flag manifolds 
have a major role to play to resolving this problem. These quantum homogeneous spaces, 
which $q$-deform the coordinate rings $\mathcal{O}(G/L_S)$ of the classical flag manifolds, 
have a noncommutative geometric structure much closer to the classical situation 
than quantum groups themselves. For those quantum flag manifolds 
of irreducible type, Heckenberger and Kolb showed that $\mathcal{O}_q(G/L_S)$ 
comes endowed with an essentially unique $q$-defomed de Rham complex \cite{HK,HK1}  
such that the homogeneous components have the same dimension as in the classical case, 
thus directly extending Podle\'s' calculus construction 
for the standard Podle\'{s} sphere. 
 
The existence of such a canonical deformation is one of the most important results 
in the noncommutative geometry of quantum groups, and has served as a solid base for 
the construction of Dirac operators for these quantum spaces. 
The first to notice it was Kr\"ahmer. 
In \cite{UKFlagDirac}, he introduced an influential algebraic 
Dirac operator for the irreducible quantum flag manifolds giving 
a commutator realization of their Heckenberger--Kolb calculi.  
However, the question of a compact resolvent of the Dirac operator remained 
unanswered as the computation of the spectrum is a difficult task. 

The problem of computing the spectrum was solved only in the 
particular case of $\mathcal{O}_q(\C\mathrm{P}^2)$ by D\c{a}browski,  D'Andrea and Landi 
\cite{SISSACP2}. As for the standard Pod\-le\'s sphere $\mathcal{O}_q(\C\mathrm{P}^1)$, 
the spectrum had exponential growth. 
Then D\c{a}browski and D'Andrea proved in \cite{SISSACPn}
that the eigenvalues of the Dolbeault--Dirac operator 
for any quantum projective space grow exponentially. 
This was achieved by relating the square of the Dolbeault--Dirac operator 
to a Casimir operator of the corresponding quantized universal enveloping algebra 
$\cU_q(\mathfrak{g})$ and without explicitly computing the spectrum.  
Later Krähmer and Tucker-Simmons gave a new procedure for the construction of quantum 
Clifford algebras and Dolbeault--Dirac operators on quantized symmetric spaces \cite{KTS}. 
However, Matassa \cite{MatassaCliff} showed that there does not hold in general a 
Parthasarathy-type formula which relates the square of these operators to a quadratic 
Casimir. Despite that, Matassa obtained an asymptotic (up to bounded operators) 
Parthasarathy-type formula for the quantum projective spaces \cite{MatassaCPn} 
and for the quantum Lagrangian Grassmannian of rank 2 \cite{Ma}, following the 
construction of Krähmer and Tucker-Simmons 
and verifying in that way the compact resolvent condition. 


In recent years, an alternative approach to the construction of Dolbeault--Dirac operators 
was initiated in \cite{OB}, where the notion of a noncommutative K\"ahler structure 
was introduced. A covariant K\"ahler structure for quantum projective space 
was constructed and shown to be the unique such structure. This construction 
was later extended by Matassa \cite{MarcoConj} to all the irreducible quantum flag manifolds, 
for all but a finite number of values of $q$. The rich structure associated to 
K\"ahler structures allows one to construct associated metrics, Hilbert space completions, 
and to verify all the requirements of a spectral triple with the crucial exception 
of the compact resolvent condition, the most challenging spectral triple axiom. 
This question was answered for the case of 
quantum projective space  $\mathcal{O}_q(\mathbb{CP}^n)$ in \cite{DOS1}. 
However, the proof drew heavily on the multiplicity free nature 
of the anti-holomorphic forms of quantum projective space and has proved difficult 
to extend to other examples.


In this paper we adopt yet another approach for the construction of a 
Dolbeault--Dirac operator related to the Heckenberger--Kolb calculus. 
It is based on a quantum group version of the Bernstein--Gelfand--Gelfand resolution for 
quantized irreducible flag manifolds established by Heckenberger and Kolb 
in \cite{HK2,HK3}. Taking the locally finite dual of the quantum 
Bernstein--Gelfand--Gelfand resolution gives an explicit description of the 
anti-holomorphic complex $(\Omega^{(0,\bullet)},\overline{\partial})$ 
by means of associated quantum vector bundles 
and actions of elements from the quantum tangent space \cite{HK3}. 
Using the Haar state of the compact quantum group $\mathcal{O}_q(G)$, 
we define an inner product on $\Omega^{(0,\bullet)}$ so that 
we can consider the adjoint operator $\overline{\partial}^\dagger$ and  
hence the $q$-deformed Dolbeault--Dirac operator 
$D = \overline{\partial} + \overline{\partial}^\dagger$. 
The bounded commutator condition follows from the fact that the 
quantum tangent space acts on the invariant algebra $\mathcal{O}_q(G/L_S)$  
by derivations satisfying the (non-twisted) Leibniz rule. 

Computing the spectrum of the $q$-deformed Dolbeault--Dirac operator 
in the general case remains a difficult task. In this paper, we accomplish it 
for the lowest dimensional irreducible quantum flag manifold 
of the B-series, i.\,e., for a quantum version of 
$\mathrm{SO}(5)/(\mathrm{SO}(2)\times \mathrm{SO}(3))$. 
The difficulties are due to the lack of an explicit description of the Haar state 
and of a convenient orthonormal basis, e.\,g.\ a quantum version of the 
Gelfand--Tsetlin basis, which only exists for the A-series. 
We circumvent these problems by expressing the actions of 
$\overline{\partial}$ and $\overline{\partial}^\dagger$ 
in purely algebraic terms by elements from $\Uqg$ and working solely 
on highest weight vectors. By equivariance, these computation 
determine completely the eigenvalues and multiplicities 
of the  Dolbeault--Dirac operator. It turns out that, similar to the 
quantum projective spaces, the eigenvalues grow exponentially. 
This offers strong support for the conjecture that the eigenvalues of 
the Dolbeault--Dirac operators of all the irreducible quantum flag manifolds 
have exponential growth.

The paper is organized as follows. In the Section \ref{DN},  
we describe the Hopf \mbox{*-al}\-ge\-bras \hs$\Uqg$ and $\OGq$, 
the Levi factor $\Uql\subset\Uqg$, 
the invariant algebra $B:= \OGq^{\inv(\Uql)}$, 
the representation theory of $\Uqg$ and $\Uql$, 
and some useful relations of the generators of $\OGq$. 
In the Section \ref{secBGG}, we give an explicit description 
of the anti-holomorphic $k$-forms and the $\overline{\partial}$-operator 
by determining the Bernstein--Gelfand--Gelfand resolution for our example. 
Our strategy to compute the spectrum of the equivariant Dolbeault--Dirac operator 
is to describe its action on highest weight vectors. As a prerequisite, 
we study in Section \ref{sec-bl} the branching rules, i.\,e., 
the problem of determining the multiplicity with which an 
irreducible representation of $\Uql$ of given highest weight 
occurs as a subrepresentation of \hs$\Uqg$. Since in our case the 
highest weight vectors of \hs$\Uql$ are given by the kernel of a unique ladder operator, 
we solve this problem by applying Kostant’s multiplicity formula. 
Then we present an explicit description of these highest weight vectors in terms 
of the generators of $\OGq$. 

The $q$-deformed Dolbeault--Dirac operator 
$D = \overline{\partial}+\overline{\partial}^\dagger$ is introduced in Section~\ref{secD}. 
For the definition of the Hilbert space adjoint 
$\overline{\partial}^\dagger$ of $\overline{\partial}$, we need to equip 
the space of anti-holomorphic $k$-forms $\Omega^{(0,k)}$ with an inner product. 
This is done by using the Haar state to define an inner product on $\OGq$ and 
considering $\Omega^{(0,k)}$ as a subspace of $\OGq$ for $k=0,3$, 
and as a subspace of $\OGq\otimes\C^3$ for $k=1,2$. 
Since $\overline{\partial}$ can be expressed by the left action 
of elements from the quantum tangent space, $D$ is equivariant with respect to 
the right $\Uqg$-action on $\OGq$. This allows us to reduce the computation of 
the spectrum to the action of $D$ on highest weight vectors. 
To simplify the computations even more, we establish in Lemma \ref{equiv} 
a unitary isomorphism between the space the highest weight vectors in $\Omega^{(0,k)}$ 
and the corresponding highest weight vectors in $\OGq$, 
and consider then the unitarily equivalent action of 
$\overline{\partial}$ between subspaces of $\OGq$ with dimension at most 2. 
Although taking Hilbert space adjoints with respect to the inner product defined 
by the Haar state coincides with the involution of \hs$\Uqg$, we cannot simply apply 
the involution to obtain an expression  for $\partial^\dagger$ because these 
operators don't leave the subspace in $\OGq$ corresponding to $\Omega^{(0,k)}$ invariant. 
Therefore it is necessary to apply subsequently an orthogonal projection 
onto the correct subspaces. In order to avoid the computation of inner products 
with respect to the Haar state, we use the fact that the positive and negative 
simple root vectors act as ladder operators to give an algebraic description of 
these projections. From the equivariance of 
$\overline{\partial}$ and $\overline{\partial}^\dagger$, it follows immediately 
that the highest weight vectors belonging to weight spaces with multiplicity 1 are 
eigenvectors of the Laplacians $\overline{\partial}^\dagger_k \overline{\partial}_k$
and $\overline{\partial}_k \overline{\partial}^\dagger_k$. 
Now the corresponding eigenvalues can be computed directly 
by applying the unitarily equivalent operators, 
expressed in terms of elements from $\Uqg$, 
on the highest weight vectors given by polynomials in the generators of $\OGq$. 
We present a few  illustrative cases of these calculations to show our 
strategy of avoiding the use of complicated commutation rules in $\OGq$, 
the other calculations can be found in~\cite{D}. 
For the calculation of the eigenvalues on the 0-forms and on the top forms, 
we relate the Laplacians to a Casimir operator and employ a 
Parthasarathy-type formula. The  knowledge of the eigenvalues 
of the Laplacians on the spaces of highest weight vectors with multplicity 1 
suffices to deduce the complete spectrum  
and the multiplicities of the eigenvalues of Dolbeault--Dirac operator. 

Our main result is the last theorem, where we prove that the 
Dolbeault--Dirac operator on the anti-holomorphic forms yields an 
equivariant, even, $0^+$-summable spectral triple. We are confident 
that the techniques and strategies presented in this paper can be used 
to compute further examples of spectral triples on 
irreducible quantum flag manifolds. 
Of course, the long-standing goal is to understand quantization in physics. 

\section{Preliminaries} \label{DN} 

We begin by recalling some general facts about the complex simple Lie algebra $\mathrm{so}(5)$. 
Since $\mathrm{so}(5)$ has rank 2, 
a chosen Cartan subalgebra $\mathfrak{h}\subset \mathrm{so}(5)$ 
is 2-dimensional. Let $\{\ha_1,\ha_2\}$ be a set of simple roots corresponding to 
$\mathfrak{h}$ 
such that $(\ha_1,\ha_1)= 2$ and $(\ha_2,\ha_2)=1$, where $(\cdot,\cdot)$ 
denotes the symmetric bilinear 
on the dual $\mathfrak{h}'$ of $\mathfrak{h}$ induced by the Killing form. 
This symmetric bilinear can be deduced from the Cartan matrix of $\mathrm{so}(5)$:  
\[  \label{cm} 
\left(  \efrac{2(\ha_i,\ha_j)}{(\ha_j,\ha_j)}  \right)_{i,j=1,2}   
\ = \  \left( \begin{array}{rr}  2 & -2\\
                                \!\!-1 & 2 \end{array}\right). 
\]
A basis of fundamental weights $\{\omega_1,\omega_2\}\subset \mathfrak{h}'$ 
satisfying $\frac{2(\omega_i,\ha_j)}{(\ha_j,\ha_j)} = \delta_{ij}$ (Kronecker delta) 
is given by $\omega_1 = \ha_1+\ha_2$ and $\omega_2 = \frac{1}{2}(\ha_1+2\ha_2)$. 
We will use Cartesian coordinates to describe the weight lattice 
$P:=\Z \hs\omega_1 + \Z\hs \omega_2$, i.\,e., $(n_1,n_2) \in \Z^2$ stands for 
$n_1\hs \omega_1 + n_2\hs \omega_2 \in P$. 
In these coordinates, the set of positive roots reads 
as follows: 
\[  \label{Rp} 
R^+:=\{\ha_1\!=\!(2,-2),\ \,\ha_2\!=\!(-1,2), \ \,
\ha_1\!+\!\ha_2\!=\!(1,0), \ \,\ha_1\!+\!2\ha_2\!=\!(0,2)\}.  
\] 
We write $P^+:= \N_0 \hs \omega_1 + \N_0\hs\omega_2 \cong \N_0^2$ for the set dominant weights. 
It contains the fundamental Weyl chamber $P^{++} := \N \hs \omega_1 + \N \hs\omega_2$. 
A weight $\mu\in P$ is said to be higher than $\nu\in P$ if 
$\mu -\nu \in \N_0\hs \ha_1 +\N_0 \ha_2$. 

The Weyl group $W$ is generated the by the reflections $\w_\ha :  P \to P$, $\ha\in R^+$, 
where $\w_\ha (\mu) :=\mu - 2\frac{(\mu,\ha)}{(\ha,\ha)} \ha$. 
For instance, the reflection $\w_{\ha_2}$ is given by 
\[ \label{r2} \w_{\ha_2} (x,y)  = (x+y, -y), \quad (x,y)\in P.\]

Throughout the paper,  $q$ stands for a real number from the interval $(0,1)$. 
Set $q_1:= q$, \,$q_2:= \sqrt{q}$, and  
\[
[n]_1 := \mbox{$\frac{q_1^n - q_1^{-n}}{q_1- q_1^{-1}}$} 
= \mbox{$\frac{q^n - q^{-n}}{q- q^{-1}}$}, \quad 
[n]_2 := \mbox{$\frac{q_2^n - q_2^{-n}}{q_2- q_2^{-1}}$}
= \mbox{$\frac{q^{n/2} - q^{-n/2}}{q^{1/2}- q^{-1/2}}$}, \quad n\in\R. 
\] 
The quantized enveloping algebra $\Uqg$ is defined as the complex *-algebra 
generated by $K_1^{\pm 1}, K_2^{\pm 1}, E_1,E_2, F_1,F_2$ with relations \cite{KS} 
\begin{align}
& K_1E_1 = q_1^2 E_1K_1,\ \ K_1E_2 = q_1^{-1} E_2K_1,\ \ 
K_1F_1 = q_1^{-2} F_1K_1,\ \ K_1F_2 = q_1 F_2K_1, \nonumber \\
\begin{split}
& K_2E_1 = q_2^{-2}E_1K_2,\ \ K_2E_2 = q_2^{2} E_2K_2,\ \ 
K_2F_1 = q_2^2 F_1K_2,\ \ K_2F_2 = q_2^{-2}F_2K_2,  \label{EFFE}  \\
&K_iK_j\!=\! K_jK_i, \ \ K_iK_i^{-1} \!=\! K_i^{-1}K_i=1,  \ \ 
E_1F_2\!=\! F_2E_1, \ \ E_2F_1\!=\! F_1E_2,\ \ i,j\!=\!1,2,
\end{split} \\ \nonumber 
&E_1F_1\! -\!F_1 E_1= \mbox{$\frac{1}{q_1\hsp -\hsp q_1^{-1}}$}(K_1\! -\!K_1^{-1}), 
E_2F_2\! -\! F_2 E_2= \mbox{$\frac{1}{q_2\hsp  -\hsp q_2^{-1}}$}(K_2\! -\!K_2^{-1}),  
\end{align} 
the quantum Serre relations
\begin{align}  \begin{split} \label{Sr1} 
0&= E_2^3E_1\! -\! [3]_2  E_2^2E_1E_2\hsp +\hsp [3]_2  E_2E_1E_2^2 \!- \! E_1E_2^3, \quad 
0= E_1^2E_2 \!- \![2]_1 E_1E_2E_1\! +\!E_2E_1^2,  \\
0&= F_2^3F_1\hsp - [3]_2  F_2^2F_1F_2 + [3]_2  F_2F_1F_2^2\hsp - F_1F_2^3,\quad 
0= F_1^2F_2 \!- \![2]_1 F_1F_2F_1\! +\!F_2F_1^2, 
\end{split} 
\end{align} 
and involution 
\begin{align}
\begin{split} \label{in} 
&K_1^* = K_1, \quad  E_1^*  = q_1K_1F_1,   \quad 
F_1^* = q_1^{-1}E_1K_1^{-1}, \\
&K_2^* = K_2,  \quad  E_2^* =q_1K_2F_2, \quad 
F_2^* = q_1^{-1}E_2K_2^{-1}.
\end{split}
\end{align}
It is a Hopf *-algebra with coproduct $\Delta$, counit $\vare$ and antipode $S$ determined by 
\begin{align}  \nonumber 
&\Delta(K_i^{\pm 1}) = K_i^{\pm 1}\!\otimes\! K_i^{\pm 1},  
\ \ \Delta(E_i) =  E_i\!\otimes\! K_i + 1\!\otimes\! E_i, \ \ 
\Delta(F_i) = F_i\!\otimes\! 1+ K_i^{-1}\!\otimes\! F_i, \label{DE} \\
&\vare(K_i)=1,  \quad  \vare(E_i)=\vare(F_i) =0, \\
&S(K_i)=K_i^{-1}, \quad  S(E_i)=-E_iK_i^{-1}, \quad  S(F_i)= -K_iF_i, \quad i=1,2.  \nonumber 
\end{align}

Note that $q_1$ and $q_2$  were chosen in such a way that 
the powers in the commutation relations of 
$K_1$ and $E_j$ in \eqref{EFFE}  
coincide with the Cartan matrix \eqref{cm}.  
Since the Cartan matrix of $\mathrm{so}(5)$ 
is the transpose of the Cartan matrix of $\mathrm{sp}(4)$, 
it follows from the definitions in \cite[Section 6]{KS} that $q$-deformations of 
$\cU(\mathrm{so}(5))$ and  $\cU(\mathrm{sp}(4))$ are isomorphic. 
Applying the *-isomorphism 
$\Phi: \cU_{q_2}(\mathrm{sp}(4)) \rightarrow \cU_{q}(\mathrm{so}(5))$ given by 
\begin{align*}
&\Phi(K_1) = K_2^{-1}, \ \ \Phi(E_1) = q_2F_2,\ \ \Phi(F_1) = q_2^{-1}E_2, \\
&\Phi(K_2) = K_1^{-1},\ \   \Phi(E_2) = q_2F_1, \ \ \Phi(F_2) = q_2^{-1}E_1,
\end{align*}
a central element $\cC$ in $\Uqg$ can be read off from the quadratic Casimir operator 
in \cite[Proposition 6.6]{Ma}. Defining 
\begin{align*} 
& \K_1 := K_2^{-2} K_1^{-1}, \quad \K_2 := K_1^{-1},\quad \F_1:=F_2,\quad \F_4 := F_1, \\ 
& \F_2:= F_2^2 F_1\hsp - \hsp (1+q^{-1}) F_2F_1F_2\hsp + \hsp q^{-1} F_1 F_2^2, \quad 
\F_3 := F_2F_1 \hsp - \hsp  q^{-1} F_1 F_2,   
\end{align*} 
we get 
\begin{align}  \nonumber
\cC &= \efrac{1}{(q_2-q_2^{-1})^2}(q^{-2}\K_1^{-1} + q^{-1}\K_2^{-1} + q\K_2 + q^{2}\K_1)\\ 
\nonumber
&+ (q_2^{-3} \K_1^{-1} \hsp + \hsp q_2^{3} \K_2) \F_1^*\F_1 +\K_1^{-1} \F_2^*\F_2 
+ (q_2^{-3} \K_1^{-1}  \hsp + \hsp q_2^3 \K_2^{-1}) \F_3^*\F_3 
+ q^{-1}[2]_2^2 \K_2^{-1}\F_4^*\F_4 \\ \nonumber
& + (1\! -\! q) [2]_2^2  \K_1^{-1} (\F_1^*\F_3^* \F_2 +\F_2^*\F_3\F_1) 
+ (q^{-2} \!-\! q^{-1}) [2]_2  \K_1^{-1} (\F_1^*\F_4^* \F_3 +\F_3^*\F_4\F_1) \\ 
& + (1\!-\! q)^2 \K_1^{-1} \F_1^*\F_3^*\F_3\F_1 
+ (1 \!-\! q)^2q_2^{-5}  [2]_2^2 \K_1^{-1} \F_1^*\F_4^* \F_4\F_1.           \label{C}   
\end{align}  
Another version of a quadratic Casimir operator is obtained by applying the 
antipode to $\cC$. 
Setting 
\begin{align} 
\begin{split} \label{Es} 
& \K_1 = K_2^{-2} K_1^{-1}, \quad \K_2 = K_1^{-1},\quad  \E_1=E_2, \quad \E_4 = E_1,  \\ 
& \E_2=  E_1E_2^2\hsp - \hsp (1+q) E_2E_1E_2\hsp + \hsp q E_2^2 E_1 , \quad 
\E_3 = E_1E_2 \hsp - \hsp  q E_2 E_1.  
\end{split} 
\end{align} 
the central element $S(\cC)$ can be written as 
\begin{align}\nonumber
&S(\cC) = \efrac{1}{(q_2-q_2^{-1})^2}
(q^{-2} \K_1 + q^{-1} \K_2 + q \K_2^{-1} + q^{2} \K_1^{-1}) \\ \nonumber
& + q^{-1}(q_2^{-3} \K_1\! + \hsp q_2^{3} \K_2^{-1} ) \E_1^*\E_1 + q^{-4}\K_1 \E_2^*\E_2 
+ q^{-3}(q_2^{-3} \K_1 \! + \hsp q_2^3 \K_2) \E_3^*\E_3 + q^{-1}[2]_2^2 \K_2\E_4^*\E_4 \\ 
\nonumber
& - (q^{-4}\! -\! q^{-3}) [2]_2^2  \K_1 (\E_1^*\E_3^* \E_2 +\E_2^*\E_3\E_1) 
- (q^{-4} \!-\! q^{-3}) [2]_2  \K_1 (\E_1^*\E_4^* \E_3 +\E_3^*\E_4\E_1) \\ 
&  + (q^{-2}\!-\! q^{-1})^2 \K_1 \E_1^*\E_3^*\E_3\E_1 
+ (q^{-2} \!-\! q^{-1})^2q_2  [2]_2^2 \K_1\E_1^*\E_4^* \E_4\E_1  .      \label{SC}  
\end{align} 

Let $\Uqg^{\circ}$ denote the dual Hopf algebra of $\Uqg$ \cite[Section 1.2.8]{KS}. 
Then $\OGq \subset \Uqg^{\circ}$ is the Hopf subalgebra generated 
by the matrix coefficients of the 5-dimensional irreducible $\Uqg$-representation 
of highest weight $(1,0)$. Denoting the matrix coefficients  by $u^i_j$, \,$i,j=1,\ldots,5$,  
they satisfy the relations 
\begin{align}   \begin{split}  \label{Rm}
&\sum_{k,l=1,...,5} \hspace{-8pt} R^{ij}_{kl}u^k_mu^l_n \ = \!\!
\sum_{k,l=1,...,5}\hspace{-8pt}  u^j_ku^i_lR^{lk}_{mn}, \ \ 
i,j,m,n=1, \ldots,5, \quad \D_q=1, \\
&\sum_{i,j,k =1,...,5}\hspace{-12pt} C^{\,i}_j (C^{-1})^{k}_m  u^n_{i} u^{k}_{j} \ = \!\!
\sum_{i,j,k =1,...,5}\hspace{-12pt}  C^{n}_i (C^{-1})^{j}_k  u^j_{i} u^{k}_{m}  
= \delta_{nm}, \ \ m,n=1, \ldots,5,
\end{split} 
\end{align} 
with the R-matrix coefficients $R^{ij}_{kl} \in \R$, the C-matrix coefficients $C^{\,i}_j\in \R$
and the quantum determinant $\D_q \in \mathrm{alg}\{ u^k_j : j,k =1,...,5\} $
given in \cite[Section 9.3]{KS}, 
and $\delta_{nm}$ denotes the Kronecker delta. 
$\OGq$ is a Hopf *-algebra, where the structure maps 
are determined by 
$$
\Delta(u^k_l) = \mbox{$\sum$}_i\, u^k_i \ot u^i_l, \ \ \vare(u^k_l) = \delta_{kl}, \ \ 
S(u^k_l) = \mbox{$\sum$}_{i,j} C^{\,k}_i (C^{-1})^{j}_l u^j_i, \ \ u^{k*}_l = S(u^l_k). 
$$

There exist a left action 
$\Uqg \otimes \OGq\ni X\ot a \mapsto X\ltr \hs a  \in \OGq$  
and a right action 
$ \OGq\otimes\Uqg\ni a\ot X \mapsto a\rtr X  \in \OGq$
such that $\OGq$ is a left and a right $\Uqg$-module *-algebra. 
Here, module *-algebra  means that these actions satisfy
\begin{align} \begin{split}  \label{ma} 
X\ltr (ab) &= (X_{(1)}\ltr a)(X_{(2)}\ltr b), \qquad (X\ltr a)^*= S(X)^*\ltr a^*, \\ 
(ab)\rtr X &= (a\rtr X_{(1)})(b\rtr X_{(2)}), \qquad  (a\rtr X)^*= a^* \rtr S(X)^* 
\end{split} 
\end{align} 
for all $a,b\in \OGq$ and $X\in \Uqg$.  
We employ Sweedler's notation for co\-pro\-ducts, 
i.\,e., $\Delta(X)= X_{(1)}\!\ot\! X_{(2)}$. 
As a consequence of \eqref{ma}, it suffices to specify the actions on generators, 
which are given as follows: 
\begin{align}
&K_1 \ltr u^k_1 = q^{-1} u^k_1, \ \ K_1 \ltr u^k_2 = q u^k_2,\ \ K_1 \ltr u^k_3 = u^k_3,\ \ 
K_1 \ltr u^k_4 = q^{-1} u^k_4, \ \ K_1 \ltr u^k_5 = q u^k_5, \nonumber \\ 
&E_1\ltr u^k_1 = u^k_2,\ \ E_1\ltr u^k_4= - u^k_5,\ \ F_1\ltr u^k_2 = u^k_1,\ \ 
F_1\ltr u^k_5= - u^k_4,\nonumber \\ 
&K_2 \ltr u^k_1 =  u^k_1, \ \ K_2 \ltr u^k_2 = q^{-1} u^k_2,\ \ 
K_2 \ltr u^k_3 = u^k_3,\ \     \label{lac}
K_2 \ltr u^k_4 = q u^k_4, \ \ K_2 \ltr u^k_5 = u^k_5, \\ 
&E_2 \ltr u^k_2 = [2]_2^{1/2}u^k_3, \ \ E_2 \ltr u^k_3= -q_2 [2]_2^{1/2}u^k_4, \nonumber \\ 
&F_2 \ltr u^k_3 = [2]_2^{1/2} u^k_2, \ \ 
F_2 \ltr u^k_4= -q_2^{-1}[2]_2^{1/2} u^k_3, \nonumber         
\end{align}
\begin{align}
&u_k^1\rtr K_1 = q^{-1} u_k^1, \ \ u_k^2\rtr K_1  = q u_k^2,\ \ u_k^3\rtr K_1  = u_k^3,\ \ 
u_k^4\rtr K_1 = q^{-1} u_k^4, \ \ u_k^5 \rtr K_1 = q u_k^5, \nonumber \\ 
&u_k^2\rtr E_1 = u_k^1,\ \  u_k^5\rtr E_1 = - u_k^4,\ \ 
u_k^1\rtr F_1 = u_k^2,\ \ u_k^4\rtr F_1 = - u_k^5,\nonumber \\ 
&u_k^1 \rtr K_2 =  u_k^1, \ \ u_k^2 \rtr K_2 = q^{-1} u_k^2,\ \ u_k^3 \rtr K_2 = u_k^3,\ \ 
\label{rac}
u_k^4 \rtr K_2 = q u_k^4, \ \ u_k^5 \rtr K_2 = u_k^5, \\ 
&u_k^3 \rtr E_2 = [2]_2^{1/2}u_k^2, \ \ u_k^4 \rtr E_2 = -q_2 [2]_2^{1/2}u_k^3, \nonumber \\ 
&u_k^2\rtr F_2 = [2]_2^{1/2} u_k^3, \ \ u_k^3 \rtr F_2 = -q_2^{-1}[2]_2^{1/2} u_k^4,  \nonumber
\end{align}
for $k=1,\ldots,5$, and zero in all other cases. 

The Levi factor $\Uql \subset \Uqg$ is the Hopf-*-subalgebra generated by 
$K_1^{\pm 1}$, $K_2^{\pm 1}$, $E_2$, $F_2$. 
The subalgebra $B\subset\OGq$ of $\Uql$-invariant elements, i.\,e.,  
\[   \label{B}
B:= \OGq^{\inv(\Uql)}\!=\{ a\in\OGq : X\ltr a= \vare(X)\hs a 
\text{ for all } X\hsp\in\hsp \Uql\}, 
\] 
is viewed as the coordinate ring of the 
quantum quadratic SO(5)/(SO(2)$\times$SO(3)). 

It is known that $\OGq$ admits a Haar state, that is, 
a positive, faithful and $\Uqg$-invariant state
\[ \label{Haar}
h: \OGq \lra \C,
\]
where $\Uqg$-invariance means that 
\[ \label{hinv} 
h(X\ltr a)= h(a\rtr X) = \vare(X)\hs h(a) , \quad X\in \Uqg, \ \ a\in \OGq. 
\] 
As $h$ is positive and faithful, it equips $\OGq$ with an inner product by setting 
\[    \label{iph} 
\ip{a}{b}_h:= h(a^*b), \qquad a,b\in \OGq. 
\] 
The left regular representation $\pi_L$ and the right regular representation 
$\pi_R$ of \hs$\Uqg$ 
on $\OGq$ are defined by 
\[ \label{regrep}
\pi_L(X)(a) := X\ltr a,\ \ \pi_R(X)(a) := a \rtr S^{-1}(X), \ \ X\in \Uqg, \ \, a\in \OGq, 
\] 
These are *-representations with respect to the inner product defined in \eqref{iph}. 
To see this, let us show the compatibility of $\pi_R$ with the involution. 
Using \eqref{ma} and \eqref{hinv}, we compute for $a,b\in\OGq$ and $X\in\Uqg$  that 
\begin{align} \nonumber
&\ip{\pi_R(X)(a)}{b}_h = h\big((a^*\rtr X^*)b\big) = \vare(S^{-1}(X_{(1)}^*))\,  
h\big((a^*\rtr X_{(2)}^*)b\big) \\
&= h\big(((a^*\rtr X_{(2)}^*)b)\rtr S^{-1}(X_{(1)}^*)\big)  \label{srep}
= h\big((a^*\rtr X_{(3)}^*S^{-1}(X_{(2)}^*))  \,(b\rtr S^{-1}(X_{(1)}^*))\big) \\  \nonumber 
&= h\big(a^* \hs(b\rtr S^{-1}(X_{(1)}^*\vare(X_{(2)}^*)))\big) 
= h\big(a^* \hs(b\rtr S^{-1}(X^*))\big) 
= \ip{a}{\pi_R(X^*)(b)}_h. 
\end{align} 
By a slight abuse of notation, we will omit the symbol $\pi_L$ and either write $X(a)$ or 
continue to use the notation $X\ltr a$ for the left regular representation.

A weight $(n,l)\in \Z\times \N_0$ defines an irreducible representation of the 
Levi factor $\Uql$ 
on a $(l\hsp +\hsp 1)$-dimensional space 
$V^{(n,l)} := \lin\{ v_{(n,l)-j(-1,2)} : j=0,1,\ldots,l\}$ by setting 
\begin{align}  \begin{split}   \label{irreps}
&K_1v_{(\mu_1,\mu_2)}=q^{\mu_1}_1v_{(\mu_1,\mu_2)}, \ \ 
K_2v_{(\mu_1,\mu_2)}=q^{\mu_2}_2v_{(\mu_1,\mu_2)}, \ \ 
E_2(v_{(\mu_1,\mu_2)}) = v_{(\mu_1-1,\mu_2+2)}, \\ 
&F_2(v_{(\mu_1,\mu_2)}) = [\half(l+\mu_2) ]_2\, [\half(l-\mu_2)+1 ]_2\, v_{(\mu_1+1,\mu_2-2)}. 
\end{split} 
\end{align} 
For instance, on 
$V^{(n,2)}=\lin\{v^{-1}\!:=\!v_{(n+2,-2)},\ v^0\!:=\!v_{(n+1,0)},\ v^1\!:=\!v_{(n,2)}\}$, 
we have 
\[  \label{rep1} 
K_1v^j = q^{n-j+1} v^j, \ \ K_2 v^j = q^j v^j, \ \ E_2 v^j = v^{j+1}, \ \ 
F_2 v^j = [2]_2  v^{j-1}, 
\]
and on $V^{(n,0)}=\lin\{u^0\!:=\!v_{(n,0)}\}$, we get 
\[ \label{rep0} 
K_1 u^0 = q^{n} u^0,\ \ K_2 u^0 = u^0, \ \ E_2 u^0 = 0, \ \ F_2 u^0 = 0. 
\]

For $(n,l)\in\Z\times\N_0$ and $(\mu_1,\mu_2)\in\Z\times \Z$, 
let $V^{(n,l)}_{(\mu_1,\mu_2)}$ denote the subspace of vectors of weight 
$(\mu_1,\mu_2)$ in $V^{(n,l)}$, i.\,e., 
$$
V^{(n,l)}_{(\mu_1,\mu_2)} := \{ v\in V^{(n,l)}: K_i(v)= q_i^{\mu_i} v,\ \,i=1,2\}. 
$$
Suppose that 
$v^{(n,2)}_{(n,2)} \oplus v^{(n-1,4)}_{(n,2)} \in V^{(n,2)}_{(n,2)} 
\oplus V^{(n\!-\!1,4)}_{(n,2)}$. 
Then,  by \eqref{irreps}, $F_2E_2(v^{(n,2)}_{(n,2)})=0$ and 
$F_2E_2(v^{(n-1,4)}_{(n,2)}) = [4]_2v^{(n-1,4)}_{(n,2)}$, thus 
\[  \label{pr2} 
\pr_{(n,2)}:= (1- \ebox{\frac{1}{[4]_2}}F_2E_2) \,:\, 
V^{(n,2)}_{(n,2)} \oplus V^{(n\!-\!1,4)}_{(n,2)} \lra V^{(n,2)}_{(n,2)}, 
\]
$\pr_{(n,2)}\big(v^{(n,2)}_{(n,2)} \oplus v^{(n-1,4)}_{(n,2)}\big)= v^{(n,2)}_{(n,2)}$, 
agrees with the (orthogonal) projection onto the first component. 

Similarly, elementary calculations using \eqref{irreps} show that 
the orthogonal projection 
$\mathrm{pr}_{(n,0)} : 
V^{(n,0)}_{(n,0)} \oplus V^{(n\hsp -\hsp 1,2)}_{(n,0)} \oplus V^{(n\hsp -\hsp 2,4)}_{(n,0)} 
\lra V^{(n,0)}_{(n,0)}$, 
\ $\mathrm{pr}_{(n,0)}\big(v^{(n,0)}_{(n,0)} \oplus v^{(n\hsp -\hsp 1,2)}_{(n,0)} 
\oplus v^{(n\hsp -\hsp 2,4)}_{(n,0)}\big)=v^{(n,0)}_{(n,0)}$, 
can be given by  
\[  \label{Pr0} 
\mathrm{pr}_{(n,0)} = (1- \ebox{\frac{1}{[2]_2}}F_2E_2 
+ \ebox{\frac{1}{[2]_2^2\hs [3]_2}}F_2^2E_2^2) 
= (1- \ebox{\frac{1}{[2]_2}}E_2F_2 + \ebox{\frac{1}{[2]_2^2\hs [3]_2}}E_2^2F_2^2). 
\]

If $a\in B$ and $x\in V^\hl_\mu\subset \OGq$,  then it follows from 
$E_2(a)=F_2(a)=0$ and $K_1(a)=K_2(a)=a$ that 
$ax$ and $xa$ belong to a representation space of the same weights. 
Moreover, using the explicit expression in \eqref{pr2} for $\pr_{(n,2)}$, 
a straightforward computation shows that 
\[  \label{pr020}
\pr_{(n,2)}(avb) = a\hs\pr_{(n,2)}(v) \hs b, \quad 
v\in V^{(n,2)}_{(n,2)} \oplus V^{(n-1,4)}_{(n,2)}\subset \OGq,\ \,a,b\in B. 
\] 

For the convenience of the reader, we finish this section with a lemma that collects 
those relations of the generators of $\OGq$ which are most frequently used 
in this paper.

\begin{lem}  \label{lem3}  
The generators $u^i_j$, \,$i,j=1,\ldots,5$, of $\OGq$ satisfy the following relations: 
\begin{align}
& u^i_1u^i_5= q^2u^i_5u^i_1, \ \ i\neq 3, \ \ u^i_lu^i_k=q u^i_ku^i_l, \ \ 
l<k, \ l\neq k', \ i\neq 3, \label{ukluRqulku} \hspace{74pt}\mbox{ } \\
& u^j_lu^i_k = u^i_ku^j_l, \ \ i<j, \ l<k, \  i\neq j', \ k\neq l', \label{uijukl} \\
& u^i_5u^j_1 = q^{-1}u^j_1u^i_5, \ \ u^i_1 u^j_5 = qu^j_5 u^i_1 + (q^2 - 1)u^i_5 u^j_1 \ \ 
i<j, \ i\neq j', \label{u15u51} \\
& u^i_k u^j_l =  u^j_l u^i_k - (q^{-1} - q)u^j_k u^i_l, \ \ 
i<j, \ k<l, \ i\neq j', \ k\neq l', \label{uij'kl'}\\
& u^i_4 u^i_2 = q^{-2}u^i_2 u^i_4 + (q^{-3} - q^{-1})u^i_1 u^i_5, \ \ i\neq 3, \label{u1142} 
\\[-24pt] \nonumber 
\end{align} 
\begin{align} 
& u^i_3u^i_3 = -[2]_2 q^{-1}u^i_2 u^i_4 - [2]_2q^{-2}u^i_1 u^i_5, \ \ i\neq 3,\\
& u^1_3u^2_3 = u^2_3u^1_3 +(q^{-1} -q)q^{-1/2} u^2_2u^1_4+(q^{-1} -q)q^{-3/2} u^2_1u^1_5, 
\label{Ru1233}\\
& u^2_3u^1_3+(q -q^{-1})u^1_3u^2_3 = u^1_3u^2_3 +(q^{-1} -q)q^{-1/2} u^1_2u^2_4
+(q^{-1} -q)q^{-3/2} u^1_1u^2_5, \label{Ru2133} \\ 
& u^1_3u^2_3 -q u^2_3u^1_3 
= (q^{1/2}\!-\!q^{-1/2})(u^1_2 u^2_4 -qu^2_2u^1_4) 
+ (q^{-1/2}\!-\!q^{-3/2})(u^1_1 u^2_5 - qu^2_1u^1_5), 
\label{u1233} \\
& u^i_ku^j_k = qu^j_ku^i_k \ \ k\neq 3, \ \, i<j, \ \ i\neq j', \label{uijukk} 
\end{align} 
\begin{align} 
& u^1_1(u^1_1 u^2_5 - qu^2_1u^1_5)=q(u^1_1 u^2_5 - qu^2_1u^1_5)u^1_1, \label{hw2}\\
& u^1_5(u^1_1 u^2_5 - qu^2_1u^1_5)=q^{-1}(u^1_1 u^2_5 - qu^2_1u^1_5)u^1_5, 
\hspace{154pt} \mbox{\ } \label{hw1} \\
& u^1_1(u^1_1 u^2_a - qu^2_1u^1_a)=(u^1_1 u^2_a - qu^2_1u^1_a)u^1_1, \label{u11Fz2} \\
& u^1_5(u^1_a u^2_5 - qu^2_au^1_5)=(u^1_a u^2_5 - qu^2_au^1_5)u^1_5, \label{uuu515} \\
& u^1_1(u^1_a u^2_5 - qu^2_au^1_5) = q^2(u^1_a u^2_5 - qu^2_a u^1_5)u^1_1,  \label{u11u25} \\
& u^1_5(u^1_1 u^2_a - qu^2_1u^1_a)=q^{-2}(u^1_1 u^2_a - qu^2_1u^1_a)u^1_5,  \label{uhw} \\ 
& u^2_5(u^1_1 u^2_a - qu^2_1u^1_a)=q^{-2}(u^1_1 u^2_a - qu^2_1u^1_a)u^2_5,   \label{u212514} \\
& u^2_5(u^1_a u^2_5 - qu^2_au^1_5) = (u^1_au^2_5 -qu^2_au^1_5)u^2_5,   \label{u25E1z2} \\
& (u^1_1u^2_a -qu^2_1u^1_a)u^1_5 
= q(u^1_1 u^2_5 - qu^2_1u^1_5)u^1_a - q^3  (u^1_a u^2_5 - qu^2_au^1_5)u^1_1,  
\label{uuu145} 
\end{align}
\begin{align}
& u^1_1 u^1_5(u^1_1 u^2_5 - qu^2_1u^1_5 ) = (u^1_1 u^2_5 - qu^2_1u^1_5)u^1_1 u^1_5, 
\label{hw3} \\
& (u^1_1u^1_a)(u^1_1u^1_5) = q^2(u^1_1u^1_5)(u^1_1u^1_a), \label{u14u15} \\ 
& (u^1_1u^1_5)\hs (u^1_au^1_5) = q^2(u^1_au^1_5)\hs (u^1_1u^1_5), \label{u15u12u15} \\
& u^1_1 u^1_5(u^1_1 u^2_a - qu^2_1u^1_a ) = q^{-2}(u^1_1 u^2_a - qu^2_1u^1_a)u^1_1 u^1_5, 
\label{u1u5Fz2} \\
& u^1_1 u^1_5 (u^1_a u^2_5 - qu^2_au^1_5) = q^2 (u^1_a u^2_5 - qu^2_au^1_5) u^1_1 u^1_5, 
\label{u1u5Ez2} \\
& (u^1_1 u^2_a - qu^2_1u^1_a)(u^1_1 u^2_5 - qu^2_1u^1_5 )
= q^{2}(u^1_1 u^2_5 - qu^2_1u^1_5 )(u^1_1 u^2_a - qu^2_1u^1_a), \label{hw4} \\
& (u^1_1 u^2_5 - qu^2_1u^1_5)\hs (u^1_a u^2_5 - qu^2_au^1_5) 
= q^2 (u^1_a u^2_5 - qu^2_au^1_5)(u^1_1 u^2_5 - qu^2_1u^1_5),   \label{z2Ez2} \\
& u^1_1 u^1_a (u^1_1 u^2_5 - qu^2_1u^1_5 )
= (u^1_1 u^2_5 - qu^2_1u^1_5 ) u^1_1 u^1_a -(1-q^2)u^1_1 u^1_5 (u^1_1 u^2_a - qu^2_1u^1_a), 
\label{hw5} \\
& u^1_au^1_5(u^1_1 u^2_5 - qu^2_1u^1_5 ) 
=(u^1_1 u^2_5 - qu^2_1u^1_5 )u^1_au^1_5 + (1-q^2)(u^1_a u^2_5 - qu^2_a u^1_5)u^1_1 u^1_5, 
\label{z2Fz2} 
\end{align} 
where $ a=2,3,4$. 
\end{lem}
\begin{proof}
Equations \eqref{ukluRqulku}-\eqref{Ru2133} follow directly from 
the R-matrix relations \eqref{Rm}. 
Equation \eqref{u1233} is obtained by subtracting $q$ times \eqref{Ru1233}
from \eqref{Ru2133} and then dividing both sides by $-(q^{-1}+1)$. 
Choosing $k=l+1\in\{2,5\}$ in \eqref{uijukl} and acting on both sides by $E_1$ or $F_1$ 
from the left 
shows \eqref{uijukk}. The proven relations imply \eqref{hw2}, \eqref{hw1}, 
\eqref{u14u15} and \eqref{u15u12u15} by elementary computations.  
For instance, \eqref{ukluRqulku},  \eqref{u15u51} and \eqref{uijukk} give  
$$
u^1_1(u^1_1 u^2_5 - qu^2_1u^1_5) = u^1_1(q u^2_5 u^1_1 - q^{-1} u^2_1u^1_5) 
=q(u^1_1 u^2_5 - qu^2_1u^1_5)u^1_1, 
$$
which yields \eqref{hw2}. 

Note that, if we prove one of the relations \eqref{u11Fz2}-\eqref{uuu145} 
and  \eqref{u14u15}-\eqref{z2Fz2} for $a=2$,  then the validity of this 
equation for $a=3$ and $a=4$ follows by acting on 
both sides with $E_2$ from the left. Similarly, if one of these relations  
holds for $a=4$, then acting on both sides with $F_2$ from the left 
yields the result for $a=3$ and $a=2$. 

From \eqref{hw2}, we obtain \eqref{u11Fz2} for $a=4$ by acting on 
both sides from the left by $F_1$.  
Analogously, we get \eqref{uuu515} for $a=2$ from \eqref{hw1} by acting 
on both sides by $E_1$. 
Next, using \eqref{ukluRqulku}-\eqref{uij'kl'},  we compute that 
\begin{align*}
&u^2_5(u^1_1 u^2_4 - qu^2_1u^1_4) 
= q^{-1}\big( u^1_1u^2_5 +(q-q^{-1})u^2_1u^1_5\big)u^2_4 -q^{-1}u^2_1u^2_5 u^1_4\\
&= q^{-1}\big( u^1_1u^2_5 \!+\!(q\!-\!q^{-1})u^2_1u^1_5\big)u^2_4 -q^{-1}u^2_1
\big(u^1_4u^2_5\! +\!(q\!-\!q^{-1})u^2_4u^1_5\big) = q^{-2}(u^1_1 u^2_4 - qu^2_1u^1_4) u^2_5. 
\end{align*} 
This shows \eqref{u212514} for $a=4$.  
The proofs of  \eqref{u11u25}, \eqref{uhw} and \eqref{u25E1z2} are similar.   
Much in the same way,  \eqref{uuu145}  follows from 
\begin{align*}
&(u^1_1u^2_5 - q u^2_1u^1_5)u^1_4 
= u^1_1 u^1_4 u^2_5 + (q^{-1} -q) u^1_1 u^2_4 u^1_5 -  u^2_1 u^1_4 u^1_5 \\
&= u^1_1 (u^1_4 u^2_5\! -\! q u^2_4 u^1_5) 
\! + \! q^{-1}  (u^1_1 u^2_4 \! -\!  q u^2_1 u^1_4) u^1_5  
\hsp = \hsp  q^2 (u^1_4 u^2_5\!-\!q u^2_4 u^1_5)u^1_1
\!+\! q^{-1}  (u^1_1 u^2_4\! -\!  q u^2_1 u^1_4) u^1_5, 
\end{align*}
where we applied 
\eqref{ukluRqulku}  and  \eqref{uij'kl'} in first equality, and \eqref{u11u25} in the third.

Equations \eqref{hw2}-\eqref{uhw}  imply immediately  
\eqref{hw3},  \eqref{u1u5Fz2}  and \eqref{u1u5Ez2}. The commutation relations of $u^2_1$ with 
$u^1_1 u^2_a - qu^2_1u^1_a$ and $u^1_a u^2_5 - qu^2_au^1_5$ 
can easily be deduced from \eqref{ukluRqulku}-\eqref{u15u51} and \eqref{uijukk}, 
and the commutation relations of $u^1_1$, $u^1_5$ and $u^2_5$  with 
$u^1_1 u^2_a - qu^2_1u^1_a$ and $u^1_a u^2_5 - qu^2_au^1_5$ 
are shown in \eqref{u11Fz2}-\eqref{u25E1z2}. Combining these results yields 
\eqref{hw4} and \eqref{z2Ez2}. 
Acting on both sides of \eqref{hw1} from the left by $F_1$  gives 
$$
u^1_4 (u^1_1 u^2_5 - qu^2_1u^1_5)  +   u^1_5(u^1_1 u^2_4 - qu^2_1u^1_4) 
= (u^1_1 u^2_4 - qu^2_1u^1_4) u^1_5  + q^{-1}(u^1_1 u^2_5 - qu^2_1u^1_5) u^1_4. 
$$ 
Therefore, by \eqref{uhw}, 
$$
u^1_4 (u^1_1 u^2_5 - qu^2_1u^1_5) 
= q^{-1}(u^1_1 u^2_5 - qu^2_1u^1_5) u^1_4 - (1-q^{2}) u^1_5 (u^1_1 u^2_4 - qu^2_1u^1_4). 
$$
Multiplying this relation by $u^1_1$ from the left and applying \eqref{hw2} 
yields \eqref{hw5} for $a=4$. 
Equation~\eqref{z2Fz2} can be proven in the same way by acting with $E_1$ on \eqref{hw2}. 
\end{proof}

\section{Differential calculus from a quantum Bernstein--Gelfand--Gelfand resolution} 
\label{secBGG}
Let 
$$
M^{(0,0)}\hsp:=\hsp  \C, \  \,M^{(0,2)}\hsp:=\hsp  \lin\{\hu^{-1}, \hu^0, \hu^1\}, \ \,
M^{(1,2)}\hsp:=\hsp \lin\{\ho^{-1}, \ho^0, \ho^1\}, \ \, M^{(3,0)}\hsp:=\hsp  \C\hs \nu^0
$$ 
be left $\Uql$-modules with the left actions described in \eqref{rep1}  
and \eqref{rep0}. From \cite{HK2} and \cite{HK3}, 
we conclude that there 
exists an exact sequence of left $\Uqg$-modules 
\begin{align} \nonumber
0 \ \lla\ \C \ 
&\overset{\ \,\vare\ot\id}{\lla}\ \Uqg \ot_{\Uql} M^{(0,0)}\ \overset{\ \,\varphi_0}{\lla}\ 
\Uqg \ot_{\Uql} M^{(0,2)}  \\ 
& \overset{\ \,\varphi_1}{\lla}\ \Uqg \ot_{\Uql} M^{(1,2)}\ \overset{\ \,\varphi_2}{\lla}\ 
\Uqg \ot_{\Uql} M^{(3,0)}\ \lla \ 0,                          \label{BGG} 
\end{align} 
and that the $\Uqg$-module maps $\varphi_0$, $\varphi_1$ and $\varphi_2$ 
are unique up to a non-zero scaling factor. 
Our aim is to find explicit expressions for these maps,  
which en passant proves the existence of the exact sequence \eqref{BGG} 
in our setting. 

For shortness of notation, we write $x\ut y$ for $x\hsp\ot_{\Uql}\hsp y$ and  
set  
\[   \label{X} 
X_{-1}:= E_1, \quad X_{0}:= E_2E_1, \quad X_{1}:= E_2^2E_1. 
\] 
These elements are known to generate the (anti-)holomorphic quantum tangent space of 
the invariant algebra $B$, see \cite{HK1}. 
By equivariance and preservation of weights, we deduce that 
$$
\varphi_2(1 \ut \nu^0) = 
\ha_{-1}X_1 \ut \ho^{-1} + \ha_{0}X_0 \ut \ho^{0} + \ha_1 X_{-1}\ut \ho^{1}. 
$$
To determine the constants, we apply $E_2$ to the highest weight vector and obtain 
\begin{align*}
0&= \varphi_2(1 \ut E_2 \nu^0) = E_2 \varphi_2(1 \ut \nu^0) 
= \ha_{-1} E_2^3E_1 \ut \ho^{-1} + \ha_{0}E_2^2E_1 \ut \ho^{0} + \ha_1 E_2E_1\ut \ho^{1}\\ 
&=\ha_{-1}([3]_2E^2E_1\ut E_2 \ho^{-1} - [3]_2EE_1\ut E_2^2 \ho^{-1})
+ \ha_{0}E_2^2E_1 \ut \ho^{0} + \ha_1 E_2E_1\ut \ho^{1}\\
&= (\ha_0+\ha_{-1}[3]_2) E_2^2E_1 \ut \ho^{0} + (\ha_1 -\ha_{-1}[3]_2) E_2E_1\ut \ho^{1}, 
\end{align*} 
where we used the quantum Serre relation \eqref{Sr1} and 
the third equation in \eqref{rep1} 
together with $E_2^3 \ho^{-1} =0$. Setting $\ha_1:=1$, 
we get $\ha_0= -1$ and $\ha_{-1} = [3]_2^{-1}$, thus 
\[
\varphi_2(Z \ut \nu^0) 
= \ebox{\frac{1}{\,[3]_2}} ZX_1 \ut \ho^{-1} - ZX_0 \ut \ho^{0} + Z X_{-1}\ut \ho^{1}, 
\quad Z\in\Uqg.  
\] 

Similarly, to find an expression for $\varphi_1$, we make the ansatz 
$$
\varphi_1(1 \ut \ho^{1}) = \hb_0 X_1 \ut \hu^0 + \hb_1 X_0 \ut \hu^1
$$
and apply $E_2$ to the highest weight vector $\ho^1$ obtaining 
$$
0= \varphi_1(1 \ut E_2 \ho^1) =\hb_0 E_2^3E_1 \ut \hu^0 + \hb_1 E_2^2E_1 \ut \hu^1 
= (\hb_1 + [3]_2\hb_0) E_2^2E_1 \ut \hu^1. 
$$

Setting $\hb_1:=[3]_2$ and $\hb_0:= - 1$ gives 
\[
\varphi_1(Z \ut \ho^{1}) = - Z X_1 \ut \hu^0 + [3]_2\hs Z X_0 \ut \hu^1, \quad Z\in\Uqg. 
\] 
Using the equivariance of $\varphi_1$ and the defining relations of $\Uqg$, we compute 
\begin{align*}
\varphi_1(1 \ut \ho^0) &= \ebox{\frac{1}{\,[2]_2}}\varphi_1(1 \ut F_2 \ho^1) 
= \ebox{\frac{-1}{\,[2]_2}} F_2 E_2^2E_1\ut \hu^0 
+\ebox{\frac{[3]_2}{[2]_2}} F_2 E_2 E_1 \ut \hu^1 \\
&= \ebox{\frac{-1}{[2]_2}}\Big(E_2^2E_1\ut F_2\hu^0 
- E_2E_1 \ut\ebox{\frac{K_2-K_2^{-1}}{q_2-q_2^{-1}}}\hu^0 
- E_2E_1 \ut\ebox{\frac{q^{-1}K_2-qK_2^{-1}}{q_2-q_2^{-1}}}\hu^0 \Big) \\ 
&\quad \ + \ebox{\frac{[3]_2}{[2]_2}}\Big(E_2E_1\ut F_2\hu^1 
- E_1 \ut \ebox{\frac{q^{-1}K_2-qK_2^{-1}}{q_2-q_2^{-1}}}\hu^1\Big) \\
&= - E_2^2 E_1\ut \hu^{-1} + [2]_1E_2E_1 \ut \hu^0, 
\end{align*} \
where we used \eqref{rep1} and $[3]_2-1 = q+q^{-1} = [2]_1$. Hence 
\[
\varphi_1(Z \ut \ho^{0}) = -  Z X_1 \ut \hu^{-1} +[2]_1 Z X_0 \ut \hu^0, \quad Z\in\Uqg. 
\] 
Analogously, 
\begin{align*}
\varphi_1(1 \ut \ho^{-1}) &= \ebox{\frac{1}{\,[2]_2}}\varphi_1(1 \ut F_2 \ho^0) 
= -\ebox{\frac{1}{[2]_2}} F_2E_2^2 E_1\ut \hu^{-1} 
+  \ebox{\frac{[2]_1}{[2]_2}} F_2E_2E_1 \ut \hu^0\\
&= (-1-\ebox{\frac{[4]_2}{[2]_2}} +[2]_1) E_2E_1 \ut \hu^{-1} + [2]_1  E_1 \ut \hu^0. 
\end{align*} 
Since $[2]_1-\ebox{\frac{[4]_2}{[2]_2}} = q+q^{-1} - \frac{q^2-q^{-2}}{q-q^{-1}} =0$, we obtain 
\[
\varphi_1(Z \ut \ho^{-1}) = - Z X_0 \ut \hu^{-1} +[2]_1 Z X_{-1} \ut \hu^0, \quad Z\in\Uqg. 
\]

Finally, we set 
$$
\varphi_0(1\ut \hu^1) := X_1 \ut 1
$$
since $1\ut \hu^1 \in \Uqg \ot_{\Uql} M^{(0,2)}$ and $X_1\ot 1\in  \Uqg \ot_{\Uql} M^{(0,0)}$ 
are both highest weight vectors of irreducible $\Uql$-modules with the same weights. 
In the latter case, this follows from the quantum Serre relation \eqref{Sr1}.  
Furthermore, 
\begin{align*}
E_2(X_0\ut 1) &= X_1 \ut 1 = \varphi_0(1\ut \hu^1) 
= \varphi_0(1\ut E_2 \hu^0) = E_2  \varphi_0(1\ut \hu^0), \\ 
E_2(X_{-1}\ut 1) &= X_0 \ut 1 = \varphi_0(1\ut \hu^0) 
= \varphi_0(1\ut E_2 \hu^{-1}) = E_2  \varphi_0(1\ut \hu^{-1}), 
\end{align*} 
from which we conclude that, for all $Z_{-1}, Z_{0}, Z_{1} \in \Uqg$, 
\[
\varphi_0(Z_{-1}\ut \hu^{-1} + Z_0 \ut\hu^0 + Z_1 \ut \hu^1 ) 
= (Z_{-1}X_{-1} + Z_{0}X_{0}+Z_{1}X_{1}) \ut 1.  
\]
One can now verify directly that $\varphi_{j-1}\circ\varphi_{j}=0$.  
Moreover, $(\vare \ot \id)\circ \varphi_0=0$ which proves that \eqref{BGG} is a complex. 
In fact, by \cite[section 3.4]{HK2}, it is exact. 

As in \cite{HK3}, 
we obtain a differential calculus over $B= \OGq^{\inv(\Uql)}$ by considering the 
locally finite dual of the complex \eqref{BGG}. First recall that the dual $M'$ 
of a left $\cU$-module $M$ 
becomes a right $\cU$-module and the dual $N'$ of a right $\cU$-module $N$ 
becomes a left $\cU$-module
by setting 
\[ \label{da} 
 (f\rtr X)(m):= f(Xm) \quad \text{and} \quad X(g)(n):= g(nX) 
\]
respectively, where $m\in M$, $f\in M'$, $n\in N$, $g\in N'$ and $X\in\cU$. 
For $k=0,\ldots,3$, define 
\[ \label{Ok}
\hO^{(0,k)}:=\big\{ f\in \big(\Uqg \ot_{\Uql} M^{\hl_k}\big)^\prime : 
\dim\big(f\rtr\, \Uqg\big) < \infty\big\}, 
\]
where 
\[     \label{lk}
\hl_0:=(0,0),\quad \hl_1:=(0,2),\quad \hl_2:=(1,2),\quad \hl_3:=(3,0).
\]  
The finiteness condition in \eqref{Ok} is what is referred to as the locally finite part. 
From \cite[Lemma 6.4]{HK3}, it follows that 
$$
\hO^{(0,k)} \subset \OGq \ot {M^{\hl_k}}', \quad k=0,1,2,3, 
$$ 
and 
\[
\hO^{(0,0)}\cong 
\{ b\in \OGq\,:\, X\ltr b = \vare(X)\hs b \ \,\text{for all } X\in\Uql\} = B.  \label{O0} 
\]
Let $\{\nu_0\}$, $\{\ho_{-1}, \ho_0, \ho_1\}$ and $\{\hu_{-1}, \hu_0, \hu_1\}$,  
denote the dual bases of $\{\nu^0\}$, $\{\ho^{-1}, \ho^0, \ho^1\}$ and 
$\{\hu^{-1}, \hu^0, \hu^1\}$, respectively. 
Then an element $\msum{j=-1,0,1}{} a_j\ot \hu_j\in \OGq \ot {M^{\hl_k}}'$ 
belongs to $\hO^{(0,k)}$ if and only if 
$\msum{j=-1,0,1}{} (X\ltr a_j)\ot \hu_j = \msum{j=-1,0,1}{} a_j\ot (\hu_j\rtr X)$. 
From \eqref{rep1} and \eqref{da}, we get 
$\hu_j\rtr K_1 = q^{1-j} \hu^j, \ \hu_j\rtr K_2 = q^j \hu^j, \ \hu_j\rtr E_2 = \hu_{j-1}$ 
and $\hu_j\rtr F_2 = [2]_2 \hu_{j+1}$, 
therefore 
\begin{align} \nonumber 
\hO^{(0,1)} &= \Big\{ \msum{j=-1}{1} a_j\ot \hu_j \in \OGq \ot {M^{(0,2)}}': \\ &\hspace{22pt} 
K_1\ltr a_j\hsp =\hsp q^{1-j} a_j, \ K_2\ltr a_j\hsp =\hsp  q^{j} a_j,\ 
E_2\ltr a_j\hsp =\hsp a_{j+1},\ F_2\ltr a_j\hsp =\hsp [2]_2a_{j-1}\big\}. \label{O1} 
\end{align} 
Analogously, 
\begin{align} \nonumber      
\hO^{(0,2)} &= \Big\{ \msum{j=-1}{1} b_j\ot \ho_j \in \OGq \ot {M^{(1,2)}}': \\ &\hspace{22pt} 
K_1\ltr b_j\hsp =\hsp q^{2-j} b_j, \ K_2\ltr b_j\hsp =\hsp  q^{j} b_j,\ 
E_2\ltr b_j\hsp =\hsp b_{j+1},\ F_2\ltr b_j\hsp =\hsp [2]_2b_{j-1}\big\}  \label{O2} 
\end{align} 
and, under the isomorphism $a\ot \nu_0 \mapsto a$,  
\[      \label{O3} 
\hO^{(0,3)} \cong \{ a\in \OGq :  K_1\ltr a\!=\!q^3a,\ K_2\ltr a\!=\!a, \ 
E_2\ltr a\!=\!F_2\ltr a\!=\!0\}.  
\]

Note that $\hO^{(0,k)}$, $k=1,2,3$, are $B$-bimodules with the left 
and the right $B$-module structure 
given by left and right multiplication on the left tensor factor, respectively. 
Since $K_i\ltr (b_1ab_2)= b_1(K_i\ltr a)b_2$, $E_2\ltr (b_1ab_2)= b_1(E_2\ltr a)b_2$,  
$F_2\ltr (b_1ab_2)= b_1(F_2\ltr a)b_2$ for all $a\in \OGq$, $b_1,b_2\in B$ and $i=1,2$,  
which follows from \eqref{ma} and 
$K_i\ltr b=b$, $E_2\ltr b =0= F_2\ltr b$ for all $b\in B$, 
the $B$-module structure is well defined. 

Of course, we could have described $\hO^{(0,k)}$ more concisely by a cotensor product 
as in \cite{HK3}. However, 
the explicit description using a fixed basis of ${M^{\hl_k}}'$ 
will be useful for the computation of the spectrum of the Dirac operator. 

Next, consider the pull-backs 
$$
\db_k:= \varphi_k^*\ : \ \hO^{(0,k)}\ \lra \ 
\big(\Uqg \ot_{\Uql} M(\hl_{k+1})\big)^\prime,\quad 
\db_k(f)(x)=f(\varphi_k(x)). 
$$
Since we are interested in explicit formulas, we compute for all $Z_{-1},Z_0,Z_{1} \in \Uqg$ 
\begin{align*}
&\db_0(b)(Z_{-1}\ut \hu^{-1} + Z_{0}\ut \hu^{0} + Z_{1}\ut \hu^{1}) 
= b(Z_{-1}X_{-1}+Z_{0}X_{0}+Z_{1}X_{1})   \\
&= (X_{-1}\ltr b \ot \hu_{-1} + X_{0}\ltr b \ot \hu_{0} + X_{1}\ltr b \ot \hu_{1}) 
(Z_{-1}\ut \hu^{-1} + Z_{0}\ut \hu^{0} + Z_{1}\ut \hu^{1}), 
\end{align*}
hence, for all $b\in B$, 
\[ \label{db0}
\db_0(b) = X_{-1}\ltr b \ot \hu_{-1} + X_{0}\ltr b \ot \hu_{0} + X_{1}\ltr b \ot \hu_{1}. 
\]
Similar computations show that 
\begin{align}  \nonumber
\db_1(a_{-1} \ot \hu_{-1} + a_{0} \ot \hu_{0}+ a_{1} \ot \hu_{1}) 
=(-\ebox{\frac{1}{\,[3]_2}} X_0\ltr a_{-1} 
+ \ebox{\frac{[2]_1}{[3]_2}} X_{-1}\ltr a_{0}) \ot \ho_{-1} &\\
+ (-\ebox{\frac{1}{\,[3]_2}} X_1\ltr a_{-1} 
+ \ebox{\frac{[2]_1}{[3]_2}} X_{0}\ltr a_{0}) \ot \ho_{0}
+ (-\ebox{\frac{1}{\,[3]_2}} X_1\ltr a_{0} + X_{0}\ltr a_{1}) \ot \ho_{1}, &      \label{db1}
\end{align}
and, after composing it with the isomorphism $a\ot \nu_0 \mapsto a$, 
\[  \label{db2}
\db_2(b_{-1} \ot \hu_{-1} + b_{0} \ot \hu_{0}+ b_{1} \ot \hu_{1}) = 
\ebox{\frac{1}{\,[3]_2}} X_1\ltr b_{-1} - X_0\ltr b_0 + X_{-1}\ltr b_1. 
\] 
In particular, since the left and the right $\Uqg$-action on $\OGq$ commute, 
it follows from \eqref{Ok} and \eqref{db0}-\eqref{db2} that 
$\db_k$ maps $\hO^{(0,k)}$ into $\hO^{(0,k+1)}$. Therefore we obtain a complex 
\[  \label{dbc}
0 \ \lra \ \hO^{(0,0)}\ \overset{\db_0}{\lra}\ \hO^{(0,1)}\ \overset{\db_1}{\lra}\ \hO^{(0,2)} 
\ \overset{\db_2}{\lra}\ \hO^{(0,3)} \ \lra \ 0. 
\] 
Moreover, the analogue of \cite[Lemma 7.4]{HK3} shows that
\[ 
\hO^{(0,k)}=\lin\{ a_0\db(a_1)\wedge\ldots\wedge\db_k(a_k)\,:\, a_0,\ldots,a_k\in B\},  
\]
where the wedge product is determined by 
\begin{align*}
&a\db_0(b)\wedge \omega_k 
= a\db_k(b\omega_k)- ab\db_k(\omega_k), \quad w_k\in \hO^{(0,k)}, \ \  k=0,1,2, \\
&\big(a_0\db_0(a_1)\wedge\db_0(a_2)\big)\wedge a_3\db(a_4) 
= a_0\db_0(a_1)\wedge \big(\db_0(a_2)\wedge a_3\db(a_4)\big), 
\quad a_0,\ldots, a_4\in B, 
\end{align*} 
and $\omega_k\wedge \omega_j=0$ for $\omega_k\in \hO^{(0,k)}$ 
and $\omega_j\in \hO^{(0,j)}$ with $k+j>3$. 
It follows now from \cite[Section 7.3]{HK3}, that the complex \eqref{dbc} is isomorphic 
to the differential calculus 
$\Gamma_{\db,\mathfrak{u}}^\wedge$ from \cite[Section 3.3.2]{HK1}. 

\section{Branching rules}  \label{sec-bl} 

To distinguish between irreducible highest weight representations of $\Uqg$ and $\Uql$, 
we will use double parentheses to designate the irreducible highest weight 
representations of \hs$\Uqg$. In what follows, 
\hsp$\up{\cc m,j\CC}V$ denotes a subspace of $\OGq$ carrying an irreducible 
$\Uqg$-re\-pre\-sen\-ta\-tions 
of weight $(m,j)$ with respect to the right regular representation $\pi_R$, and 
$V^{\cc m,j\CC} \subset \OGq$ denotes a subspace of $\OGq$ 
carrying an irreducible $\Uqg$-re\-pre\-sen\-ta\-tions 
of weight $(m,j)$ with respect to the left regular representation $\pi_L$, see \eqref{regrep}. 
From the Peter-Weyl theorem for compact quantum groups, it is known that 
\[ \label{PW} 
 \OGq \  \cong 
 \underset{(m,j)\in \N_0\otimes 2\N_0}{\bigoplus} \up{\cc m,j\CC}V \ot V^{\cc m,j\CC}. 
\] 
The subspace of  $\OGq$ isomorphic to 
$\up{\cc m,j\CC}V \ot V^{\cc m,j\CC}$ is spanned by 
the matrix coefficients $t^{(m,j)}_{\hl,\mu}$  
of an irreducible $\OGq$-corepresentation on itself such that the vector space 
generated by the matrix coefficients $t^{(m,j)}_{\hl,\mu}$ 
for fixed $\mu$ is isomorphic to $\up{\cc m,j\CC}V$ and  
the vector space generated by the matrix coefficients $t^{(m,j)}_{\hl,\mu}$ 
for fixed $\hl$ is isomorphic to $V^{\cc m,j\CC}$. 

The descriptions of $\hO^{(0,0)},\ldots, \hO^{(0,3)}$ in Equations \eqref{O0}-\eqref{O3} show 
that there exists a nontrivial element in $\hO^{(0,k)}$ whenever the representation of $\Uql$ 
on $V^{\cc m,j\CC}$ contains an (irreducible) subrepresentation determined by 
the corresponding highest weight 
$\hl_k$ from \eqref{lk}, 
and that all $(0,k)$-forms are given by linear combinations of elements of this form. 
This leads us to the so-called branching rule, the problem of determining if a certain 
irreducible representation of $\Uql$  occurs as a subrepresentation of 
$\Uqg$ on $V^{\cc m,j\CC}$ and 
with which multiplicity. 

Let $V^{\cc m,j\CC}_\mu\subset V^{\cc m,j\CC}$ denote the subspace 
of vectors of weight $\mu=(\mu_1,\mu_2)$, 
i.\,e., $v_\mu \in V^{\cc m,j\CC}_\mu$ if and only if $v_\mu\in V^{\cc m,j\CC}$ and 
$K_i\hs \ltr\hs v_\mu = q^{\mu_i} v_{\mu}$, $i\!=\!1,2$. 
Then the branching rule for $\hl\in\Z\times \N_0$ can be solved by determining the space 
of highest weight vectors in 
$V^{\cc m,j\CC}_{\,\hl}$ for the $\Uql$-representation, that is, the vectors 
$v_{\hl} \in V^{\cc m,j\CC}_{\,\hl}$ such that $E_2\ltr v_{\hl}= 0$. 
In particular, the multiplicities are given by 
$\dim\!\big(\ker(E_2\uhr{V^{\cc m,j\CC}_{\,\hl}})\big)$. 
Equation~\eqref{irreps} shows that $F_2 v_{(\mu_1,\mu_2)} \neq 0$ 
for $v_{(\mu_1,\mu_2)}\in V^{\cc m,j\CC}_{(\mu_1,\mu_2)}\setminus \{0\}$ whenever $\mu_2>0$ 
and that, in this case, $E_2F_2 v_{(\mu_1,\mu_2)} \neq 0$. Therefore 
$E_2: F_2\big(V^{\cc m,j\CC}_{\hl + \alpha_2}\big)\lra V^{\cc m,j\CC}_{\hl + \alpha_2}$ 
is an isomorphism, hence the multiplicities are equal to  
$\dim(V^{\cc m,j\CC}_{\hl})-\dim(V^{\cc m,j\CC}_{\hl + \ha_2})$. 
This difference can easily be computed by Kostant's multiplicity formula (see e.\,g.~\cite{JEH}) 
$$  
\dim (V^{\mu}_{\lambda}) 
= \sum_{\w \in W}(-1)^{\ell(\w)} \# \Pi \big( \w(\mu + \delta ) - (\lambda + \delta) \big),
$$
where $\delta$ denotes the half-sum of positive roots,  
$W$ the Weyl group, $\ell(\w)$ the  length of the Weyl group element $w\in W$,  
and $\Z^2 \ni \nu \mapsto \# \Pi(\nu) $ corresponds to the partition function, i.\,e.\  the 
number of all possible ways of writing  $\nu$ as a sum of positive roots. 
In our situation, $\delta = (1,1)$ and 
$\# \Pi(\nu)$ equals the number of elements of the set 
\[
\Pi(\nu) := 
\{ (n_1,n_2,n_3,n_4)\in \mathbb{N}_0^4 \ | \ 
\nu = n_1 \alpha_1+n_2 \alpha_2 +n_3 (\alpha_1+\alpha_2) + n_4 (\alpha_1+2\alpha_2) \},  
\] 
see \eqref{Rp}. Setting 
$$ 
d(\nu):= \# \Pi(\nu) - \# \Pi(\nu-\alpha_2), \quad \nu \in \Z^2, 
$$
we obtain 
\begin{align}\label{Kf}
\dim(V^{(\mu)}_{\hl})-\dim(V^{(\mu)}_{\hl + \alpha_2}) 
= \sum_{\w \in W}(-1)^{\ell(\w)}  
d \big( \w(\mu + \delta ) - (\lambda + \delta) \big),\ \ k=0,\ldots,3,  
\end{align}
where the parenthesis in the expression $(\mu)$ indicate that 
$V^{(\mu)}$ denotes the vector space of an irreducible $\Uqg$-representation of weight $\mu$. 
The value of $d(\nu)$ is easily derived from combinatorial considerations. 
In the following lemma, we compute those cases that are of interest to us. 

\begin{lem}\label{lemmapartis} 
Let $\lfloor x \rfloor$ denote the integer part of $x\geq 0$. Then 
\begin{enumerate}[label=(\roman*)]
\item $d(k,2l) =1+ \lfloor \efrac{k}{2} \rfloor $ for \,$(k,l)\in \N_0\times 2\N_0$. \label{d1} 
\item $d(k,-2l) = 1+\lfloor \efrac{k}{2}\rfloor -l$ 
for \,\,$(k,l)\in \N_0\times 2\N_0$ and $k\geq 2l$. \label{d2} 
\item $d(k,2l) = 0$ if $(k,2l)\in \Z\times 2\Z$ such that either $k<0$,  
or $k\geq 0$ and $-2l > k$. \label{d3} 
\end{enumerate}  
\end{lem}
\begin{proof}
Obviously, for all $\xi\in \Pi(\nu-\alpha_2)$, we have $\xi +\alpha_2 \in \Pi(\nu)$. 
On the other hand, if $\kappa=(n_1,n_2,n_3,n_4)\in \Pi(\nu)$ such that 
$\kappa -\alpha_2 \in \Pi(\nu -\alpha_2 )$, then $n_2>0$. Therefore,  
\begin{align*}
d(\nu) = \# \Pi(\nu) - \# \Pi(\nu-\alpha_2) = \# \{(n_1,n_2,n_3,n_4)\in \Pi(\nu) \ | \ n_2=0 \}.
\end{align*}
Hence $d(\nu)$ equals the number of solutions of 
$\nu =n_1 \alpha_1 +n_3 (\alpha_1+\alpha_2) + n_4 (\alpha_1+2\alpha_2)$ 
with $n_1,n_2,n_3\in \N_0$. 
For $\nu=(x,2y)\in \Z\times 2\Z$, $\alpha_1= (2,-2)$, $\alpha_1+\alpha_2=(1,0)$ 
and $\alpha_1+2\alpha_2=(0,2)$, we need to determine 
the number of solutions of  
$$   
x=2n_1+n_3, \quad y= -n_1+n_4, \quad n_1,n_2,n_3\in \N_0. 
$$
Clearly, there is no solution if $x<0$ or if $x+2y<0$. 
This implies \ref{d3}. 

If $x\geq 0$, then $n_3= x- 2 n_1$ and $n_4 = y +n_1$, so that $n_3$ and $n_4$ are uniquely 
determined by $n_1$. Thus it now suffices to count the number of possible $n_1\in \N_0$. 
For $x\in \N_0$ and $y\in \N_0$, there is a unique solution if and only if 
$n_1\leq \lfloor \efrac{x}{2} \rfloor$. This proves~\ref{d1}. 
Given $x\in \N_0$ and $y\in \Z\setminus \N_0$,  
there is a unique solution if and only if $|y| \leq  n_1\leq \lfloor \efrac{x}{2} \rfloor$, 
which implies \ref{d2}.
\end{proof}

The multiplicities of irreducible representations of the Levi factor on $V^{\cc m,j\CC}$ 
can now be computed 
by applying Lemma \ref{lemmapartis} to \eqref{Kf}.  

\begin{prop} \label{br} 
For $(n,l)\in\N_0^2$, let $V^{\cc n,l\CC}$ be a vector space carrying an 
irreducible $\Uqg$-re\-pre\-sen\-ta\-tion 
of highest weight $(n,l)$, and let $V^{(n,l)}$ denote the  vector space of an irreducible 
$\Uql$-re\-pre\-sen\-ta\-tion as described in \eqref{irreps}. 
\begin{enumerate} 
\item The trivial representation $V^{(0,0)}$ 
of the Levi factor $\Uql$ occurs in $V^{\cc n,l\CC}$ 
with the following multiplicities:  \label{00} 
\begin{enumerate}[label=(\roman*)] 
\item multiplicity 1 in $V^{\cc 2n,2l\CC}$,     \label{1i} 
\item multiplicity 0 in all other cases.   \label{1ii} 
\end{enumerate}
\item The irreducible representation $V^{(0,2)}$ 
of the Levi factor $\Uql$ occurs in $V^{\cc n,l\CC}$ 
with the following multiplicities:  \label{02} 
\begin{enumerate}[label=(\roman*)]
\item multiplicity 2 in $V^{\cc 2n+2,2l+2\CC}$,    \label{2i} 
\item multiplicity 1 in $V^{\cc 2n+2,0\CC}$,    \label{2ii} 
\item multiplicity 1 in $V^{\cc 0,2l+2\CC}$,    \label{2iii} 
\item multiplicity 1 in $V^{\cc 2n+1,2l+2\CC}$,   \label{2iv} 
\item multiplicity 0 in all other cases.    \label{2v}  
\end{enumerate}
\item The irreducible representation $V^{(1,2)}$ 
of the Levi factor $\Uql$ occurs in $V^{\cc n,l\CC}$ 
with the following multiplicities:   \label{12} 
\begin{enumerate}[label=(\roman*)]
\item multiplicity 1 in $V^{\cc 2n+2,2l+2\CC}$,    \label{3i} 
\item multiplicity 1 in $V^{\cc 2n+3,0\CC}$,   \label{3ii} 
\item multiplicity 1 in $V^{\cc 1,2l+2\CC}$,     \label{3iii} 
\item multiplicity 2 in $V^{\cc 2n+3,2l+2\CC}$,   \label{3iv} 
\item multiplicity 0 in all other cases.     \label{3v}  
\end{enumerate}
\item The irreducible representation $V^{(3,0)}$ of the Levi factor 
$\Uql$ occurs in $V^{\cc n,l\CC}$ 
with the following multiplicities:     \label{30} 
\begin{enumerate}[label=(\roman*)]
\item multiplicity 1 in $V^{\cc 2n+3,2l\CC}$,     \label{4i} 
\item multiplicity 0 in all other cases.            \label{4ii} 
\end{enumerate}
\end{enumerate}
\end{prop} 
\begin{proof}
Clearly, in order that a non-zero multiplicity of $V^{\hl}$ in $V^{\cc m,j\CC}$ 
exists, $(m,j)$ must be higher than $\hl$. 
Moreover, since $\hl_k\in \N_0\times 2\N_0$ for all $\hl_k$ from \eqref{lk}, 
the multiplicity of $V^{\hl_k}$ in $V^{\cc m,j\CC}$ 
is 0 whenever $j\notin 2\N_0$.

Let $P^+$, $P^{++}$ and $\w_{\ha_2} \in W$ be as defined in Section \ref{DN}. 
Note that the sets in Lemma \ref{lemmapartis} \ref{d1} and \ref{d2} 
are contained in $P^+$ and $\w_{\ha_2}(P^+)$, respectively. 
It thus follows from Lemma \ref{lemmapartis} that, for $\nu \in \Z\times 2\Z$,  
$d(\nu)=0$ if $\nu \notin P^+ \cup \w_{\ha_2}(P^+)$. 
As $\delta \in P^{++}$, we have $\mu + \delta \in P^{++}$ for all $\mu \in P^+$ 
and $\hl_k + \delta \in P^{++}$ for $k=0, \ldots, 3$. By  Lemma~\ref{lemmapartis}~\ref{d3}, 
the latter implies that $d\big(\nu - (\hl_k + \delta)\big) = 0$ if $d(\nu)=0$. 
Therefore the formula in \eqref{Kf} reduces to 
$$
\dim(V^{(\mu)}_{\hl_k})-\dim(V^{(\mu)}_{\hl_k + \alpha_2}) 
= d \big(\mu - \lambda_k \big) 
- d \big( \w_{\ha_2}(\mu + \delta ) - (\lambda_k + \delta) \big) .   
$$
As explained before Lemma~\ref{lemmapartis}, this difference 
computes the multiplicity of the irreducible $\Uql$-representation of highest weight $\hl_k$ 
in an irreducible $\Uqg$-re\-pre\-sen\-ta\-tion of 
with highest weight $\mu\in \N_0\times 2\N_0$. 
These numbers are easily computed by applying Equation \eqref{r2} and Lemma~\ref{lemmapartis}. 
We present the results in the following tables. \\[-6pt] 
%

\noindent
\begin{tabular}{ |c|l| } 
 \hline
 \ref{00}\ref{1i} & 
 $d\big((2n,2l)\big) - d\big((2n\!+\!2l\!+\!1,-\!2l\!-\!2)  \big) 
 =n \! -\!(n\!+\!l)\!+\!(l\!+\!1) =1 $  \\ 
  \hline
 \phantom{i}\ref{00}\ref{1ii}  & $d\big((2n\!+\!1,2l)\big) 
 - d\big((2n\!+\!2l\!+\!2,-\!2l\!-\!2)  \big) 
=n \!-\! (n\!+\!l\!+\!1)\!+\!(l\!+\!1) =0 $  \\ 
 \hline
\end{tabular}
%
\mbox{ } \\[8pt] 

\noindent
\begin{tabular}{ |c|l| } 
 \hline
 \ref{02}\ref{2i} & 
 $d\big((2n,2l\!-\!2)\big) - d\big((2n\!+\!2l\!+\!1,-\!2l\!-\!4)  \big) 
 =n \! -\!(n\!+\!l)\!+\!(l\!+\!2) =2 $  \\ 
  \hline
 \ref{02}\ref{2ii} & $d\big((2n,- 2)\big) - d\big((2n\!+\!1,- 4)  \big) = n -1 -n +2= 1 $  \\ 
 \hline
 \ref{02}\ref{2iii} & 
 $d\big((0,2l\!-\!2)\big) - d\big((2l\!+\!1,-2l\!-\!4)  \big) = 1-0   =1$  \\ 
  \hline
 \ref{02}\ref{2iv}  & $d\big((2n\!+\!1,2l\!-\!2)\big) 
 - d\big((2n\!+\!2l\!+\!2,-\!2l\!-\!4)  \big) 
=n \!-\! (n\!+\!l\!+\!1)\!+\!(l\!+\!2) =1 $  \\ 
 \hline
\end{tabular}
%
\mbox{ } \\[8pt] 

\noindent
\noindent
\begin{tabular}{ |c|l| } 
 \hline
 \ref{12}\ref{3i} & $d\big((2n\!-\!1,2l\!-\!2)\big) - d\big((2n\!+\!2l,-\!2l\!-\!4) \big) 
 =n \! -\!1 \! -\!(n\!+\!l)\!+\!(l\!+\!2) =1 $  \\ 
  \hline
 \ref{12}\ref{3ii} & $d\big((2n,- 2)\big) - d\big((2n\!+\!1,- 4)  \big) 
 = n\! -\!1 - n  \!+\! 2= 1 $  \\ 
 \hline
 \ref{12}\ref{3iii} & 
 $d\big((0,2l\!-\!2)\big) - d\big((2l\!+\!1,-2l\!-\!4)  \big) = 1-0   =1$  \\ 
  \hline
 \ref{12}\ref{3iv}  & $d\big((2n,2l\!-\!2)\big) - d\big((2n\!+\!2l\!+\!1,-\!2l\!-\!4)  \big) 
                                =n \!-\! (n\!+\!l)\!+\!(l\!+\!2) =2 $  \\ 
 \hline
\end{tabular} 
%
\mbox{ } \\[8pt] 

\noindent
\begin{tabular}{ |c|l| } 
 \hline
 \ref{30}\ref{4i} & $d\big((2n\!-\!2,2l)\big) - d\big((2n\!+\!2l\!-\!1,-\!2l\!-\!2)  \big) 
 =n\! -\!1 \! -\!(n\!+\!l\! -\!1)\!+\!(l\!+\!1) =1 $  \\ 
  \hline
 \phantom{i}\ref{30}\ref{4ii} & $d\big((2n\!-\!3,2l)\big) 
 - d\big((2n\!+\!2l\!-\!2,-\!2l\!-\!2)  \big) 
 =n\! -\!2 \! -\!(n\!+\!l\! -\!1)\!+\!(l\!+\!1)=0 $  \\ 
 \hline
\end{tabular}
%
\mbox{ } \\[6pt] 

The formula in  \ref{12}\ref{3ii}  assumes $n>1$. 
For $n=1$, we have $(3,- 4)\notin P^+ \cup \w_{\ha_2}(P^+)$ so that 
$d\big((2,- 2)\big) - d\big((3,- 4)  \big) =d\big((2,- 2)\big)=1$. 
This completes the proof.  
\end{proof} 

Our next aim is to provide an explicit description of the elements in $\OGq$ which 
generate the subspaces that satisfy the branching rules of Proposition \ref{br}. 
Since the vector space of an irreducible highest weight representation of 
$\Uqg$ (resp.\ $\Uql$) is generated by acting with elements from $\Uqg$ (resp.\ $\Uql$) 
on a non-zero highest weight vector, 
it suffices to restrict our attention to highest weight vectors. 

By a slight abuse of notation, 
which amounts to identifying elements under the isomorphism in \eqref{PW}, let 
$$
\up{\cc n,l\CC}V^{(m,j)} \,\subset\, \up{\cc n,l\CC}V \ot V^{\cc n,l\CC}  \,\subset\, \OGq 
$$
denote the vector space of all elements in $\OGq$ belonging to 
a representation of highest weight $(n,l)$ 
with respect to the right {\it and} left regular representations of $\Uqg$ on $\OGq$, 
and to a representation of highest weight $(m,j)$ 
with respect to the left $\Uql$-action on $\OGq$. 
Note that we do not assume that the representation of  $\Uql$ on $\up{\cc n,l\CC}V^{(m,j)}$ 
is irreducible. 
Furthermore, let   
\begin{align} 
\up{\cc n,l\CC}V^{(m,j)}_{(s_1,s_2)} &:= \{ v\in \up{\cc n,l\CC}V^{(m,j)} : 
 K_i\ltr v= q_i^{s_i} v, \ \ i=1,2 \},  \label{Vlhw} \\
{}^{\ \cc n,l\CC}_{(r_1,r_2)}V^{(m,j)}_{(s_1,s_2)} &:= \{ v\in \up{\cc n,l\CC}V^{(m,j)} : 
\pi_R(K_i)(v) = q_i^{r_i} v, \ \ K_i\ltr v= q_i^{s_i} v, \ \ i=1,2 \}  \label{Vhw}
\end{align}
denote the vector spaces of weight vectors of weight $(s_1,s_2)$ with respect 
to the left $\Uql$-action on $\OGq$  and, in the second case, also 
of weight $(r_1,r_2)$ with respect to the right regular representation of $\Uqg$. 
In particular, ${}^{\cc n,l\CC}_{\,\,(n,l)}V^{(m,j)}_{(m,j)}$ 
is the vector space of highest weights with respect to both representations. 
These highest weight vectors determine the whole representation space $\up{\cc n,l\CC}V^{(m,j)}$
since 
$$
\up{\cc n,l\CC}V^{(m,j)}
= \pi_R(\Uqg)\big(\Uql\ltr {}^{\cc n,l\CC}_{\,\,(n,l)}V^{(m,j)}_{(m,j)}\big) 
=\Uql\ltr\big(\pi_R(\Uqg){}^{\cc n,l\CC}_{\,\,(n,l)}V^{(m,j)}_{(m,j)}\big).
$$ 
A complete list of highest weight vectors in the cases of interest 
will be given in the next proposition. 
\begin{prop} \label{hw} 
Let $u^i_j$, $i,j=1,\ldots,5$, denote the generators of $\OGq$. 
Set 
\[ \label{z1z2}
 z_1:= u^1_1u^1_5, \qquad z_2:= u^1_1 u^2_5 - qu^2_1u^1_5. 
\] 
Then, for $n,l\in\N_0$, 
\begin{enumerate}
\item 
\begin{enumerate}[label=(\roman*)]
    \item ${}^{\cc 2n,2l\CC}_{\,\,(2n,2l)}V^{(0,0)}_{(0,0)}= \lin\{z_2^l z_1^n\}.$ 
\end{enumerate} 
\item  \label{hwV2} 
\begin{enumerate}[label=(\roman*)]
    \item \label{hwV2i} 
    ${}^{\cc 2n+2,2l+2\CC}_{\,\,(2n+2,2l+2)}V^{(0,2)}_{(0,2)}
    =\lin\{ z_2^{l+1}u^1_4u^1_5z_1^{n}, \ \,z_2^{l}(u^1_4 u^2_5 - qu^2_4u^1_5)z_1^{n+1} \}$, 
    \item
  ${}^{\cc 2n+2,0\CC}_{\,\,(2n+2,0)}V^{(0,2)}_{(0,2)}\,=\,\lin\{u^1_4u^1_5z_1^{n}\}$, 
    \item
  ${}^{\cc 0,2l+2\CC}_{\,\,(0,2l+2)}V^{(0,2)}_{(0,2)} \,\,=\, \,
  \lin\{z_2^{l}(u^1_4 u^2_5 - qu^2_4u^1_5)\}$,
    \item 
${}^{\cc 2n+1,2l+2\CC}_{\,\,(2n+1,2l+2)}V^{(0,2)}_{(0,2)}
=\lin\Big\{z_2^{l}(u^1_3 u^2_4 - qu^2_3u^1_4)u^1_5z_1^{n}\Big\}$. 
\end{enumerate}
\item  \label{hwV3}  \begin{enumerate}[label=(\roman*)]
   \item
${}^{\cc 2n+2,2l+2\CC}_{\,\,(2n+2,2l+2)}V^{(1,2)}_{(1,2)}
=\lin\Big\{z_2^{l}(u^1_3 u^2_4 - qu^2_3u^1_4)(u^1_5)^2z_1^{n}  \Big\}$, 
    \item
${}^{\cc 1,2l+2\CC}_{\,\,(1,2l+2)}V^{(1,2)}_{(1,2)}
=\lin\{z_2^{l} (u^1_4 u^2_5 - qu^2_4u^1_5)u^1_5 \}$, 
    \item
${}^{\cc 2n+3,0\CC}_{\,\,(2n+3,0)}V^{(1,2)}_{(1,2)}
=\lin\{u^1_4(u^1_5)^2 z_1^{n}\}$, 
    \item   \label{hwV3iv}
${}^{\cc 2n+3,2l+2\CC}_{\,\,(2n+3,2l+2)}V^{(1,2)}_{(1,2)}
=\lin\{ z_2^{l+1} u^1_4(u^1_5)^2 z_1^{n}, \ \,
z_2^{l}(u^1_4 u^2_5 - qu^2_4u^1_5)u^1_5z_1^{n+1}\}$.  
\end{enumerate}
\item 
\begin{enumerate}[label=(\roman*)]
    \item
${}^{\cc 2n+3,2l\CC}_{\,\,(2n+3,2l)}V^{(3,0)}_{(3,0)}
=\lin\{   z_2^l (u^1_5)^3 z_1^{n} \}$. 
\end{enumerate}
\end{enumerate}
\end{prop} 
\begin{proof} 
First note that the numbers of generating vectors in Proposition \ref{hw} 
agree with the multiplicities 
in Proposition \ref{br}. Therefore it suffices to show that 
the sets of generating vectors are highest 
weight vectors of the corresponding weights and are linearly independent. 

Using the formulas given in \eqref{lac} and \eqref{rac}, elementary computations show that 
the following elements are highest weight vectors for the right regular representation 
of $\Uqg$ and 
the left action of $\Uql$: 
\begin{align*} 
z_1=u^1_1u^1_5, \quad z_2=u^1_1 u^2_5 - qu^2_1u^1_5, \quad u^1_4, \quad u^1_5, \quad 
u^1_4 u^2_5 - qu^2_4u^1_5, \quad u^1_3 u^2_4 - qu^2_3u^1_4. 
\end{align*}  
Since, by \eqref{DE} and \eqref{ma}, the product of highest weight vectors is 
a highest weight vector, 
where the resulting weight is the sum of the individual weights, it follows that all 
elements listed in Proposition \ref{hw} are highest vectors belonging to 
the corresponding weight spaces. 

We need to prove that none of the highest weight vectors is zero. 
Since it is known that $\OGq$ is a domain, it suffices to verify this for  
$u^1_1 u^2_5 - qu^2_1u^1_5$, $u^1_4 u^2_5 - qu^2_4u^1_5$ and  
$u^1_3 u^2_4 - qu^2_3u^1_4$.

Elementary considerations show that $u^1_4 u^2_5 - qu^2_4u^1_5$ is the unique 
(up to a scaling factor) 
highest weight vector of the representation $\up{\cc 0,2\CC}V \ot V^{\cc 0,2\CC}$ from the 
Peter-Weyl decomposition \eqref{PW}, so it cannot be zero. 
Since $E_2^2E_1(u^1_1 u^2_5 - qu^2_1u^1_5)$ yields 
a non-zero multiple of $u^1_4 u^2_5 - qu^2_4u^1_5$, we have necessarily 
$z_2=u^1_1 u^2_5 - qu^2_1u^1_5\neq 0$.   Analogously, 
$E_2E_1(u^1_3 u^2_4 - qu^2_3u^1_4)$ is a non-zero multiple of $u^1_4 u^2_5 - qu^2_4u^1_5$, 
hence $u^1_3 u^2_4 - qu^2_3u^1_4\neq 0$.

It remains to prove that the two generating vectors are linearly independent 
when the multiplicity is 
equal to 2. Let $\alpha, \beta \in \C$ such that 
\begin{align}  \nonumber 
0&= \alpha z_2^{l+1}u^1_4u^1_5z_1^{n} + \beta z_2^{l}(u^1_4 u^2_5 - qu^2_4u^1_5)z_1^{n+1} \\
&= z^l\big(\alpha  (u^1_1 u^2_5 - qu^2_1u^1_5) u^1_4 
+\beta (u^1_4 u^2_5 - qu^2_4u^1_5)u^1_1\big)u^1_5 z^n. 
\label{li} 
\end{align} 
Then 
$\alpha  (u^1_1 u^2_5 - qu^2_1u^1_5) u^1_4  + \beta (u^1_4 u^2_5 - qu^2_4u^1_5)u^1_1 =0$ 
since 
$\OGq$ is a domain. Acting on this expression with $F_1$ from the left yields 
\begin{align*} 
0= - \alpha q (u^1_1 u^2_4 - qu^2_1u^1_4) u^1_4  - \beta q(u^1_4 u^2_4 - qu^2_4u^1_4)u^1_1 
= - \alpha q (u^1_1 u^2_4 - qu^2_1u^1_4) u^1_4 
\end{align*} 
as $u^1_4 u^2_4 - qu^2_4u^1_4 =0$ by \eqref{uijukk}.  
On the other hand, the unique element (up to a constant) in 
$\up{\cc 0,2\CC}V \ot V^{\cc 0,2\CC}$, which 
is a highest weight vector with respect to right regular representation and a 
lowest weight vector with respect to left regular representation, is given by  
$u^1_1 u^2_2 - qu^2_1u^1_2 = -q_2 [2]_2^{-1} F_2^2(u^1_1 u^2_4 - qu^2_1u^1_4)$. 
Therefore $u^1_1 u^2_4 - qu^2_1u^1_4 \neq 0$, hence $\alpha = 0$ and consequently 
also $\beta=0$.  
This proves that the two vectors on the right-hand side of \ref{hwV2}.\ref{hwV2i}
are linearly independent. 

In Case \ref{hwV3}.\ref{hwV3iv}, if 
$\alpha z_2^{l+1} u^1_4(u^1_5)^2 z_1^{n} + \beta z_2^{l}(u^1_4 u^2_5 
- qu^2_4u^1_5)u^1_5z_1^{n+1}=0$, 
then 
$$
0= z^l\big(\alpha(u^1_1 u^2_5 - qu^2_1u^1_5) u^1_4 
+q^{-2}\beta (u^1_4 u^2_5 - qu^2_4u^1_5)u^1_1\big)(u^1_5)^2 z^n, 
$$ 
and therefore $\alpha=\beta=0$ by the same reasoning, see \eqref{li}. 
\end{proof}

\section{Dobeault--Dirac operator}    \label{secD}

The aim of this section is to study the Dobeault--Dirac operator 
$\db + \db^\dagger$ 
on the Hilbert space closure of 
$\hO^{(0,0)}\oplus \hO^{(0,1)}\oplus\hO^{(0,2)} \oplus \hO^{(0,3)}$, where 
$\db = \db_0\oplus \db_1\oplus \db_2$ maps the first three components into the last three, and 
$\db^\dagger$ denotes the Hilbert space adjoint. 

The first step is to assign a ``natural'' inner product to the space of $(0,k)$-forms. 
By \eqref{O0} and \eqref{O3}, we may view $\hO^{(0,0)}$ and 
$\hO^{(0,3)}$ as subspace of $\OGq$. 
Therefore it is natural to define 
\begin{align}  \label{hc0}
&\ip{\cdot}{\cdot}_0 : \hO^{(0,0)} \times \hO^{(0,0)} \lra \C, \quad 
\ip{b_1}{b_2}_0 := c_0\hs h(b_1^*b_2), \\ 
&\ip{\cdot}{\cdot}_3 : \hO^{(0,3)} \times \hO^{(0,3)} \lra \C, \quad 
\ip{a_1}{a_2}_3 :=c_3\hs h(a_1^*a_2),  
\label{hc3}
\end{align}
where $h$ denotes the Haar state from \eqref{Haar} and $c_0,c_3\in (0,\infty)$. 
Furthermore, by \eqref{O1} and \eqref{O2}, 
we have $\hO^{(0,k)} \subset\OGq\ot {M^{\hl_k}}'\cong \OGq\ot \C^3$, $k=1,2$. 
Taking into account that $K_1$ and $K_2$ should act as self-adjoint operators, we require that 
$\lin\{\hu_{-1}, \hu_0, \hu_1\}\subset {M^{(0,2)}}'$ and 
$\lin\{\ho_{-1}, \ho_0, \ho_1\}\subset {M^{(1,2)}}'$ are orthogonal bases. 
For this reason, we first set 
\[  \label{ipc} 
 \ip{\hu_{i}}{\hu_{j}}_{\hl_1} 
 := c_{1j}\hs \delta_{ij},\qquad \ip{\ho_{i}}{\ho_{j}}_{\hl_2} := c_{2j}\hs \delta_{ij}, 
\]
where $\delta_{ij}$ denotes the Kronecker delta and $c_{1j}, c_{2j} \in(0,\infty)$, 
and then consider the tensor product inner product on $\OGq\ot {M^{\hl_k}}'$, $k=1,2$,   
given by 
\begin{align} \label{ip1}
&\Big\langle\hs \msum{i=-1}{1}a_i\ot \hu_i\,, \msum{j=-1}{1}b_j\ot \hu_j\Big\rangle_{\hsp 1} 
= \msum{j=-1}{1} c_{1j}\hs h(a_j^*b_j), \\
&\Big\langle\hs \msum{i=-1}{1}x_i\ot \ho_i\,, \msum{j=-1}{1}y_j\ot \ho_j\Big\rangle_{\hsp 2} 
= \msum{j=-1}{1} c_{2j}\hs h(x_j^*y_j). \label{ip2}
\end{align} 
The inner product on $(0,1)$- and $(0,2)$-forms will be given by the restriction  
of $\ip{\cdot}{\cdot}_1$ and $\ip{\cdot}{\cdot}_2$ 
to $\hO^{(0,1)}\subset \OGq\ot {M^{(0,2)}}'$ and $\hO^{(0,2)}\subset \OGq\ot {M^{(1,2)}}'$, 
respectively. 

Of course, the spectrum of the Dirac operator may depend on the positive real numbers 
$c_0, c_{1j}, c_{2j}$ and $c_3$. A similar dependence on a scaling factor was observed 
in the definition of the differentials $\db_k$. 
The dependence of the spectrum of the Dirac operator 
on these parameters will be discussed in Remark \ref{scale}. 

By the definition of $\hO^{(0,k)}$ as a subspace of $\OGq\ot {M^{\hl_k}}'$, 
the right regular representation 
of \,$\Uqg$ on the first tensor factor of  $\OGq\ot {M^{\hl_k}}'$ restricts to $\hO^{(0,k)}$ 
since it acts only on the left tensor factor in the Peter-Weyl decomposition \eqref{PW}. 
The next lemma shows that the differential complex \eqref{dbc} is equivariant 
with respect to the right regular representation. 
\begin{lem} \label{dbeq}
On \hs$\hO^{(0,k)}\hsp\subset \hsp\OGq\ot {M^{\hl_k}}'$, 
consider the right regular representation 
of \hs$\Uqg$ given by $\pi_R(X)(a\ot v):= a\rtr S^{-1}(X) \ot v$ 
for $a\ot v \in \OGq\ot {M^{\hl_k}}'$ and $X\in \Uqg$. Then $\pi_R$ defines a *-representation 
of $\Uqg$ on $\hO^{(0,k)}$ such that 
$\pi_R(X) \db_k (\ho_k) =  \db_k \big(\pi_R(X) (\ho_k)\big)$ 
for all $\ho_k \in \hO^{(0,k)}$. 

On the subalgebra  $B= \OGq^{\inv(\Uql)}$ of $\OGq$, consider the left action 
$\btr :  \Uqg \times B \to B$, \,$X\btr b:= b \rtr S^{-1}(X)$. 
Then 
$$
\pi_R(X)(b\hs\ho_k) = (X_{(2)}\btr b)\pi_R(X_{(1)})(\ho_k) 
$$ 
for all $X\in  \Uqg$, $b\in B$ and $\ho_k \in \hO^{(0,k)}$. 
\end{lem} 
\begin{proof}
Since, for all $a\in \OGq$ and $X,Y\in \Uqg$, 
\[  \label{pidb} 
Y\ltr\big(\pi_R(X)(a)\big) = Y\ltr a \rtr S^{-1}(X)=\pi_R(X)( Y\ltr a), 
\]
it follows from \eqref{O0}-\eqref{O3} that $\pi_R(X) : \hO^{(0,k)} \to \hO^{(0,k)}$ 
is well defined. As mentioned before \eqref{srep}, 
$\pi_R$ is a *-representation with respect to the inner product on $\OGq$ 
defined by the Haar state.  As a consequence,  by the 
definitions in \eqref{hc0}-\eqref{ip2}, $\pi_R$ yields a *-representation on $\hO^{(0,k)}$. 
Note that the definition of $\db_k$ involves only the left $\Uqg$-action. 
Therefore \eqref{pidb} implies that $\pi_R(X)$ and $\db_k$ commute for all $X\in \Uqg$ and 
$k=0,1,2$. 

As $B= \hO^{(0,0)}$ and $\pi_R(X) : \hO^{(0,0)} \to \hO^{(0,0)}$, we know already that 
the action $\btr$ leaves $B$ invariant. Furthermore, 
for all $\sum_i a_i\ot v_i \in \hO^{(0,k)}\subset \OGq\ot {M^{\hl_k}}'$ 
and all $X\in \Uqg$, an application of \eqref{ma} gives 
\begin{align*}
\pi_R(X)(\msum{i}{}\hs ba_i\ot v_i) &=\msum{i}{} (ba_i)\rtr S^{-1}(X) \ot v_i  
=\msum{i}{} (b\rtr S^{-1}(X_{(2)})(a_i\rtr S^{-1}(X_{(1)})) \ot v_i \\
&= (X_{(2)} \btr b) \pi_R(X_{(1)})(\msum{i}{}\hs a_i\ot v_i), 
\end{align*}
which proves the second part of the lemma. 
\end{proof}

As in Section \ref{sec-bl}, let ${}^{\cc n,l\CC}\hO^{(0,k)}$ denote the space of vectors 
belonging to an irreducible representation of highest weight $(n,l)$ with respect to the 
right regular representation and let 
${}^{\cc n,l\CC}_{\,\,(r,s)}\hO^{(0,k)} \subset {}^{\cc n,l\CC}\hO^{(0,k)}$ 
denote the subspace of weight vectors of weight $(r,s)$. 
By the Peter-Weyl decomposition \eqref{PW},  
elements of $\OGq$ belonging to irreducible representations of different highest weights  
are orthogonal with respect to the inner product given by the Haar state. 
Therefore, from Equations \eqref{O0}-\eqref{O3}, \eqref{Vhw} and Proposition \ref{hw}, 
we obtain 
the orthogonal decompositions
\[ \label{orthoO} 
\hO^{(0,k)} = \bigoplus \big\{ {}^{\cc n,l\CC}\hO^{(0,k)} : n\in\N_0,\ l\in 2\N_0, \ 
\dim\big({}^{\cc n,l\CC}_{\,\,(n,l)}V^{\hl_k}_{\hl_k}\big)\neq 0\big\}, 
\]
where
\[\label{hO}
{}^{\cc n,l\CC}\hO^{(0,k)}  = \pi_R(\Uqg)\big({}^{\cc n,l\CC}_{\,\,(n,l)}\hO^{(0,k)}\big) 
\] 
and, for $n,m\in\N_0$ and $l\in 2\N_0$, 
\begin{align}
&{}^{\cc 2n,2m\CC}_{\,\,(2n,2m)}\hO^{(0,0)} 
= {}^{\cc 2n,2m\CC}_{\,\,(2n,2m)}V^{(0,0)}_{(0,0)}, \quad 
{}^{\cc 2n+3,2m\CC}_{\,\,(2n+3,2m)}\hO^{(0,3)} 
= {}^{\cc 2n+3,2m\CC}_{\,\,(2n+3,2m)}V^{(3,0)}_{(3,0)}\ot\nu^0, \label{iso0} \\
&{}^{\cc n,l\CC}_{\,\,(n,l)}\hO^{(0,1)}
= \big\{ \ebox{\frac{1}{[2]_2^2}} F_2^2\ltr a_1\ot \hu_{-1}  \label{iso1} 
+ \ebox{\frac{1}{[2]_2}} F_2\ltr a_1 \ot \hu_{0}+ a_{1} \ot \hu_{1} \,:\, 
a_1\in {}^{\cc n,l\CC}_{\,\,(n,l)}V^{(0,2)}_{(0,2)} \big\} \\   \label{iso2} 
&{}^{\cc n,l\CC}_{\,\,(n,l)}\hO^{(0,2)}
= \big\{ \ebox{\frac{1}{[2]_2^2}} F_2^2\ltr b_1\ot \ho_{-1} 
+ \ebox{\frac{1}{[2]_2}} F_2\ltr b_1 \ot \ho_{0}+ b_{1} \ot \ho_{1} \,:\, 
b_1\in {}^{\cc n,l\CC}_{\,\,(n,l)}V^{(1,2)}_{(1,2)} \big\} 
\end{align} 
Note that ${}^{\cc n,l\CC}_{\,\,(n,l)}\hO^{(0,k)}$ is uniquely determined by 
${}^{\cc n,l\CC}_{\,\,(n,l)}V^{\hl_k}_{\hl_k}$. 
The following lemma establishes an explicit unitary isomorphisms between these spaces.  

\begin{lem}  \label{equiv} 
On the linear subspaces ${}^{(\mu)}_{\,\,\,\mu}V^{\hl_k}_{\hl_k}\subset \OGq$
from Proposition~\ref{hw}, 
consider the inner product \eqref{iph} given by the Haar state on $\OGq$. 
Let $c_0$ and $c_3$ denote the constants from \eqref{hc0} and \eqref{hc3}, respectively.  
With the constants $c_{kj}$ from \eqref{ipc}, set
\[  \label{ck} 
c_k:= \sqrt{ \ebox{\frac{c_{k,-1}}{q[2]_2^{2}}}  + \ebox{\frac{c_{k,0}}{q[2]_2}}   
+ c_{k,1}}, \qquad k=1,2, 
\] 
Then the linear operators 
\begin{align*}
&J_0:{}^{(\mu)}_{\,\,\,\mu}V^{(0,0)}_{(0,0)} \lra{}^{(\mu)}_{\,\,\,\mu}\hO^{(0,0)}, \quad 
J_0(a):= \ebox{\frac{1}{c_0}}\hs a, \\
&J_1 :{}^{(\mu)}_{\,\,\,\mu}V^{(0,2)}_{(0,2)} \lra{}^{(\mu)}_{\,\,\,\mu}\Omega^{(0,1)}, \quad 
J_1(a):= \ebox{\frac{1}{c_1}}\Big(\ebox{\frac{1}{[2]_2^2}} F_2^2\ltr a\ot \hu_{-1} 
+ \ebox{\frac{1}{[2]_2}} F_2\ltr a \ot \hu_{0}+ a \ot \hu_{1}\Big), \\
&J_2 :{}^{(\mu)}_{\,\,\,\mu}V^{(1,2)}_{(1,2)} \lra{}^{(\mu)}_{\,\,\,\mu}\Omega^{(0,2)}, \quad 
J_2(a):= \ebox{\frac{1}{c_2}}\Big(\ebox{\frac{1}{[2]_2^2}} F_2^2\ltr a\ot \ho_{-1} 
+ \ebox{\frac{1}{[2]_2}} F_2\ltr a \ot \ho_{0}+ a \ot \ho_{1}\Big), \\
&J_3:{}^{(\mu)}_{\,\,\,\mu}V^{(3,0)}_{(3,0)} \lra{}^{(\mu)}_{\,\,\,\mu}\hO^{(0,3)}, \quad 
J_3(a):= \ebox{\frac{1}{c_3}}\hs a \ot \nu^0. 
\end{align*} 
are unitary isomorphisms. 
\end{lem} 
\begin{proof}
We prove the lemma for $J_1$. For $J_2$, the proof is similar and for $J_0$ and $J_3$, 
the proof is even more elementary. From \eqref{iso1}, 
it follows that $J_1$ defines a linear isomorphism. 
Note that, by \eqref{rep1}, 
$\ebox{\frac{1}{[2]_2^2}} E_2^2F_2^2\ltr a = \ebox{\frac{1}{[2]_2}} E_2F_2\ltr a = a$ 
for all $a\in{}^{(\mu)}_{\,\,\,\mu}V^{(0,2)}_{(0,2)}$. 
Let $a_1,a_2\in{}^{(\mu)}_{\,\,\,\mu}V^{(0,2)}_{(0,2)}$.  
Then 
\begin{align*}
&\ip{J_1(a_1)}{J_1(a_2)}_1 \\
&\quad = \ebox{\frac{1}{c_1^2}}\Big( 
\ebox{\frac{c_{1,-1}}{[2]_2^{4}}} h\big((F_2^2\ltr a_1)^*(F_2^2\ltr a_2)\big)  + 
\ebox{\frac{c_{1,0}}{[2]_2^2}} h\big((F_2\ltr a_1)^*(F_2\ltr a_2)\big) 
+ c_{1,1} h(a_1^*a_2)\Big) \\
&\quad = \ebox{\frac{1}{c_1^2}}\Big( 
\ebox{\frac{c_{1,-1}}{[2]_2^{4}}} h\big(a_1^*((K_2^{-1}E_2)^2F_2^2\ltr a_2)\big)  + 
\ebox{\frac{c_{1,0}}{[2]_2^2}} h\big(a_1^*(K_2^{-1}E_2F_2\ltr a_2)\big) 
+ c_{1,1} h(a_1^*a_2)\Big) \\
&\quad = \ebox{\frac{1}{c_1^2}}
\big( \ebox{\frac{c_{1,-1}}{q[2]_2^{2}}}  + \ebox{\frac{c_{1,0}}{q[2]_2}}   
+ c_{1,1}\big)h(a_1^*a_2) = h(a_1^*a_2), 
\end{align*}
where we used the defintion of $J_1$ and \eqref{ip1} in the first equality, 
the fact that $\pi_L$ from \eqref{regrep} is a *-representation in the second equality, 
the defining relations of $\Uqg$ and \eqref{rep1} in the third equality, and \eqref{ck} 
in the las equality. 
Hence the linear isomorphism $J_1$ is also isometric which proves its unitarity. 
\end{proof} 

Lemma \ref{equiv} and Equations \eqref{orthoO}-\eqref{iso2} show that the operator $\db_k$ 
is uniquely determined by the unitarily equivalent action 
on the spaces of highest weight vectors 
${}^{(\mu)}_{\,\,\,\mu}V^{\hl_k}_{\hl_k}$. The next corollary gives explicit formulas.  
\begin{cor} \label{coreth}
The restrictions of \,$\db_k$ to the spaces of highest weight vectors, 
\[  \label{muO} 
\db_k \,:\,{}^{(\mu)}_{\,\,\,\mu}\hO^{(0,k)} \,\lra\,{}^{(\mu)}_{\,\,\,\mu}\hO^{(0,k+1)},
\quad \mu\in\N_0\times 2\N_0,\quad  k=0,1,2, 
\] 
are unitarily equivalent to 
\begin{align}   \label{eth0} 
&\eth_0:= J_1^*\hsp\circ\hsp \db_0\hsp\circ\hsp J_0 \,:\,  
{}^{(\mu)}_{\,\,\,\mu}V^{(0,0)}_{(0,0)} \lra{}^{(\mu)}_{\,\,\,\mu}V^{(0,2)}_{(0,2)}, \quad 
\eth_0(v)=\ebox{\frac{c_1}{c_0}}\hs  X_1(v), \\ 
& \eth_1:= J_2^*\hsp\circ\hsp \db_1\hsp\circ\hsp J_1 \,:\,  
{}^{(\mu)}_{\,\,\,\mu}V^{(0,2)}_{(0,2)} \lra{}^{(\mu)}_{\,\,\,\mu}V^{(1,2)}_{(1,2)}, \quad 
\eth_1(v)= \ebox{\frac{c_2}{c_1}}\hs 
\ebox{\frac{[3]_2-1}{[3]_2}}\hs \pr_{(1,2)}\hsp\circ\hsp X_0(v), \\ 
& \eth_2:= J_3^*\hsp\circ\hsp \db_2\hsp\circ\hsp J_2 \,:\,  
{}^{(\mu)}_{\,\,\,\mu}V^{(1,2)}_{(1,2)} \lra{}^{(\mu)}_{\,\,\,\mu}V^{(3,0)}_{(3,0)}, \quad 
\eth_2(v)=  \ebox{\frac{c_3}{c_2}}\hs \pr_{(3,0)}\hsp\circ\hsp X_{-1}(v),
\end{align} 
where $\pr_{(1,2)}$ and $\pr_{(3,0)}$ are defined in \eqref{pr2} and \eqref{Pr0}, 
respectively. 
\end{cor} 
\begin{proof}
By the equivariance property shown in Lemma \ref{dbeq}, 
the operator $\db_k$ maps highest weight vectors of weight $\mu\in\N_0\times 2\N_0$ 
into highest weight vectors of the same weight, thus \eqref{muO} holds. 
From Lemma \ref{equiv} and \eqref{db0}, it follows that, for all 
$v\in{}^{(\mu)}_{\,\,\,\mu}V^{(0,0)}_{(0,0)}$, 
$$
\eth_0(v) =\ebox{\frac{1}{c_0}}\hs 
J_1^{-1}( X_{-1}\ltr v \ot \hu_{-1} + X_{0}\ltr v \ot \hu_{0} + X_{1}\ltr v \ot \hu_{1}  ) 
= \ebox{\frac{c_1}{c_0}}\hs X_{1}\ltr v.
$$
Similarly, Lemma \ref{equiv} and \eqref{db1} give 
\begin{align*}
\eth_1(v) &= \ebox{\frac{c_2}{c_1}}\hs\big(
\ebox{\frac{-1}{[2]_2[3]_2}} X_1F_2\ltr v + X_{0}\ltr v\big)
= \ebox{\frac{c_2}{c_1}}\hs\big(
\ebox{\frac{-1}{[2]_2 [3]_2}} E_2^2E_1F_2\ltr v + E_2E_1\ltr v\big)\\
&= \ebox{\frac{c_2}{c_1}}\hs\Big(\ebox{\frac{-1}{[2]_2 [3]_2}} \big(
F_2E_2^2E_1\ltr v + E_2\ebox{\frac{K_2-K_2^{-1}}{q_2-q_2^{-1}}}E_1\ltr v 
+\ebox{\frac{K_2-K_2^{-1}}{q_2-q_2^{-1}}}E_2E_1\ltr v \big)
+ E_2E_1\ltr v\Big)\\
&= \ebox{\frac{c_2}{c_1}}\big(
1-\ebox{\frac{1}{[3]_2}}-\ebox{\frac{1}{[2]_2 [3]_2}} F_2E_2 \big) E_2E_1\ltr v 
= \ebox{\frac{c_2}{c_1}}\ebox{\frac{[3]_2-1}{[3]_2}}
\big(1 - \ebox{\frac{1}{[4]_2}} F_2E_2 \big)E_2E_1\ltr v \\
&= \ebox{\frac{c_2}{c_1}}\ebox{\frac{[3]_2-1}{[3]_2}}\pr_{(1,2)}\hsp\circ\hsp X_0(v), 
\end{align*}
where we applied the defining relations of $\Uqg$, $K_2\ltr v = q_2^2 v$,
$[2]_2([3]_2-1) = [4]_2$ and \eqref{pr2}. 
Before using \eqref{pr2}, 
it is important to note that 
$F_2^3 E_1\ltr v\hsp =\hsp E_1F_2^3\ltr v\hsp =\hsp 0$ for all 
$v\hsp \in\hsp{}^{(\mu)}_{\,\,\,\mu}V^{(0,2)}_{(0,2)}$, 
so that $E_1\ltr v \in{}^{(\mu)}_{\,\,\,\mu}V^{(2,0)}_{(2,0)} 
\oplus{}^{(\mu)}_{\,\,\,\mu}V^{(1,2)}_{(2,0)} 
\oplus{}^{(\mu)}_{\,\,\,\mu}V^{(0,4)}_{(2,0)}$\hsp, 
thus $E_2E_1\ltr v \in{}^{(\mu)}_{\,\,\,\mu}V^{(1,2)}_{(1,2)} 
\oplus{}^{(\mu)}_{\,\,\,\mu}V^{(0,4)}_{(1,2)}$. 

Much in the same way, by Lemma \ref{equiv} and \eqref{db2}, 
\begin{align*}
\eth_2(v) &= \ebox{\frac{c_3}{c_2}}\big(\ebox{\frac{1}{[2]_2^2[3]_2}} X_1F_2^2\ltr v 
- \ebox{\frac{1}{[2]_2}}X_0F_2\ltr v + X_{-1}\ltr v\big) \\
&= \ebox{\frac{c_3}{c_2}}\big(\ebox{\frac{1}{[2]_2^2[3]_2}} E_2^2E_1F_2^2\ltr v 
- \ebox{\frac{1}{[2]_2}}E_2E_1F_2\ltr v + E_1\ltr v\big) \\
&= \ebox{\frac{c_3}{c_2}}\big(\ebox{\frac{1}{[2]_2^2[3]_2}} E_2^2F_2^2
- \ebox{\frac{1}{[2]_2}}E_2F_2 + 1\big)E_1\ltr v \\
&=   \ebox{\frac{c_3}{c_2}}\, \pr_{(3,0)}\hsp\circ\hsp X_{-1}(v), 
\end{align*}
where we used the fact that $F_2^3 E_1 \ltr v = 0$ for 
$v\in{}^{(\mu)}_{\,\,\,\mu}V^{(1,2)}_{(1,2)}$ 
which implies that  $E_1 \ltr v$ belongs to 
${}^{(\mu)}_{\,\,\,\mu}V^{(3,0)}_{(3,0)} \oplus{}^{(\mu)}_{\,\,\,\mu}V^{(2,2)}_{(2,0)} 
\oplus{}^{(\mu)}_{\,\,\,\mu}V^{(1,4)}_{(2,0)}$ 
so that \eqref{Pr0} can be applied. 
\end{proof} 

\begin{rem} \label{rem} 
The reason why the projection $\pr_{(0,2)}$ does not appear in \eqref{eth0} is that 
$F_2E_1\ltr v = E_1F_2\ltr v =0$ which implies that 
$E_1\ltr v\in {}^{(\mu)}_{\,\,\,\mu}V^{(0,2)}_{(2,-2)}$ 
is a lowest weight vector for $\Uql$,  
thus $X_1\ltr v=E_2^2E_1\ltr v$ belongs already to ${}^{(\mu)}_{\,\,\,\mu}V^{(0,2)}_{(0,2)}$. 
\end{rem} 

To determine the spectrum of the Dirac operator defined as 
the operator closure of $\db + \db^\dagger$, 
we begin with computing the eigenvalues of 
$\eth_j^\dagger\eth_j$ and $\eth_j\eth_j^\dagger$. For this, we need an explicit description 
of the Hilbert space adjoint $\eth_j^\dagger$ of $\eth_j$. Recall that the operators $\eth_j$ 
in Corollary~\ref{coreth} are described by actions of elements from $\Uqg$ 
and that the inner product on 
${}^{(\mu)}_{\,\,\,\mu}V^{\hl_k}_{\hl_k}$ is given by the Haar state. 
As mentioned in Section \ref{DN}, 
the Hilbert space adjoint with respect to the 
left regular representation of $\Uqg$ on $\OGq$ coincides with the 
conjugate element in the *-algebra $\Uqg$. 
However, we cannot simply apply the involution of $\Uqg$ to $\eth_j$ 
since, in general, these adjoint elements from $\Uqg$ do not leave the 
orthogonol complement $({}^{(\mu)}_{\,\,\,\mu}V^{\hl_k}_{\hl_k})^\bot \subset \OGq$ 
invariant. In this case, 
it will be necessary to project onto the correct weight space. The following lemma will state 
explicit formulas for the Hilbert space adjoints of $\eth_j$, $j=0,1,2$, 
using the projection from \eqref{pr2} 
and  \eqref{Pr0}. 

\begin{lem}  \label{adjoints} 
The Hilbert space adjoints of the operators $\eth_0$, $\eth_1$ and $\eth_2$ 
from Corollary~\ref{coreth} 
are given by 
\begin{align}   \label{deth0} 
&\eth_0^\dagger\,:\, {}^{(\mu)}_{\,\,\,\mu}V^{(0,2)}_{(0,2)} \lra
{}^{(\mu)}_{\,\,\,\mu}V^{(0,0)}_{(0,0)}, \quad 
\eth_0^\dagger(v)=q^2\hs\ebox{\frac{c_1}{c_0}}\hs  \pr_{(0,0)}\hsp\circ\hsp F_1F_2^2\ltr v,\\ 
& \eth_1^\dagger \,:\,  {}^{(\mu)}_{\,\,\,\mu}V^{(1,2)}_{(1,2)} \lra  \label{deth1} 
{}^{(\mu)}_{\,\,\,\mu}V^{(0,2)}_{(0,2)}, \quad  
\eth_1^\dagger(v)= q^2\hs \ebox{\frac{c_2}{c_1}}\hs 
\ebox{\frac{[3]_2-1}{[3]_2}}\hs \pr_{(0,2)}\hsp\circ\hsp F_1F_2\ltr v, \\ 
& \eth_2^\dagger \,:\, {}^{(\mu)}_{\,\,\,\mu}V^{(3,0)}_{(3,0)} \lra 
{}^{(\mu)}_{\,\,\,\mu}V^{(1,2)}_{(1,2)} , \quad 
\eth_2^\dagger(v)=  q^2\hs \ebox{\frac{c_3}{c_2}}\hs F_1\ltr v. \label{deth2} 
\end{align} 
\end{lem} 
\begin{proof}
From the fact that the left $\Uqg$-action on $\OGq$ 
defines a *-re\-pre\-sen\-ta\-tion of $\Uqg$ with respect to the inner product \eqref{iph}, 
the defining relations of $\Uqg$ and \eqref{rep1}, we obtain for all 
$v_0\in{}^{(\mu)}_{\,\,\,\mu}V^{(0,0)}_{(0,0)}$ and 
$v_1\in {}^{(\mu)}_{\,\,\,\mu}V^{(0,2)}_{(0,2)}$
\[ \label{ipV0} 
\ip{v_1}{X_1\ltr v_0}_h 
= \ip{X_1^*\ltr v_1}{v_0}_h= q^{-1}\ip{F_1K_1(F_2K_2)^2\ltr v_1}{v_0}_h 
= q^2 \ip{F_1F_2^2\ltr v_1}{v_0}_h. 
\]
Since $E_2^3F_1F_2^2\ltr v_1=F_1E_2^3F_2^2\ltr v_1=[2]_2^2\hs F_1E_2\ltr v_1=0$
by \eqref{EFFE} and \eqref{rep1}, it follows that 
$F_1F_2^2\ltr v_1\in {}^{(\mu)}_{\,\,\,\mu}V^{(0,0)}_{(0,0)} 
\oplus{}^{(\mu)}_{\,\,\,\mu}V^{(-1,2)}_{(0,0)} \oplus{}^{(\mu)}_{\,\,\,\mu}V^{(-2,4)}_{(0,0)}$, 
where the orthogonality of the decomposition results from the Peter-Weyl theorem. 
To pick the element in ${}^{(\mu)}_{\,\,\,\mu}V^{(0,0)}_{(0,0)}$ 
without changing the inner product \eqref{ipV0}, 
we apply the orthogonal projection $\pr_{(0,0)}$ from \eqref{Pr0} and obtain 
$\ip{v_1}{\eth_0v_0}_h 
=\ebox{\frac{c_1}{c_0}}\hs q^2\ip{\pr_{(0,0)}\hsp\circ\hsp F_1F_2^2\ltr v_1}{v_0}_h$
which proves \eqref{deth0}. 

From the orthogonality of the Peter-Weyl decomposition and the definition of $\pr_{(1,2)}$ 
below Equation \eqref{pr2}, it follows that
$\ip{v_2}{\pr_{(1,2)}\hsp\circ\hsp X_0\ltr v_1}_h = \ip{v_2 }{X_0\ltr v_1}_h$ 
for all $v_1\in{}^{(\mu)}_{\,\,\,\mu}V^{(0,2)}_{(0,2)}$ and 
$v_2\in{}^{(\mu)}_{\,\,\,\mu}V^{(1,2)}_{(1,2)}$. 
As in the previous paragraph, we need to combine the action of $X_0^* = F_1F_2K_1K_2$ on $v_2$ 
with the orthogonal projection onto ${}^{(\mu)}_{\,\,\,\mu}V^{(0,2)}_{(0,2)}$. Since 
$E_2^2F_1F_2\ltr v_2 =[2]_2 F_1E_2v_2=0$, we conclude that 
$F_1F_2\ltr v_2\in{}^{(\mu)}_{\,\,\,\mu}V^{(0,2)}_{(0,2)} 
\oplus{}^{(\mu)}_{\,\,\,\mu}V^{(-1,4)}_{(0,2)}$ so that we 
can apply the orthogonal projection $\pr_{(0,2)}$ from \eqref{pr2}. This yields 
$$
\ip{v_2}{\eth_1 v_1}_h = \ebox{\frac{c_2}{c_1}}\hs \ebox{\frac{[3]_2-1}{[3]_2}}\hs q^2 \hs 
\ip{\pr_{(0,2)}\hsp\circ\hsp F_1F_2\ltr v_2}{ v_1}_h 
$$
so that \eqref{deth1} holds. 

Analogously, we have for all $v_2\in{}^{(\mu)}_{\,\,\,\mu}V^{(1,2)}_{(1,2)}$ and 
$v_3\in{}^{(\mu)}_{\,\,\,\mu}V^{(3,0)}_{(3,0)}$,  
\[ \label{ipV3}
\ip{v_3}{\eth_2 v_2}_h =
\ebox{\frac{c_3}{c_2}}\hs\ip{v_3}{\pr_{(3,0)}\hsp\circ\hsp X_{-1}\ltr v_2}_h 
= \ebox{\frac{c_3}{c_2}}\hs\ip{v_3}{X_{-1}\ltr v_2}_h 
= q^2 \hs\ebox{\frac{c_3}{c_2}}\hs\ip{F_1\ltr v_3}{v_2}_h. 
\]
Since $E_2F_1\ltr v_3 = F_1E_2\ltr v_3 =0$, it follows that 
$F_1\ltr v_3\in{}^{(\mu)}_{\,\,\,\mu}V^{(1,2)}_{(1,2)}$, 
so we do not need to apply an orthogonal projection 
and consequently \eqref{ipV3} yields \eqref{deth2}. 
\end{proof} 

Our next aim is to compute the eigenvalues 
of self-adjoint operators $\eth_k^\dagger\eth_k$ and 
$\eth_k\eth_k^\dagger$ acting on the 1 or 2 dimensional Hilbert spaces 
of highest weight vectors. 
Since the projections, and therefore all operators appearing in 
Corollary \ref{coreth} and Lemma \ref{adjoints}, can be expressed by elements from $\Uqg$, 
we can compute the eigenvalues directly by acting with these operators on 
the highest weight vectors 
from Proposition \ref{hw}. 
The next lemma shows that it suffices to compute the eigenvalues of $\eth_j \eth_j^\dagger$ and 
$\eth_j^\dagger \eth_j$ on a certain collection of 1-dimensional spaces 
from Proposition \ref{hw}, 
the other eigenvalues can be deduced from these. 

\begin{lem} \label{d*d} 
a) 
In the notation of Proposition \ref{hw} 
and with $\hl_k$ defined in \eqref{lk}, let $\mu\in \N_0\times 2\N_0$ such that 
$\dim({}^{(\mu)}_{\,\,\,\mu}V^{\hl_k}_{\hl_k})
=\dim({}^{(\mu)}_{\,\,\,\mu}V^{\hl_{k+1}}_{\hl_{k+1}}) =1$.  
Then the 1-dimensional operators 
$\eth_k^\dagger \eth_k :{}^{(\mu)}_{\,\,\,\mu}V^{\hl_k}_{\hl_k}\lra 
{}^{(\mu)}_{\,\,\,\mu}V^{\hl_k}_{\hl_k}$ and 
\,$\eth_{k}\eth_{k}^\dagger :{}^{(\mu)}_{\,\,\,\mu}V^{\hl_{k+1}}_{\hl_{k+1}}\lra 
{}^{(\mu)}_{\,\,\,\mu}V^{\hl_{k+1}}_{\hl_{k+1}}$ 
have the same eigenvalue. 

b) Let $\mu\in \N_0\times 2\N_0$ and $k\in\{1,2\}$ such that 
$\dim({}^{(\mu)}_{\,\,\,\mu}V^{\hl_k}_{\hl_k})=2$. 
Assume that $\eth_{k-1}^\dagger \eth_{k-1} $ and \,$\eth_k \eth_k^\dagger$ 
have the non-zero eigenvalues 
$c^{k-1}_\mu$ and $d^{k+1}_\mu$ on ${}^{(\mu)}_{\,\,\,\mu}V^{\hl_{k-1}}_{\hl_{k-1}}$ 
and ${}^{(\mu)}_{\,\,\,\mu}V^{\hl_{k+1}}_{\hl_{k+1}}$, 
respectively. Then there exists an orthogonal basis 
$\{a,b\} \subset{}^{(\mu)}_{\,\,\,\mu}V^{\hl_k}_{\hl_k}$ such that 
$$
\eth_{k-1}\eth_{k-1}^\dagger(a) = c^{k-1}_\mu\hs a, \quad 
\eth_{k-1}\eth_{k-1}^\dagger(b) = 0,\qquad 
\eth_k^\dagger \eth_k(a)=0, \quad \eth_k^\dagger \eth_k(b)= d^{k+1}_\mu\hs b. 
$$
\end{lem} 
\begin{proof} a) 
It is known that $\spec(T^*T)\hsp\setminus\hsp\{0\} = \spec(TT^*)\hsp\setminus\hsp\{0\}$ 
for all operators $T$ from the C*-algebra 
$\rmB(\hH)$. By considering the Hilbert space $\hH:=\hH_1\oplus \hH_2$, 
it is readily seen that the same holds true for operators 
$T\in B(\hH_1,\hH_2)$. Since $\dim({}^{(\mu)}_{\,\,\,\mu}V^{\hl_k}_{\hl_k})
=\dim({}^{(\mu)}_{\,\,\,\mu}V^{\hl_{k+1}}_{\hl_{k+1}}) =1$, 
the eigenvalues of $\eth_k^\dagger \eth_k$ and $\eth_{k}\eth_{k}^\dagger$ are either both zero, 
or both are non-zero and equal. 

b) If $\dim({}^{(\mu)}_{\,\,\,\mu}V^{\hl_k}_{\hl_k})=2$, then
$\dim({}^{(\mu)}_{\,\,\,\mu}V^{\hl_{k-1}}_{\hl_{k-1}})
=\dim({}^{(\mu)}_{\,\,\,\mu}V^{\hl_{k+1}}_{\hl_{k+1}}) =1$
by Proposition \ref{hw}. Considering the 1-dimensional operators 
$\eth_{k-1}:{}^{(\mu)}_{\,\,\,\mu}V^{\hl_{k-1}}_{\hl_{k-1}} \lra 
\ran(\eth_{k-1})\subset{}^{(\mu)}_{\,\,\,\mu}V^{\hl_k}_{\hl_k}$ and 
$\eth_{k}^\dagger:{}^{(\mu)}_{\,\,\,\mu}V^{\hl_{k+1}}_{\hl_{k+1}} \lra 
\ran(\eth_{k}^\dagger)\subset{}^{(\mu)}_{\,\,\,\mu}V^{\hl_k}_{\hl_k}$ 
it follows from the same arguments as in the proof of part a) that 
$\eth_{k-1}\eth_{k-1}^\dagger(a) = c^{k-1}_\mu\hs a$ for all $a \in \ran(\eth_{k-1})$ and 
$\eth_k^\dagger \eth_k(b)= d^{k+1}_\mu\hs b$ for all $b \in \ran(\eth_{k}^\dagger)$. 
Moreover,  $\ran(\eth_{k-1})\neq \{0\}$ and $\ran(\eth_{k}^\dagger)\neq \{0\}$ 
since the eigenvalues 
$c^{k-1}_\mu$ and $d^{k+1}_\mu$ are non-zero. As \eqref{dbc} is a complex, 
we have $\ran(\eth_{k-1})\subset \ker(\eth_k)= \ran(\eth_k^\dagger)^\bot$. 
Therefore any pair of non-zero elements $a\in \ran(\eth_{k-1})$ 
and $b\in \ran(\eth_{k}^\dagger)$ 
yields the required orthogonal basis. 
\end{proof} 

The following lemma shows how the eigenvalues and eigenvectors of the 
``Laplacians'' $\eth_k^\dagger \eth_k$ and $\eth_k\eth_k^\dagger$ determine the 
eigenvalues and eigenvectors of the Dolbeault--Dirac operator $\db + \db^\dagger$.

\begin{lem} \label{lev} 
a) Let $\mu\in\N_0\times 2\N_0$ and $k\in\{ 0,1,2\}$ be given such that 
$\dim({}^{(\mu)}_{\,\,\,\mu}V^{\hl_k}_{\hl_k})=1$, 
and let $\eth_k^\dagger \eth_k(u) = c^k_\mu u$ 
for  $u\in{}^{(\mu)}_{\,\,\,\mu}V^{\hl_k}_{\hl_k}$, where $c^k_\mu>0$. 
Then there exists an ortho\-normal basis 
$\{ u_-,u_+\} \subset{}^{(\mu)}_{\,\,\,\mu}\Omega^{(0,k)} 
\oplus \db_k\big({}^{(\mu)}_{\,\,\,\mu}\Omega^{(0,k)}\big)$ 
such that $Du_- = -\csqrt{c^k_\mu}\hs u_-$ and $Du_+ = \csqrt{c^k_\mu}\hs u_+$. 

b) Let $\mu\in\N_0\times 2\N_0$ and $k\in\{1,2,3\}$ be given such that 
$\dim({}^{(\mu)}_{\,\,\,\mu}V^{\hl_k}_{\hl_k})=1$, and let $\eth_{k-1} \eth_{k-1}^\dagger(v) 
= d^k_\mu v$ 
for  $v\in{}^{(\mu)}_{\,\,\,\mu}V^{\hl_k}_{\hl_k}$, where $d^k_\mu>0$. 
Then there exists an ortho\-normal basis 
$\{ v_-,v_+\}\hsp \subset\hsp\db^\dagger_{k-1}\big({}^{(\mu)}_{\,\,\,\mu}\Omega^{(0,k)}\big) 
\oplus{}^{(\mu)}_{\,\,\,\mu}\Omega^{(0,k)}$ 
such that $Dv_-\hsp  =\hsp  -\csqrt{d^k_\mu}\hs v_-$ and $Dv_+ 
\hsp =\hsp  \csqrt{d^k_\mu}\hs v_+$. 

c)  Let $\mu\in\N_0\times 2\N_0$ and $k\in\{1,2\}$ be given such that 
$\dim({}^{(\mu)}_{\,\,\,\mu}V^{\hl_k}_{\hl_k})=2$ and assume that 
$\eth_{k-1}^\dagger \eth_{k-1}(u) = c^{k-1}_\mu u$ for 
$u\in{}^{(\mu)}_{\,\,\,\mu}V^{\hl_{k-1}}_{\hl_{k-1}}$ 
and $\eth_{k} \eth_{k}^\dagger(v) 
= d^{k+1}_\mu v$ for $v\in{}^{(\mu)}_{\,\,\,\mu}V^{\hl_{k+1}}_{\hl_{k+1}}$ with 
$c^{k-1}_\mu >0$ and $d^{k+1}_\mu >0$. Then there exists an orthonormal basis 
$\{ u_-,u_+,v_-,v_+\}$ in 
\,${}^{(\mu)}_{\,\,\,\mu}\Omega^{(0,k-1)} \oplus{}^{(\mu)}_{\,\,\,\mu}\Omega^{(0,k)} 
\oplus{}^{(\mu)}_{\,\,\,\mu}\Omega^{(0,k+1)}$ 
such that $Du_\pm = \pm\sqrt{c^{k-1}_\mu}\hs u_\pm$ and $Dv_\pm 
= \pm\sqrt{d^{k+1}_\mu}\hs v_\pm$. 
\end{lem} 
\begin{proof} 
a) Let $u\in{}^{(\mu)}_{\,\,\,\mu}V^{\hl_k}_{\hl_k}\setminus \{0\}$ be given. 
With the unitary operator $J_k$ from Lemma~\ref{equiv}, 
set $\tilde{u}:= J_k(u)\in{}^{(\mu)}_{\,\,\,\mu}\Omega^{(0,k)}$ and 
$\tilde{w}:= \db_k(\tilde{u})\in{}^{(\mu)}_{\,\,\,\mu}\Omega^{(0,k+1)}$. Then
$$
\db_k^\dagger(\tilde{w}) = J_k J_k^*\db_k^\dagger \db_k J_k(u)= J_k \eth_k^\dagger \eth_k(u)
= c^k_\mu J_k(u) =c^k_\mu\tilde{u}. 
$$
In particular, $\tilde{w}\hsp\neq\hsp 0$ since $c^k_\mu\neq 0$ and $\tilde{u}\hsp\neq\hsp 0$. 
Moreover, $\pm \csqrt{c^k_\mu}\hs\tilde{u} + \tilde{w}\hsp\neq\hsp 0$ 
since $\tilde{u}\hs \bot\hs \tilde{w}$, and 
$ \db_{k+1}(\tilde{w})=\db_{k+1} \db_k(\tilde{u})=0 
= \frac{1}{c^k_\mu} \db_{k-1}^\dagger\db_{k}^\dagger(\tilde{w})= \db_{k-1}^\dagger(\tilde{u})$ 
since \eqref{dbc} is a complex and 
$\db_{k-1}^\dagger\db_{k}^\dagger = (\db_{k}\db_{k-1})^\dagger$. 
Therefore, 
$$
D(\pm \csqrt{c^k_\mu}\hs\tilde{u} + \tilde{w}) 
= \pm \csqrt{c^k_\mu}\hs \db_{k}(\tilde{u}) + \db_k^\dagger(\tilde{w}) 
=\pm \csqrt{c^k_\mu}\hs\tilde{w} + c^k_\mu\tilde{u}
= \pm \csqrt{c^k_\mu}\big(\pm\! \csqrt{c^k_\mu}\hs\tilde{u} + \tilde{w} \big). 
$$
In particular, 
$\csqrt{c^k_\mu}\hs\tilde{u} + \tilde{w}$ and $-\csqrt{c^k_\mu}\hs\tilde{u} + \tilde{w}$ 
are orthogonal 
because they are eigenvectors corresponding to distinct eigenvalues of the symmetric operator 
$D=\db + \db^\dagger$. 
Finally, setting 
$u_\pm := \frac{1}{\|\pm \csqrt{\,c^k_\mu}\hs\tilde{u} + \tilde{w} \|}
\big(\pm \csqrt{c^k_\mu}\hs\tilde{u} + \tilde{w}\big)$ 
gives the required orthonormal basis. 

b) For $v\hsp\in\hsp{}^{(\mu)}_{\,\,\,\mu}V^{\hl_k}_{\hl_k}\setminus \{0\}$, set 
$v_\pm \hsp := \hsp \frac{1}{\|\db_{k-1}^\dagger(J_k(v)) \pm \csqrt{\hs d^k_\mu}\hs J_k(v)  \|}
\big( \db_{k-1}^\dagger J_k(v))\pm \csqrt{d^k_\mu}\hs J_k(v) \big). $
By similar arguments as in a), it can be shown that 
$\{ v_-,v_+\}$ is an othonormal basis of 
$\db^\dagger\big({}^{(\mu)}_{\,\,\,\mu}\Omega^{(0,k)}\big) 
\oplus{}^{(\mu)}_{\,\,\,\mu}\Omega^{(0,k)}$ 
such that $Dv_\pm = \pm \csqrt{d^k_\mu}\hs v_\pm$.  

c) If $\dim({}^{(\mu)}_{\,\,\,\mu}V^{\hl_k}_{\hl_k})=2$, 
then $\dim({}^{(\mu)}_{\,\,\,\mu}V^{\hl_{k-1}}_{\hl_{k-1}})
=\dim({}^{(\mu)}_{\,\,\,\mu}V^{\hl_{k+1}}_{\hl_{k+1}})=1$
by Proposition~\ref{hw}. 
Applying a) and b) to $u\in{}^{(\mu)}_{\,\,\,\mu}V^{\hl_{k-1}}_{\hl_{k-1}}\setminus \{0\}$ and 
$v\in{}^{(\mu)}_{\,\,\,\mu}V^{\hl_{k+1}}_{\hl_{k+1}}\setminus \{0\}$ gives an orthonormal bases 
$\{ u_-,u_+\} \subset{}^{(\mu)}_{\,\,\,\mu}\Omega^{(0,k-1)} 
\oplus \db_{k-1}\big({}^{(\mu)}_{\,\,\,\mu}\Omega^{(0,k-1)}\big)$
corresponding to the eigenvalues $\pm\sqrt{c^{k-1}_\mu}$ 
and an orthonormal bases  
$\{ v_-,v_+\}\hsp \subset\hsp\db^\dagger_{k}\big({}^{(\mu)}_{\,\,\,\mu}\Omega^{(0,k+1)}\big)
 \oplus{}^{(\mu)}_{\,\,\,\mu}\Omega^{(0,k+1)}$ 
corresponding to the eigenvalues $\pm\sqrt{d^{k+1}_\mu}$. 
As $\ran(\db_{k-1})\subset \ker(\db_k) \,\bot\, \ran(\db_k^\dagger)$, 
we conclude that $\{ u_-,u_+\}\bot \{ v_-,v_+\} $, hence 
the union $\{ u_-,u_+,v_-,v_+\}$ is an orthonormal basis of the 4-dimensional space 
${}^{(\mu)}_{\,\,\,\mu}\Omega^{(0,k-1)} \oplus{}^{(\mu)}_{\,\,\,\mu}\Omega^{(0,k)} 
\oplus{}^{(\mu)}_{\,\,\,\mu}\Omega^{(0,k+1)}$. 
\end{proof} 

The previous lemma allows us to deduce a complete list of 
eigenvalues and their multiplicities for the Dolbeault--Dirac operator
from the eigenvalues of the ``Laplacians'' $\eth_k^\dagger \eth_k$ and $\eth_k\eth_k^\dagger$. 
Moreover, by Lemma \ref{d*d} and  Lemma \ref{lev}.(c), it suffices to compute the eigenvalues
on 1-dimensional highest weight spaces. This can be done directly 
and without the need to compute inner products 
by applying the algebraic expressions of 
$\eth_j$ and $\eth_j^\dagger$ from Corollary \ref{coreth} and Lemma \ref{adjoints} 
to the generating vectors in Proposition~\ref{hw}. 
We pursue this strategy in Proposition \ref{ev} to compute the eigenvalues 
of  $\eth_1^\dagger \eth_1$ and $\eth_1\eth_1^\dagger$ 
on 1-dimensional highest weight spaces. 
The eigenvalues of $\eth_0^\dagger \eth_0$ and $\eth_2\eth_2^\dagger$ 
will be obtained in the next proposition by using the Casimir operators 
from Section \ref{DN}. 

\begin{prop}          \label{evC} 
For $n,l\in \N_0$,  
\begin{align}  \label{dee0} 
&\eth_0^\dagger \eth_0 (u) = \efrac{c_1^2}{c_0^2} \,
\efrac{q^2\hs [2]_2^2\,( [n+l+2]_1[n+l]_1+ [l+1]_1[l]_1)}{[3]_2} \, u, & 
&u\in {}^{\cc 2n,2l\CC}_{\,\,(2n,2l)}V^{(0,0)}_{(0,0)}, \\     \label{ede2} 
&\eth_2\eth_2^\dagger(v) =\efrac{c_3^2}{c_2^2}\, 
\efrac{q^2\hs ( [n+l+3]_1[n+l+2]_1+ [n+3]_1[n+1]_1)}{[3]_2}\, v,  & 
&v\in {}^{\cc 2n+3,2l\CC}_{\,\,(2n+3,2l)}V^{(3,0)}_{(3,0)} .  
\end{align} 
\end{prop} 
\begin{proof}
Let $u\in {}^{\cc 2n,2l\CC}_{\,\,(2n,2l)}V^{(0,0)}_{(0,0)} $ such that $\|u\|=1$. 
As $\eth_0^\dagger \eth_0 (u) \in {}^{\cc 2n,2l\CC}_{\,\,(2n,2l)}V^{(0,0)}_{(0,0)} 
=\lin\{u\}$, 
it is clear that $\eth_0^\dagger \eth_0 (u) = \lambda u$ with 
$\lambda = \ip{u}{\eth_0^\dagger \eth_0 (u)}$. 
Since $\pr_{(0,0)}$ is an orthogonal projection, we have, 
for all $v\in {}^{\cc 2n,2l\CC}_{\,\,(2n,2l)}V^{(0,2)}_{(0,2)}$,   
\[ \label{pre}
\ebox{\frac{c_0}{q^2c_1}}\ip{u}{\eth_0^\dagger(v)} 
= \ip{u}{\pr_{(0,0)}(F_1F_2^2(v))}= \ip{\pr_{(0,0)}(u)}{F_1F_2^2(v)}=\ip{u}{F_1F_2^2(v)}. 
\] 
By \eqref{EFFE} and \eqref{rep0}, $F_2E_1(u)= E_1F_2(u)= 0$, hence  
$E_1(u)$ is a lowest weight vector of weight 
$\ha_1=(2,-2)$ and consequently 
$E_1(u) \in {}^{\cc 2n,2l\CC}_{\,\,(2n,2l)}V^{(0,2)}_{(2,-2)}$.  
Applying now \eqref{rep1}  and \eqref{pre}, we get  
\[ \label{F2E2}
\ebox{\frac{c_0^2}{q^2c_1^2}}\ip{u}{\eth_0^\dagger \eth_0 (u)}  
= \ip{u}{F_1F_2^2E_2^2E_1(u)}  = [2]_2\ip{u}{F_1F_2E_2E_1(u)}   
= [2]_2^2\ip{u}{F_1E_1(u)}.   
\] 
Again by \eqref{rep0},  $\K_i^*(u)=\K_i(u)=u$ for $i=1,2$ and 
$\E_1(u) = E_2(u) =F_2^*(u)=0$. Hence  
\[ \label{ZK} 
\ip{u}{\K_i^{\pm 1} Z(u)}=\ip{u}{K_i^{\pm 1} Z(u)}  =  \ip{u}{Z(u)}, \quad i=1,2, 
\]
and
\[ \label{ZE} 
\ip{u}{Z\E_1(u)}=\ip{u}{\E_1^*Z(u)}=\ip{u}{ZE_2(u)}=\ip{u}{F_2Z(u)}=0 
\]  
for all $Z\in \Uqg$. 
Equations \eqref{EFFE}, \eqref{in}, \eqref{Es}, \eqref{F2E2}, \eqref{ZK} and \eqref{ZE} 
imply 
\begin{align}   \label{E2E2}
&\ip{u}{\E_2^*\E_2(u)} = q^4\ip{u}{F_1F_2^2E_2^2E_1(u)} 
= \ebox{\frac{q^2c_0^2}{c_1^2}}\ip{u}{\eth_0^\dagger \eth_0 (u)},   \\
&\ip{u}{\E_3^*\E_3(u)} = q^3\ip{u}{F_1F_2E_2E_1(u)}  
= \ebox{\frac{qc_0^2}{[2]_2c_1^2}}\ip{u}{\eth_0^\dagger \eth_0 (u)},   \\ 
& \ip{u}{\E_4^*\E_4(u)}  = q \ip{u}{F_1E_1(u)}  
=  \ebox{\frac{c_0^2}{q[2]_2^2c_1^2}}\ip{u}{\eth_0^\dagger \eth_0 (u)}.  \label{E4E4} 
\end{align} 
Inserting  \eqref{SC} into $\ip{u}{S(\cC)(u)}$ and using \eqref{ZK}-\eqref{E4E4} yields  
\begin{align}  \label{lSc}
\ip{u}{S(\cC)(u)} =  \efrac{q^{-2} +q^{-1} + q + q^{2}}{(q_2-q_2^{-1})^2} 
+  \ebox{\frac{c_0^2}{q^2[2]_2c_1^2}}\left([2]_2 + (q_2^{-3} + q_2^{3}) + [2]_2 \right) 
\ip{u}{\eth_0^\dagger \eth_0 (u)}. 
\end{align} 
On the other hand, $S(\cC)= \lambda_{(2n,2l)} \id$ on  
${}^{\cc 2n,2l\CC}V^{(0,0)} \subset {}^{\cc 2n,2l\CC}V\ot V^{\cc 2n,2l\CC}$. 
Choosing a unit vector 
$u_0\in {}^{\cc 2n,2l\CC}_{\,\,(2n,2l)}V \ot V^{\cc 2n,2l\CC}_{\,\,(2n,2l)}$, 
which is a highest weight vector of weight $(2n,2l)$ 
for the left regular representation of $\Uqg$, 
we have $\K_1(u_0) = q^{-2n-2l} u_0$, $\K_2(u_0) = q^{-2n} u_0$ and $E_1(u_0)=E_2(u_0)=0$,   
so that 
\[    \label{nl} 
\lambda_{(2n,2l)} = 
\ip{u}{S(\cC)(u)}= \ip{u_0}{S(\cC)(u_0)} = 
\efrac{q^{-2n-2l-2} + q^{-2n-1}  + q^{2n+1} + q^{2n+2l+2}}{(q_2-q_2^{-1})^2}. 
\] 
Straightforward calculations show that 
\[
[2]_2 + (q_2^{-3} + q_2^{3})+[2]_2 = [4]_2+[2]_2 =[3]_2[2]_2, \label{Eq32} 
\]
\[
\efrac{q^{-2n-2l-2} + q^{-2n-1}  + q^{2n+1} + q^{2n+2l+2}  -q^{-2} -q^{-1} - q 
- q^{2}}{(q_2-q_2^{-1})^2} 
= [2]_2^2 ([n\hsp+\hsp l\hsp +\hsp 2]_1[n\hsp+\hsp l]_1+ [n\hsp +\hsp 1]_1[n]_1).  \label{Eqnl} 
\]
From \eqref{lSc}-\eqref{Eqnl}, we finally obtain 
\begin{align} \nonumber 
\lambda = \ip{u}{\eth_0^\dagger \eth_0 (u)} 
&= \ebox{\frac{q^2 c_1^2}{[3]_2c_0^2}} 
\left( \ip{u}{S(\cC)(u)}  - \efrac{q^{-2} +q^{-1} + q + q^{2}}{(q_2-q_2^{-1})^2} \right) \\ 
&=  \efrac{q^2\hs [2]_2^2\hs([n+l+2]_1[n+l]_1+ [n+1]_1[n]_1)}
         {[3]_2} \,\efrac{c_1^2}{c_0^2}.       \label{ev0} 
\end{align}  

The eigenvalues of $\eth_2\eth_2^\dagger$ can be determined in the same way. 
Let $v\in {}^{\cc 2n+3,2l\CC}_{\,\,(2n+3,2l)}V^{(3,0)}_{(3,0)}$ such that $\|v\|=1$. 
Since ${}^{\cc 2n+3,2l\CC}_{\,\,(2n+3,2l)}V^{(3,0)}_{(3,0)}$ is 1-dimensional, 
we have $\eth_2\eth_2^\dagger(v) = \mu\hs v$ with 
$\mu = \ip{v}{\eth_2\eth_2^\dagger(v) }$. As in \eqref{pre}, 
\[ \label{pre3}
\ebox{\frac{c_2}{c_3}}\ip{v}{\eth_2(w)} 
= \ip{v}{\pr_{(3,0)}(E_1(w))}= \ip{\pr_{(3,0)}(v)}{E_1(w)}=\ip{v}{E_1(w)} 
\] 
for all $w\in {}^{\cc 2n+3,2l\CC}_{\,\,(2n+3,2l)}V^{(1,2)}_{(1,2)}$. 
By \eqref{deth2}, 
$\eth_2^\dagger(v)=  q^2\hs \ebox{\frac{c_3}{c_2}}\hs F_1(v) 
\in {}^{\cc 2n+3,2l\CC}_{\,\,(2n+3,2l)}V^{(1,2)}_{(1,2)}$
is a highest weight vector of weight $(1,2)$, therefore 
\[ \label{E1F1}
\ebox{\frac{c_2^2}{q^2c_3^2}}\ip{v}{\eth_2 \eth_2^\dagger (v)}  
= \ip{v}{E_1F_1(v)}  = \efrac{1}{[2]_2}\ip{v}{E_1E_2F_2F_1(v)}   
= \efrac{1}{[2]_2^2}\ip{v}{E_1E_2^2F_2^2F_1(v)}.   
\] 
Next, $K_1(v) = q^3v$, $K_2(v)=v$ and $K_i^*=K_i$ imply 
\begin{align}
&\ip{v}{K_1^{\pm 1}Z(v)} = q^{\pm 3}\ip{v}{Z(v)},\quad  \ip{v}{K_2^{\pm 1}Z(v)} 
= \ip{v}{Z(v)}, \\ 
&\ip{v}{\K_1^{\pm 1}Z(v)} = \ip{v}{\K_2^{\pm 1}Z(v)} =q^{\mp 3}\ip{v}{Z(v)},   \label{KZv}
\end{align} 
and $\F_1(v)=F_2(v) = E_2^*(v)=0$ yields 
\[ \label{ZF} 
\ip{v}{Z\F_1(v)}=\ip{v}{\F_1^*Z(v)}=\ip{v}{ZF_2(v)}=\ip{v}{E_2Z(v)}=0 
\]  
for all $Z\in \Uqg$. Furthermore, analogous to \eqref{E2E2}-\eqref{E4E4}, 
we deduce from the previous relations that 
\begin{align}   \label{F2F2}
&\ip{v}{\F_2^*\F_2(v)} = \ip{v}{K_1^{-1}K_2^{-2}E_1E_2^2F_2^2F_1(v)} 
= \ebox{\frac{[2]_2^2c_2^2}{q^5 c_3^2}}\ip{v}{\eth_2 \eth_2^\dagger  (v)},   \\
&\ip{v}{\F_3^*\F_3(v)} = \ip{v}{K_1^{-1}K_2^{-1} E_1E_2F_2F_1(v)}  
= \ebox{\frac{[2]_2c_2^2}{q^5 c_3^2}}\ip{v}{\eth_2 \eth_2^\dagger  (v)},   \\ 
& \ip{v}{\F_4^*\F_4(v)}  = \ip{v}{qK_1^{-1} E_1F_1(v)}  
=  \ebox{\frac{c_2^2}{q^4 c_3^2}}\ip{v}{\eth_2 \eth_2^\dagger (v)}.  \label{F4F4} 
\end{align} 
From \eqref{KZv}-\eqref{F4F4} and \eqref{Eq32}, it follows that 
\begin{align}  \label{vCv}
\ip{v}{\cC(v)} =  \efrac{ q + q^{2}+q^{-2}+q^{-1}}{(q_2-q_2^{-1})^2} 
+\ebox{\frac{[3]_2\hs [2]_2^2\hs  c_2^2}{q^2\hs c_3^2}} \,\ip{v}{\eth_2 \eth_2^\dagger (v)}. 
\end{align} 
To determine $\ip{v}{\cC(v)}$, we use again the fact that $\cC = c\,\id$ on 
 ${}^{\cc 2n+3,2l\CC}V\ot V^{\cc 2n+3,2l\CC}$ and evaluate $\cC$ on a 
lowest weight vector 
$v_0\in {}^{\cc 2n+3,2l\CC}_{\,\,(2n+3,2l)}V \ot V^{\cc 2n+3,2l\CC}_{(-(2n+3),-2l)}$
for which $F_1(v_0)=F_2(v_0)=0$. Assuming $\|v_0\|=1$ and inserting \eqref{C} 
into $\ip{v_0}{\cC(v_0)}$ 
gives 
\[    \label{mu} 
\ip{v}{\cC(v)} = \ip{v_0}{\cC(v_0)} = 
\efrac{q^{-2n-2l-5} + q^{-2n-4}  + q^{2n+4} + q^{2n+2l+5}}{(q_2-q_2^{-1})^2}. 
\] 
Similarly to \eqref{Eqnl}, we have 
\[
\efrac{q^{-2n-2l-5} + q^{-2n-4}  + q^{2n+4} + q^{2n+2l+5} - q 
- q^{2}-q^{-2}-q^{-1}}{(q_2-q_2^{-1})^2 } 
= [2]_2^2 ([n\!+\!l\!+\!3]_1[n\!+\!l\!+\!2]_1+ [n\!+\!3]_1[n\!+\!1]_1).  \label{Eqnl5}
\] 
Combining \eqref{vCv}-\eqref{Eqnl5} yields 
$$
\ip{v}{\eth_2 \eth_2^\dagger (v)} = 
\efrac{q^2([n+l+3]_1[n+l+2]_1+ [n+3]_1[n+1]_1)}{[3]_2}\, \efrac{c_3^2}{c_2^2},  
$$
which proves \eqref{ede2}. 
\end{proof} 

Before determining in Proposition \ref{ev} the eigenvalues of $\eth_1^\dagger \eth_1$ and 
$\eth_1\eth_1^\dagger$ on 1-dimen\-sional weight spaces 
by computing directly the action on highest weight vectors, 
we collect some auxiliary results regarding the projections 
$\pr_{(0,2)}$ and $\pr_{(1,2)}$ from \eqref{pr2} in the following lemma. 
\begin{lem}  \label{lemmprojs} 
Let $b= (u^1_3 u^2_4 - qu^2_3u^1_4)u^1_5$. 
Then, with $X_0:=E_2E_1$,  
\begin{align}   \label{Pr12uX0z1}
&\pr_{(1,2)}(bX_0(z_1)) 
= \efrac{[2]_2^{1/2}}{[2]_1} \big(  
q^{4}[2]_2 (u^1_4 u^2_5 - qu^2_4 u^1_5)u^1_5z_1 
-q^{-1/2}z_2u^1_4(u^1_5)^2\big), \\
&\pr_{(1,2)}(X_0(z_2)b) 
= \efrac{[2]_2^{1/2}}{q^{5/2}[2]_1}z_2(u^1_4u^2_5 - qu^2_4u^1_5)u^1_5,   
\label{Pr12z2}\\
&\pr_{(0,2)}\big(F_1(z_2)u^1_3(u^1_5)^2\big)= -\efrac{q^5}{[2]_1} b z_1, 
\label{prz2u13uu15}\\   
&\pr_{(0,2)}\big(F_1(z_2) (u^1_3u^2_5 - qu^2_3u^1_5)u^1_5 \big)  
= - \efrac{q^3}{[2]_1}z_2 b, \label{prz2u35u15} \\
&\pr_{(0,2)}\big((u^1_3u^2_5 - qu^2_3u^1_5)u^1_4\big)    
=  \efrac{q^{-2}}{[2]_1}\,b. \label{pru35u14} 
\end{align} 
Let $e:= (u^1_3u^2_4- qu^2_3u^1_4)(u^1_5)^2$. Then 
\begin{align}  \label{1p02}
& \pr_{(0,2)}(F_1(z_2)F_2(e)) 
=\efrac{q^4 [2]_2^{1/2} }{[2]_1} \big(z_2^2u^1_4u^1_5 
\, - \, q^2z_2(u^1_4u^2_5 - qu^2_4u^1_5)z_1 \big), \\ 
& \pr_{(0,2)}(F_2(e)F_1(z_1)) 
= \efrac{q^{-4} [2]_2^{1/2} }{[2]_1}  \big(z_2u^1_4u^1_5z_1   \label{2p02}
\,-\, (q^2\hsp +\hsp q^3) (u^1_4u^2_5 - qu^2_4u^1_5)z_1^2 \big), \\ 
&\pr_{(0,2)}(F_1F_2(e)) 
= [2]_2^{1/2}\big( (q\!+\!q^{-3})z_2u^1_4u^1_5 - 
(1\!+\! q^{-1}\! +\! q^3)(u^1_4u^2_5 - qu^2_4u^1_5)z_1\big),  \label{3p02} \\
&\pr_{(1,2)}(X_0(z_2)u^1_4u^1_5)= \efrac{q^{-1}[2]_2^{1/2}}{[2]_1}\, e,  \label{p12Xz2}  \\
&\pr_{(1,2)}((u^1_4u^2_5-qu^2_4u^1_5)X_0(z_1)) 
= -\efrac{[2]_2^{1/2}}{[2]_1}\,e. \label{p12Xz1} 
\end{align}
\end{lem} 
\begin{proof} 
The lemma is proven by direct computation using the action \eqref{lac} 
and the relations in Lemma \ref{lem3}. For the convenience of the reader, 
we show \eqref{Pr12uX0z1}. The proof of the other equations is similar and 
can be found in \cite{D}. 

First note that 
$b\in {}^{\cc 1,2 \CC}_{(1,2 )}V^{(0,2)}_{(0,2)}$ and 
$X_0(z_1) \in  X_0\big({}^{\cc 2,0\CC}_{(2,0)}V^{(0,0)}_{(0,0)}\big) 
= {}^{\cc 2,0\CC}_{(2,0)}V^{(0,2)}_{(1,0)}$ by Proposition \ref{hw}. 
Therefore,  by \eqref{ma} and \eqref{rep1}, $E_2^2(bX_0(z_1))=0$, 
\,$K_1(bX_0(z_1))= q\hs bX_0(z_1)$ and $K_2(bX_0(z_1))= q\hs bX_0(z_1)$, 
which shows that 
$bX_0(z_1)\in {}^{\cc 3,2\CC}_{\,(3,2)}V^{(1,2)}_{(1,2)} 
\oplus {}^{\cc 3,2\CC}_{\,(3,2)}V^{(0,4)}_{(1,2)}$. 
As a consequence, the orthogonal projection onto ${}^{\cc 3,2\CC}_{\,(3,2)}V^{(1,2)}_{(1,2)}$ 
is given by \eqref{pr2}. 
From \eqref{DE}, \eqref{ma}, \eqref{lac}, $E_2(b)=0$, $F_2E_2X_0(z_1) =[2]_2 X_0(z_1)$ 
by \eqref{rep1}, 
and $\efrac{[2]_2}{[4]_2} = \efrac{1}{q+q^{-1}} = \efrac{1}{[2]_1}$, 
it follows that  
$$
\pr_{(1,2)}(bX_0(z_1)) = (1-\efrac{1}{[4]_2}F_2E_2)(bX_0(z_1))  
= \efrac{q}{[2]_1} bX_0(z_1)  - \efrac{1}{[4]_2}F_2(b) E_2X_0(z_1). 
$$
Inserting $X_0(z_1)= q [2]_2^{1/2} u^1_3 u^1_5$,   
\,$E_2X_0(z_1)= -q^{3/2}  [2]_2 u^1_4 u^1_5$,  
\begin{align*}  
F_2(b)&=[2]_2^{1/2} (u^1_2 u^2_4 - qu^2_2u^1_4)u^1_5 
- q^{-1/2}[2]_2^{1/2}(u^1_3 u^2_3 - qu^2_3u^1_3)u^1_5\\
&=q^{-1}[2]_2^{1/2}(u^1_2 u^2_4 - qu^2_2u^1_4)u^1_5 + [2]_2^{1/2}(q^{-2}\!-\!q^{-1}) z_2u^1_5, 
\end{align*} 
where we applied \eqref{u1233} in the last equation, 
and $u^1_5u^1_j=q^{-1}u^1_j u^1_5$, $j=3,4$,  
we obtain
\begin{align}   \nonumber 
\pr_{(1,2)}(bX_0(z_1)) &=\efrac{q[2]_2^{1/2}}{[2]_1} (u^1_3 u^2_4 - qu^2_3u^1_4)u^1_3 (u^1_5)^2 
+ \efrac{q^{-1/2}[2]_2^{1/2}}{[2]_1}(u^1_2 u^2_4 - qu^2_2u^1_4) u^1_4(u^1_5)^2 \\
&\quad +\efrac{(q^{-3/2}-q^{-1/2})[2]_2^{1/2}}{[2]_1}z_2u^1_4(u^1_5)^2.   \label{pr12u} 
\end{align}
Since $\pr_{(1,2)}(bX_0(z_1)), \, z_2u^1_4(u^1_5)^2\in {}^{\cc 3,2 \CC}_{(3,2)}V^{(1,2)}_{(1,2)}
=\lin\{(u^1_4 u^2_5 - qu^2_4u^1_5)u^1_5z_1,\,  z_2  u^1_4(u^1_5)^2\}$, 
there exist $\ha,\,\hb\in\C$ such that 
\begin{align} \nonumber
&q(u^1_3 u^2_4 - qu^2_3u^1_4)u^1_3(u^1_5)^2 
+ q^{-1/2}(u^1_2 u^2_4 - qu^2_2u^1_4) u^1_4(u^1_5)^2 \\
&\hspace{185pt}=\ha (u^1_4 u^2_5 - qu^2_4u^1_5)u^1_5z_1 
+ \hb z_2 u^1_4(u^1_5)^2.  \label{ab} 
\end{align}
As $u^1_5z_1 = q^{-2} u^1_1(u^1_5)^2$ and 
$\OGq$ has no zero divisors, \eqref{ab} is equivalent to 
\begin{align} \nonumber
&q(u^1_3 u^2_4 - qu^2_3u^1_4)u^1_3
+ q^{-1/2}(u^1_2 u^2_4 - qu^2_2u^1_4) u^1_4 \\
&\hspace{135pt}=\ha q^{-2}(u^1_4 u^2_5 - qu^2_4u^1_5)u^1_1
+ \hb (u^1_1 u^2_5 - qu^2_1u^1_5)  u^1_4 .  \label{hahb} 
\end{align}
Acting on both sides of \eqref{hahb} with $F_1$ yields 
$$
q^{-1/2}(u^1_1 u^2_4 - qu^2_1u^1_4) u^1_4=-\hb q (u^1_1 u^2_4 - qu^2_1u^1_4) u^1_4, 
$$
hence $\hb=-q^{-3/2}$. Similarly, acting on both sides of \eqref{hahb} 
with $\frac{-q^{-1/2}}{[2]_2}E_2^2E_1$ gives 
$$
(q^{3/2}+q^{5/2} - q^{1/2}) (u^1_4 u^2_5 - qu^2_4u^1_5) u^1_4 
= (q^{-2} \ha + q^2 \hb) (u^1_4 u^2_5 - qu^2_4u^1_5) u^1_4. 
$$
Thus, with $\hb=-q^{-3/2}$, we obtain $\ha= q^4[2]_2$. 
Inserting first $\ha$ and $\hb$ into \eqref{ab} and then \eqref{ab} into \eqref{pr12u} 
proves  \eqref{Pr12uX0z1}. 
\end{proof}

\begin{prop}  \label{ev} 
For $n,l\in\N_0$, 
\begin{align} 
& \eth_1^\dagger \eth_1 (u) = \efrac{c^2_2}{c_1^2}\,
\efrac{q^2\hs  [2]_2\hs ([3]_2-1)^2\, 
(\hs[n+l+3]_1 [n+l+2]_1 + [n+2]_1[n]_1)}{[3]^2_2\hs [2]^2_1} \, u, & &
u\in { }^{\cc 2n+1,2l+2\CC}_{\,\,(2n+1,2l+2)}V^{(0,2)}_{(0,2)},   \label{dV1}\\
&\eth_1\eth_1^\dagger(v) =   \efrac{c^2_2}{c_1^2}\, 
\efrac{q^2\hs [2]_2\hs ([3]_2-1)^2\, 
(\hs [n+l+4]_1[n+l+2]_1 + [n+2]_1[n+1]_1)}{[3]_2^{2}\hs [2]_1^2}\,v, & &
v\in {}^{\cc 2n+2,2l+2\CC}_{\,\,(2n+2,2l+2)}V^{(1,2)}_{(1,2)}.   \label{d*V2}
\end{align} 
\end{prop} 
\begin{proof} Recall that $X_0=E_2E_1$, $z_1= u^1_1 u^1_5$ 
and $z_2=u^1_1 u^2_5 - qu^2_1u^1_5$.   
Applying the formulas of the left action from Section \ref{DN} gives 
\begin{align} 
\begin{split} 
& X_0(z_1) = q [2]_2^{1/2} u^1_3 u^1_5, \quad  
X_0(z_2)= q [2]_2^{1/2} (u^1_3 u^2_5 - qu^2_3u^1_5),\\
& F_1(z_1)= - q u^1_1 u^1_4, \hspace{32pt} F_1(z_2)= -q (u^1_1 u^2_4 - qu^2_1u^1_4).
\end{split} 
\end{align} 
From Lemma \ref{lem3}, it follows that 
\begin{align} 
\begin{split}  \label{XFz}
  &z_1\hs X_0(z_1) = q^2 \hs X_0(z_1)\hs  z_1, \quad 
  z_2 \hs X_0(z_2) = q^2\hs X_0(z_2)\hs z_2, \\
  & F_1(z_1) \hs z_1 = q^2 \hs z_1 \hs  F_1(z_1), \quad \  
  \,F_1(z_2) \hs z_2 = q^2 \hs z_2 \hs  F_1(z_2). 
\end{split} 
\end{align} 
Note that, by  \eqref{ma}  and \eqref{B}, $X_0(ab) = X_0(a)\hs b + a\hs X_0(b)$ and 
$F_1(ab) = F_1(a)\hs b + a\hs F_1(b)$ for all $a,b\in B$. 
Therefore \eqref{XFz} gives for $n,l\in \N$ 
\begin{align} 
\begin{split}  \label{XFzn}
& X_0(z_1^n) =  \msum{k=0}{n-1} z_1^{k} X_0(z_1) z_1^{n-k-1}  
=  \msum{k=0}{n-1} q^{2k} X_0(z_1) z_1^{n-1}  = q^{n-1}[n]_1 X_0(z_1) z_1^{n-k-1}  , \\
& X_0(z_2^l) =  \msum{k=0}{l-1} z_2^{l-k-1} X_0(z_2)  z_2^{k} 
=  \msum{k=0}{l-1} q^{-2k} z_2^{l-1} X_0(z_2)  = q^{-l+1}[l]_1 z_2^{l-k-1} X_0(z_2), 
\end{split} 
\end{align} 
and analogously  
\[    \label{F1z} 
 F_1(z_1^n) = q^{-n+1}[n]_1 F_1(z_1) z_1^{n-1}, \quad  
 F_1(z_2^l) = q^{l-1}[l]_1 z_2^{l-1} F_1(z_2). 
\] 
Let $b:= (u^1_3 u^2_4 - qu^2_3u^1_4)u^1_5$. 
Then 
$z_2^{l}b z_1^{n}$ spans ${}^{\cc 2n+1,2l +2 \CC}_{(2n+1,2l +2 )}V^{(0,2)}_{(0,2)}$ 
and $E_2(b)=0$ by Proposition \ref{hw}. 
Using again \eqref{ma} and \eqref{B}, we get  
\begin{align}   \nonumber 
&\pr_{(1,2)}\circ X_0\hs(z_2^lbz_1^n) 
= \pr_{(1,2)}\big(X_0(z_2^l)\hs K_2K_1(b)\hs z_1^n + z_2^l X_0(b) z_1^n 
+ z_2^lb X_0(z_1^n)\big)  \\
&\hspace{40pt}  = \pr_{(1,2)}\big( q^{-l+2}[l]_1  z_2^{l-1}  X_0(z_2) b z_1^n 
+  z_2^l X_0(b) z_1^n 
+ q^{n-1}[n]_1  z_2^l \hs b X_0(z_1) z_1^{n-1}\big)  \nonumber \\  \nonumber 
&\hspace{40pt}   =  q^{-l+2}[l]_1  z_2^{l-1}  \hs \pr_{(1,2)}\big(X_0(z_2) b\big) z_1^n 
  +  z_2^l \hs \pr_{(1,2)}\big(X_0(b)\big) z_1^n \\  \label{Xzlbzn}
&\hspace{40pt}  \quad  + q^{n-1}[n]_1  z_2^l \, \pr_{(1,2)}\big(b X_0(z_1)\big) z_1^{n-1} , 
\end{align}
where we applied \eqref{XFzn} in the second equality and \eqref{pr020} in the third. 
By the formulas in \eqref{lac},   
$X_0(b)= E_2E_1((u^1_3 u^2_4 - qu^2_3u^1_4)u^1_5)
= q^{3/2} [2]_2^{1/2} (u^1_4 u^2_5 - qu^2_4u^1_5)u^1_5
\in {}^{\cc 1,2\CC}_{(1,2)}V^{(1,2)}_{(1,2)}$, so that 
$\pr_{(1,2)}(X_0(b)) = X_0(b)$. 
Inserting $\pr_{(1,2)}(b X_0(z_1))$ and  $\pr_{(1,2)}(X_0(z_2)b)$ from 
Lemma~\ref{lemmprojs} gives 
\begin{align}   \label{eth1z2b1z1} \nonumber
 \pr_{(1,2)}\hsp\circ\hsp X_0(z_2^lbz_1^n)
& =  \big( q^{3/2}[2]_2^{1/2} + q^{-l}[l]_1\efrac{[2]_2^{1/2}}{q^{1/2}[2]_1} 
+ q^{n+3}[n]_1  \efrac{[2]_2^{3/2}}{[2]_1} \big) 
z_2^l(u^1_4u^2_5 - qu^2_4u^1_5)u^1_5 z_1^n \\ 
& \quad \,-\, q^{n-1}[n]_1  
\efrac{[2]_2^{1/2}}{q^{1/2} [2]_1}  z_2^{l+1} u^1_4(u^1_5)^2z_1^{n-1}. 
\end{align} 
Our next aim is to compute the action of $\pr_{(0,2)} \circ F_1F_2$ on \eqref{eth1z2b1z1}. 
Acting by $F_2$ on the last equation yields 
\begin{align}  \nonumber
&F_2 \big( \pr_{(1,2)}\hsp\circ\hsp X_0(z_2^lbz_1^n) \big) 
=  q^{n-2}[n]_1  \efrac{[2]_2}{[2]_1}  z_2^{l+1} u^1_3(u^1_5)^2z_1^{n-1} \\
&\hspace{70pt}- \big( q[2]_2+ q^{-l-1}[l]_1\efrac{[2]_2}{[2]_1}  
+ q^{n+3}[n]_1  \efrac{[2]_2^{2}}{q^{1/2}[2]_1} \big) 
z_2^l(u^1_3u^2_5 - qu^2_3u^1_5)u^1_5 z_1^n.   \label{ethdageth}
\end{align} 
Using \eqref{F1z} and $u^1_5 F_1(z_1) = -q^{-1} u^1_4 z_1$, we compute that    
\begin{align} 
 &F_1\big(z_2^l(u^1_3u^2_5 - qu^2_3u^1_5)u^1_5z_1^n \big) 
=  q^{l-1} [l]_1 z_2^{l-1} F_1(z_2)(u^1_3u^2_5 - qu^2_3u^1_5)u^1_5z_1^n  \nonumber  \\
&\hspace{40pt}  - z_2^l(u^1_3u^2_4 - qu^2_3u^1_4)u^1_5z_1^n
- (q^{-1}+q^{-n-2} [n]_1 )z_2^l(u^1_3u^2_5 - qu^2_3u^1_5)u^1_4z_1^n  \label{F1z2u35z1}
\end{align}
and similarly
\begin{align}  \nonumber 
F_1(z_2^{l+1}u^1_3(u^1_5)^2z^{n-1}_1) 
& = q^{l}[l\!+\!1]_1z_2^{l} F_1(z_2) u^1_3(u^1_5)^2z_1^{n-1}  \\
& \quad -( q^{-1}[2]_1 + q^{-n-2} [n\!-\!1]_1)z_2^{l+1}u^1_3u^1_4u^1_5z^{n-1}_1. 
\label{F1z2uu15z1}
\end{align}
From $u^1_3u^1_4u^1_5 \sim F_2((u^1_4)^2u^1_5)$, it follows that  
$u^1_3u^1_4u^1_5 \in { }^{\cc 3,0 \CC}_{(3,0)}V^{(-1,4)}_{(0,2)}$.  
This and the definition of $\pr_{(0,2)}$ in \eqref{pr2}  imply 
$\pr_{(0,2)}(u^1_3u^1_4u^1_5)=0$.  
By applying $\pr_{(0,2)}$ to \eqref{F1z2uu15z1}, using \eqref{pr020}, 
and inserting \eqref{prz2u13uu15} 
and the last relation, we obtain 
\[    \label{prconst1}
\pr_{(0,2)}\circ F_1\big(z_2^{l+1}u^1_3(u^1_5)^2z^{n-1}_1\big) 
= -\efrac{q^{l+5}\hs [l+1]_1}{[2]_1}z_2^{l} b z_1^n. 
\]
To compute the action of  $\pr_{(0,2)}$ on \eqref{F1z2u35z1}, 
recall that $z_2^l(u^1_3u^2_4 - qu^2_3u^1_4)u^1_5z_1^n$ 
belongs to  ${}^{\cc 2n+1,2l+2\CC}_{\,\,(2n+1,2l+2)}V^{(0,2)}_{(0,2)}$ 
and that $\pr_{(0,2)}$ acts on this space as the identity. 
From this and Equations \eqref{prz2u35u15} and \eqref{pru35u14}, it follows that 
\[     \label{prconst2}
\pr_{(0,2)}\big( F_1(z_2^l(u^1_3u^2_5 - qu^2_3u^1_5)u^1_5z_1^n) \big) 
=- \efrac{q^{l+2} [l]_1 + [2]_1 + q^{-3} + q^{-n-4} [n]_1 }{[2]_1} 
z_2^l b z_1^n. 
\] 
Combining (\ref{ethdageth}), (\ref{prconst1}) and (\ref{prconst2}) gives 
\begin{align}   \nonumber 
& \efrac{c_1^2}{c_2^2}\hs 
\efrac{q^{-2}[3]_2^2}{([3]_2-1)^2}\hs \eth_1^\dagger \eth_1 (z_2^lbz_1^n) 
= \pr_{(0,2)} \hsp\circ\hsp F_1F_2 \hsp\circ\hsp \pr_{(1,2)}\hsp\circ\hsp X_0 (z_2^lbz_1^n)\\
&= \efrac{[2]_2\left(-\,q^{n+l+3}[n]_1[l+1]_1 
\,\hs+\,\hs (q[2]_1+ q^{-l-1}[l]_1 +   q^{n+5/2} [2]_2 [n]_1)(q^{l+2} [l]_1 
+ [2]_1 + q^{-3} + q^{-n-4} [n]_1) \right)}
{[2]_1^2} \, z_2^lbz_1^n   \nonumber \\
& = \efrac{ [2]_2\,(\hs[n+l+3]_1 [n+l+2]_1 + [n+2]_1[n]_1)}{[2]^2_1} \, z_2^lbz_1^n, 
\end{align} 
where the last equation follows by elementary calculations, see \cite{D}.  
This shows \eqref{dV1}. 

Equation \eqref{d*V2} can be proven in the same way. Let 
$e:= (u^1_3 u^2_4 - qu^2_3u^1_4)(u^1_5)^2$. Then, by  Proposition \ref{hw}, 
${}^{\cc 2n+2,2l+2\CC}_{\,\,(2n+2,2l+2)}V^{(1,2)}_{(1,2)}= \lin \{z_2^l e z_1^n\}$.  
Since $F_2(z^l)=F_2(z_1^n)=0$, we obtain from \eqref{ma}, \eqref{lac} and \eqref{F1z} that 
$$
F_1F_2(z_2^l e z_1^n) 
= q^{l-1}[l]_1 z_2^{l-1} F_1(z_2)F_2(e) z_1^n \hsp +\hsp  z_2^l F_1F_2(e) z_1^n 
\hsp +\hsp  q^{-n-1}[n]_1  z_2^l F_2(e) F_1(z_1) z_1^{n-1}\hsp.  
$$ 
Applying \eqref{pr020} and inserting \eqref{1p02}-\eqref{3p02} gives 
\begin{align}  \label{pFFe} 
 \pr_{(0,2)} \hsp\circ\hsp F_1F_2(z_2^l e z_1^n)  
 &= \efrac{[2]_2^{1/2}\left(q^{l+3} [l]_1\hs +\hs (q+q^{-3}) [2]_1 \hs 
 +\hs q^{-n-5}[n]_1\right) }{[2]_1}  
    z_2^{l+1}u^1_4u^1_5 z_1^n  \\ \nonumber 
&- \efrac{[2]_2^{1/2}\left(q^{l+5} [l]_1\hs +\hs (1+q^{-1}+ q^3) [2]_1 
\hs +\hs (1+q)q^{-n-3}[n]_1\right) } {[2]_1}  z_2^l(u^1_4u^2_5 - qu^2_4u^1_5)z_1^{n+1}. 
 \end{align} 
Next, using \eqref{ma}, \eqref{lac}, \eqref{B}, \eqref{XFzn} and 
$E_2(u^1_4u^1_5) = X_0(u^1_4u^1_5)= 0$, we compute that 
\begin{align*}  
X_0(z_2^{l+1}u^1_4u^1_5z_1^n) 
&= X_0(z_2^{l+1}) K_2K_1(u^1_4u^1_5)z_1^n  + z_2^{l+1}u^1_4u^1_5X_0(z_1^n) \\  \nonumber 
& = q^{-l+1}[l+1]_1 z_2^{l} X_0(z_2) u^1_4u^1_5z_1^n  
+ q^{n-1}[n]_1z_2^{l+1}u^1_4u^1_5X_0(z_1) z_1^{n-1}. 
\end{align*}
Since $u^1_4u^1_5X_0(z_1) \sim u^1_3u^1_4(u^1_5)^2 \sim   F_2F_1^2((u^1_5)^4) 
\in \up{\cc 4,0\CC}V \ot V^{\cc 0,4\CC}$ 
and ${}^{\cc 4,0\CC}_{\,\,(4,0)}V^{(1,2)}_{(1,2)} =\{0\}$ by Proposition~\ref{br}, 
we have $ \pr_{(1,2)} (u^1_4u^1_5X_0(z_1))=0$. 
Hence, by \eqref{pr020} and \eqref{p12Xz2}, 
\[  \label{1pX0}
\pr_{(1,2)}\circ X_0(z_2^{l+1}u^1_4u^1_5z_1^n) 
=  \efrac{q^{-l}[2]_2^{1/2} [l+1]_1}{[2]_1}\hs z_2^l e z_1^n. 
\]
Similarly, \eqref{ma}, \eqref{lac}, \eqref{B}, \eqref{XFzn} and 
$E_2(u^1_4u^2_5 - qu^2_4u^1_5) = X_0(u^1_4u^2_5 - qu^2_4u^1_5)= 0$ imply 
\begin{align*}  
&X_0(z_2^{l}(u^1_4u^2_5\hsp -\hsp qu^2_4u^1_5)z_1^{n+1}) 
= qX_0(z_2^{l}) (u^1_4u^2_5\hsp -\hsp qu^2_4u^1_5)z_1^{n+1}  
+ z_2^{l}(u^1_4u^2_5\hsp -\hsp qu^2_4u^1_5)X_0(z_1^{n+1}) \\  \nonumber 
& = q^{-l+2}[l]_1 z_2^{l} X_0(z_2)(u^1_4u^2_5\hsp -\hsp qu^2_4u^1_5)z_1^{n+1}  
+ q^{n}[n\hsp +\hsp 1]_1z_2^{l}(u^1_4u^2_5\hsp -\hsp qu^2_4u^1_5)X_0(z_1) z_1^{n}. 
\end{align*}
Note that $X_0(z_2)(u^1_4u^2_5\hsp -\hsp qu^2_4u^1_5) \sim 
(u^1_3u^2_5\hsp -\hsp qu^2_3u^1_5)(u^1_4u^2_5\hsp -\hsp qu^2_4u^1_5) \sim 
F_2((u^1_4u^2_5\hsp -\hsp qu^2_4u^1_5)^2)$ belongs to $\up{\cc 4,0\CC}V \ot V^{\cc 0,4\CC}$, 
where we used the $q$-commutation relation between the weight vectors 
$(u^1_3u^2_5\hsp -\hsp qu^2_3u^1_5)$ and $(u^1_4u^2_5\hsp -\hsp qu^2_4u^1_5)$ obtained by 
setting $a=4$ in \eqref{z2Ez2} and acting on both sides by $E_2E_1$.  
Since ${}^{\cc 0,4\CC}_{\,\,(0,4)}V^{(1,2)}_{(1,2)} =\{0\}$ by Proposition~\ref{br}, 
we conclude that 
$\pr_{(1,2)}(X_0(z_2)(u^1_4u^2_5\hsp -\hsp qu^2_4u^1_5) )=0$. 
Thus, by \eqref{pr020} and \eqref{p12Xz1},  
\[ \label{2pX0}
\pr_{(1,2)}\circ X_0(z_2^{l}(u^1_4u^2_5\hsp -\hsp qu^2_4u^1_5)z_1^{n+1}) =
-\efrac{q^{n}[2]_2^{1/2}[n +1]_1}{[2]_1}\, z_2^{l} e z_1^{n}. 
\] 
Combining \eqref{pFFe}-\eqref{2pX0} yields 
$$
 \efrac{c_1^2}{c_2^2}\hs 
 \efrac{q^{-2}[3]_2^2}{([3]_2-1)^2}\hs \eth_1 \eth_1^\dagger (z_2^l e z_1^n)  
= \pr_{(1,2)} \hsp\circ\hsp X_0 \hsp\circ\hsp \pr_{(0,2)} \hsp\circ\hsp F_1F_2(z_2^l e z_1^n)  
=\lambda\hs z_2^l e z_1^n
$$
with 
\begin{align} \nonumber
\lambda &= \efrac{[2]_2\left(\hspace{-0.5pt}
q^{-l}[l+1]_1(q^{l+3} [l]_1\hs +\hs (q+q^{-3}) [2]_1 \hs +\hs q^{-n-5}[n]_1)
\hs +\hs q^{n}[n +1]_1(q^{l+5} [l]_1\hs +\hs (1+q^{-1}+ q^3) [2]_1 \hs+\hs (1+q)q^{-n-3}[n]_1) 
\hspace{-0.5pt}\right) }{[2]_1^2} \\  \nonumber
&= \efrac{[2]_2(\hs [n+l+4]_1[n+l+2]_1 + [n+2]_1[n+1]_1)}{[2]_1^2} , 
\end{align}
where the last equation can easily be shown by elementary computations, see \cite{D}. 
This finishes the proof.
\end{proof}

We are now in a position to give a complete set of eigenvalues of the Dirac operator 
together with the corresponding multiplicities. 

\begin{cor} \label{spec}
Let $D$ denote the operator closure of the densly defined, symmetric operator  
$\db+ \db^\dagger$ on the domain 
$\hO^{(0,\bullet)} :=\hO^{(0,0)}\oplus \hO^{(0,1)}\oplus\hO^{(0,2)} \oplus \hO^{(0,3)}$. 
Let 
\begin{align*} 
c^0_{(2n,2l)}&:= \efrac{c_1^2}{c_0^2} \,
\efrac{q^2\hs [2]_2^2\,( [n+l+2]_1[n+l]_1+ [l+1]_1[l]_1)}{[3]_2}, & & n,l\in\N_0, \\  
c^1_{(2n+1,2l)} &:=\efrac{c^2_2}{c_1^2}\,\efrac{q^2\hs  [2]_2\hs ([3]_2-1)^2\, 
(\hs[n+l+2]_1 [n+l+1]_1 + [n+2]_1[n]_1)}{[3]^2_2\hs [2]^2_1}, & & n\in \N_0, \ l\in \N, \\
d^{2}_{(2n,2l)} &:= \efrac{c^2_2}{c_1^2}\,\efrac{q^2\hs [2]_2\hs ([3]_2-1)^2\, 
(\hs [n+l+2]_1[n+l]_1 + [n+1]_1[n]_1)}{[3]_2^{2}\hs [2]_1^2},& & n,l\in\N,\\
d^{3}_{(2n+1,2l)}&:= \efrac{c_3^2}{c_2^2}\,
\efrac{q^2\hs ( [n+l+2]_1[n+l+1]_1+ [n+2]_1[n]_1)}{[3]_2}, 
& & n\in\N,\  l\in\N_0, 
\end{align*}
denote the eigenvalues 
of $\eth_0^\dagger \eth_0$, $\eth_1^\dagger \eth_1$, $\eth_1\eth_1^\dagger$ 
and $\eth_2\eth_2^\dagger$ on 
${}^{\cc 2n,2l\CC}_{\,\,(2n,2l)}V^{(0,0)}_{(0,0)}$, 
${}^{\cc 2n+1,2l\CC}_{\,\,(2n+1,2l)}V^{(0,2)}_{(0,2)}$, 
${}^{\cc 2n,2l\CC}_{\,\,(2n,2l)}V^{(1,2)}_{(1,2)}$ and 
${}^{\cc 2n+1,2l\CC}_{\,\,(2n+1,2l)}V^{(3,0)}_{(3,0)}$, 
respectively. 
Then $D$ is self-adjoint and has spectrum 
\begin{align}
\spec(D)&= \big\{0\big\} \cup \left\{ \pm \sqrt{c^0_{(2n,2l)}}: 
n,l\hsp\in\hsp \N_0,\ n\hsp +\hsp l>0\right\} \cup 
 \left\{ \pm \sqrt{d^2_{(2n,2l)}}: n,l\hsp\in\hsp \N\right\}         \nonumber \\ \label{sD} 
&\hspace{0pt} \cup \left\{ 
\pm \sqrt{c^1_{(2n+1,2l)}}:n\hsp\in\hsp \N_0, \ l\hsp\in\hsp \N\right\}
 \cup 
\left\{ \pm \sqrt{d^3_{(2n+1,2l)}}: n\hsp\in\hsp \N,\  l\hsp\in\hsp \N_0\right\}
\end{align}
with the following multiplicities:  
\begin{align*} 
&\mbox{$\frac{(2n+ 1)(2l+1)(n+ l+ 1)(4n+ 2l+ 3)}{3}$ 
\,for the eigenvalues 
\mbox{$\pm \sqrt{\phantom{c^0_.\hspace{12pt}a}}\hspace{-29pt}{c^0_{(2n,2l)}}$} and 
\mbox{$\pm \sqrt{\phantom{c^0_.\hspace{12pt}a}}\hspace{-29pt}{d^2_{(2n,2l)}}$},} \\ 
&\mbox{$\frac{(n+1)(2l+1)(2n+ 2l+ 3)(4n+ 2l+ 5)}{3} $ \,for the eigenvalues 
\mbox{$\pm \sqrt{\phantom{c^0_.\hspace{21pt}a}}\hspace{-39pt}c^1_{(2n\hsp +\hsp 1,2l)}$}  and 
\mbox{$\pm \sqrt{\phantom{c^0_.\hspace{21pt}a}}\hspace{-39pt}d^3_{(2n\hsp +\hsp 1,2l)}$},} 
\end{align*}
for $n,l\in \N_0$ such that the corresponding eigenvalues are listed in \eqref{sD}. 
\end{cor} 
\begin{proof} 
From \eqref{orthoO}, \eqref{hO}, Lemma \ref{equiv}, Corollary \ref{coreth} 
and Lemma \ref{adjoints}, 
we conclude that $\hO^{(0,\bullet)}$ belongs to the domain of $\db^\dagger$. 
Hence $\db+ \db^\dagger$ is densely defined and symmetric.

Let $n,l\in\N$. From Proposition \ref{br}, Equations \eqref{iso0}-\eqref{iso2} 
and Lemma \ref{lev}, 
it follows that 
$$
{}^{\cc 2n,2l\CC}_{\,\,(2n,2l)}\big(
\hO^{(0,0)} \oplus \hO^{(0,1)} \oplus \hO^{(0,2)} \oplus \hO^{(0,3)}\big)
= {}^{\cc 2n,2l\CC}_{\,\,(2n,2l)}\hO^{(0,0)} \oplus {}^{\cc 2n,2l\CC}_{\,\,(2n,2l)}\hO^{(0,1)} 
\oplus {}^{\cc 2n,2l\CC}_{\,\,(2n,2l)}\hO^{(0,2)} 
$$
has an orthonormal basis given by 4 eigenvectors of $D$ corresponding to the eigenvalues 
\mbox{$\pm \sqrt{\phantom{c^0_.\hspace{12pt}a}}\hspace{-29pt}{c^0_{(2n,2l)}}$} and 
\mbox{$\pm \sqrt{\phantom{c^0_.\hspace{12pt}a}}\hspace{-29pt}{d^2_{(2n,2l)}}$}. 
By Lemma \ref{dbeq} and Equation \eqref{hO}, each of these eigenvectors 
generates an irreducible $\Uqg$-representation of highest weight $(2n,2l)$ 
consisting of eigenvectors  
corresponding to the same eigenvalue. Hence the multiplicity of this eigenvalue is given by 
$\dim\hsp\big({}^{\cc 2n,2l\CC}V\big)$. 
According to Weyl's dimension formula, 
$\dim\hsp\big({}^{\cc m,k\CC}V\big)= \frac{1}{6}(m+1)(k+1)(m+k+2)(2m+k+3)$. 
Inserting in this formula $(m,k)=(2n,2l)$ 
yields the stated multiplicities for the eigenvalues 
\mbox{$\pm \sqrt{\phantom{c^0_.\hspace{12pt}a}}\hspace{-29pt}{c^0_{(2n,2l)}}$} and 
\mbox{$\pm \sqrt{\phantom{c^0_.\hspace{14pt}a}}\hspace{-31pt}{d^2_{(2n,2l)}}$}. 
Moreover, all these eigenvectors span the vector space 
${}^{\cc 2n,2l\CC}\big(
\hO^{(0,0)} \oplus \hO^{(0,1)}\oplus \hO^{(0,2)} \oplus \hO^{(0,3)}\big)$ 
by \eqref{hO}. 

Analogously, for $n,l\in \N$, 
$$
{}^{\cc 2n+1,2l\CC}_{\,\,(2n+1,2l)}\big(
\hO^{(0,0)} \oplus \hO^{(0,1)} \oplus \hO^{(0,2)} \oplus \hO^{(0,3)}\big)
= {}^{\cc 2n+1,2l\CC}_{\,\,(2n+1,2l)}\hO^{(0,1)} 
\oplus {}^{\cc 2n+1,2l\CC}_{\,\,(2n+1,2l)}\hO^{(0,2)} 
\oplus {}^{\cc 2n+1,2l\CC}_{\,\,(2n+1,2l)}\hO^{(0,3)} 
$$
is spanned by 4 orthogonal eigenvectors corresponding to the eigenvalues 
\mbox{$\pm \sqrt{\phantom{c^0_.\hspace{23pt}a}}\hspace{-39pt}c^1_{(2n\hsp +\hsp 1,2l)}$} and 
\mbox{$\pm \sqrt{\phantom{c^0_.\hspace{23pt}a}}\hspace{-39pt}d^3_{(2n\hsp +\hsp 1,2l)}$}. 
By the same reasoning as above, each eigenvalue has the multiplicity 
$\dim\hsp\big({}^{\cc 2n+1,2l\CC}V\big)$ in 
${}^{\cc 2n+1,2l\CC}\big(\hO^{(0,0)}\oplus\hO^{(0,1)}\oplus\hO^{(0,2)}\oplus\hO^{(0,3)}\big)$, 
and applying Weyl's dimension formula gives the result. 

The cases $n=0$ or $l=0$ are treated similarly, the difference being that certain eigenvalues 
do not appear if the corresponding weight spaces have dimension 0 according to 
Proposition \ref{hw}. 
For $n=l=0$, we have ${}^{\cc 0,0\CC}\hO^{(0,\bullet)} =\lin\{1\}= \ker(D)$, so that 
the eigenvalue $c_{(0,0)}^{\,0}=0$ also satisfies the multiplicity formula. 

Finally, Proposition \ref{br} and Equation \eqref{orthoO} show that 
$$
\hO^{(0,0)}\oplus \hO^{(0,1)}\oplus\hO^{(0,2)} \oplus \hO^{(0,3)} 
=\bigoplus_{(m,k)\in\N_0\times 2\N_0} {}^{\cc m,k\CC}
\big(\hO^{(0,0)} \oplus \hO^{(0,1)}\oplus\hO^{(0,2)} \oplus \hO^{(0,3)}\big). 
$$
Hence the considered eigenvectors of $D$ form a complete orthonormal basis. 
This proves first that the restriction of $D$ to 
$\hO^{(0,0)} \oplus \hO^{(0,1)}\oplus\hO^{(0,2)} \oplus \hO^{(0,3)}$ 
is essentially self-adjoint and then that its spectrum is given by 
the set described in \eqref{sD} since 
the eigenvalues tend to infinity for $n,l\to \infty$ 
and therefore have no finite accumulation point. 
\end{proof}

In our final theorem, which presents the main results of this paper, 
we show that the Dobeault--Dirac operator 
$D:=\db + \db^\dagger$ defines an even spectral triple and  
summarize some of its properties. Recall that an even spectral triple 
$(\A, \hH, D, \gamma)$ is given by a *-al\-ge\-bra $\A$, 
a Hilbert space $\hH$, 
a faithful *-representation of $\A$ as bounded operators on $\hH$, 
a self-adjoint operator $D$ on $\hH$ with compact resolvent and 
a self-adjoint grading operator $\gamma$ satisfying the following conditions: 
$\gamma^2=\id$, $\gamma D=-D\gamma$, $\gamma a = a \gamma$, and 
$[D,a]$ is densely defined and bounded for all $a\in \A$, where we identify 
(by a slight abuse of notation) an element of $\A$ with the Hilbert space operator 
representing it, see e.\,g.\ \cite{GFV}. The spectral triple is called $0^+$-summable if 
$(1+|D|)^{-t}$ yields a trace class operator for all $t>0$. 
If \hs$\cU$ is a Hopf *-algebra and $B$ is a left $\cU$-module algebra with left action $\ltr$, 
then we say that the spectral triple is $\cU$-equivariant, 
if there exists a *-representation $\pi$ 
of $\cU$ on $\hH$ such that $\pi(X)D(h) = D\pi(X)(h)$ and 
$\pi(X)(ah)= (X_{(1)}\ltr a) \pi(X_{(2)})(h)$ 
for all $X\in \cU$, $a\in \A$ and 
$h$ from a dense domain in $\hH$ \cite{S}. 

Let us also recall that, for a Hopf algebra $\cU$ with invertible antipode, 
$\cU^{\mathrm{cop}}$ stands for the co-opposite Hopf algebra 
with the opposite coproduct $\Delta^{\!\mathrm{cop}}(X):= X_{(2)}\ot X_{(1)}$. 

\begin{thm} \label{st} 
Let $\hO^{(0,k)}$ and $\db_j$ be given as in Section \ref{secBGG}. 
With the inner product on $\hO^{(0,k)}$, $k=0,\hsp...\hs,3$,  
defined in  \eqref{hc0}-\eqref{ip2},  
let $\hH$ denote the Hilbert space closure of 
$\hO^{(0,0)}\oplus \hO^{(0,1)}\oplus\hO^{(0,2)} \oplus \hO^{(0,3)}$  
and extend the left module multiplication of the *-al\-ge\-bra $B:= \OGq^{\inv(\Uql)}$  
to a Hilbert space *-representation on $\hH$. 
Consider the operator $\db := \db_0 \oplus \db_1\oplus \db_2$ on the 
domain $\hO^{(0,\bullet)}:=\hO^{(0,0)}\oplus \hO^{(0,1)}\oplus\hO^{(0,2)} \oplus \hO^{(0,3)}$ 
and let 
$D$ denote the closure of the densely defined, symmetric operator $\db+ \db^\dagger$ 
on $\hO^{(0,\bullet)}$. 

Then $D$ is a self-adjoint operator with spectrum and multiplicities of eigenvalues given in 
Corollary \ref{spec},  
and $(B,\hH, D,\gamma)$ defines a $\Uqg^{\mathrm{cop}}$-equivariant,   
$0^+$-summ\-able, even spectral triple,  
where the grading operator 
$\gamma$ acts on $\hO^{(0,k)}$ by multiplication by $(-1)^k$ 
and the *-representation of \,$\Uqg$ on $\hH$ has been described in Lemma~\ref{dbeq}.  
\end{thm}

\begin{proof}
The fact that left module multiplication of $B$ defines a *-representation follows from 
the definition of the inner product and $h((bx)^*y) = h(x^*(b^*y))$ 
for all $b\in B$ and $x,y\in\OGq$. 
The $\Uqg^{\mathrm{cop}}$-equivariance on the domain $\hO^{(0,\bullet)}$ 
can be derived 
from Lemma~\ref{dbeq} and 
$\db^\dagger \pi_R(X)(\omega) =(\pi_R(X^*) \db)^\dagger(\omega)  
= (\db\, \pi_R(X^*) )^\dagger(\omega)  = \pi_R(X) \db^\dagger(\omega)$ 
for all  $\omega\in\hO^{(0,\bullet)}$. 
The statements about self-adjointness of $D$, its spectrum and multiplicities of eigenvalues 
are the content of  Corollary \ref{spec}.

To prove that $D$ has bounded commutators with the elements of $B$, we consider the embeddings 
$\hO^{(0,k)} \subset\OGq\ot {M^{\hl_k}}'$ with the inner products 
given by \eqref{hc0}-\eqref{ip2} 
and view $\db_k$ as an operator from 
$\hO^{(0,k)}$ to $\OGq\ot (M^{\hl_{k+1}})'$. 
Observing that 
\begin{align*}
&\OGq\ot {M^{\hl_k}}' \cong \OGq \ot \C \cong \OGq, & &k=0,3,\\
&\OGq\ot {M^{\hl_k}}' \cong \OGq \ot \C^3 \cong \overset{1}{\underset{i=-1}{\oplus}} \OGq, 
& &k=1,2,
\end{align*}
we can express $\db_k$ in the following matrix notation:  
\[ \label{dbmatrix} 
\db_0\hsp=\hsp \begin{pmatrix} X_{-1} \\ X_0 \\ X_1 \end{pmatrix}\!, \ \ 
\db_1 \hsp=\hsp \frac{1}{[3]_2}\!\begin{pmatrix} -X_{0} & \![2]_1 X_{-1} & \!0  
\\ -X_1 & \![2]_1 X_0 & \!0 \\ 0 & \!-X_1 &\! [3]_2 X_0 \end{pmatrix}\!,\ \ 
\db_2 \hsp=\hsp (\ebox{\frac{1}{[3]_2}} X_1, -X_0, X_{-1}).
\]
Furthermore, the action of $b\in B$ is represented by a diagonal matrix 
with $b$ in each diagonal entry. 
Viewing the operators in \eqref{dbmatrix} as a finite sum of matrix operators with exactly 
one non-zero matrix entry, 
it suffices to show the boundedness for commutators with matrices having 
at one place an $X_j$ and $0$ in all the others. 
We will show this for 
\[  \label{X1b} 
\left[\begin{pmatrix}  0 & 0 & 0 \\ X_1 & 0 & 0 \\ 0 & 0 & 0 \end{pmatrix}, 
\begin{pmatrix}  b & 0 & 0 \\ 0 & b & 0 \\ 0 & 0 & b \end{pmatrix} \right]
\begin{pmatrix} a_{-1}\\ a_0 \\ a_1\end{pmatrix} 
= \begin{pmatrix} 0 \\ X_1 \ltr (ba_{-1}) - b(X_1 \ltr a_{-1})\\ 0\end{pmatrix},  
\] 
where $a:=\sum_{i=-1}^1 a_{i}\ot \hu_i \in \hO^{(0,1)}\subset \OGq\ot M((0,2))'
\cong  \overset{1}{\underset{i=-1}{\oplus}}\OGq$. 
The other cases are similar (and simpler). 

As $a\in \hO^{(0,1)}$, we have 
$K_1\ltr a_j\hsp =\hsp q^{1-j} a_j, \ K_2\ltr a_j\hsp =\hsp  q^{j} a_j$ 
and $E_2\ltr a_j=a_{j+1}$ 
by \eqref{rep1} and \eqref{O1}. Moreover, $K_1\ltr b = b= K_2\ltr b$ and 
$E_2\ltr b=0$ for all $b\in B$. 
From this and \eqref{ma}, it follows that 
$$
X_1 \ltr (ba_{-1}) - b(X_1 \ltr a_{-1}) 
= (E_2^2E_1\ltr b)a_{-1} + (q+q^2)(E_2E_1\ltr b)a_{0} + q^2(E_1\ltr b)a_{1}.  
$$
Set $\|\cdot\|_h:=\sqrt{\ip{\cdot}{\cdot}_h}$ and let $\|\cdot\|_{GNS}$ 
denote the operator norm of the GNS-re\-pre\-sen\-tation 
of the compact quantum group $\OGq$ with respect to the Haar state $h$.  
Then we obtain in view of \eqref{ipc} that 
$$
\|X_1 \ltr (ba_{-1}) - b(X_1 \ltr a_{-1})\|_h 
\leq \sqrt{ 
\ebox{\frac{\|E_2^2E_1\hs\ltr \hs b \|_{GNS}^2}{c_{1,-1}}} + 
\efrac{(q+ q^2)^2\|E_2E_1\hs\ltr\hs b \|_{GNS}^2}{c_{1,0}}
+ \efrac{q^4\|E_1\hs\ltr \hs b \|_{GNS}^2}{c_{1,1}} }\,\hs  \|a\|. 
$$
Consequently, the commutator in \eqref{X1b} is bounded. 
Furthermore, the boundedness of $[\db, b^*]$ 
implies the boundedness of $[\db^\dagger, b] \subset -[\db, b^*]^\dagger$, 
and combining both results, 
we conclude that 
 $[D,b]=[\db, b] + [\db^\dagger, b]$ is bounded for all $b\in B$. 

Since $[n\hsp+\hsp l\hsp +\hsp k] = q^{-(n+l+k)} \efrac{1-\hs q^{2n+2l+2k}}{q^{-1}-\hs q}$ 
with $q^{2n+2l+k} \to 0$ for $n,l\to \infty$ and fixed $k\in\N_0$, 
we see that the eigenvalues grow exponentially fast, and since the multiplicities 
grow asymptotically no faster than $n^3l^3$, 
it follows that the trace of the positive operator $(1+|D|)^{-t}$ exists for all $t>0$. 
As a consequence, 
$D$ has compact resolvent and the spectral triple is $0^+$-summable. 

It is obvious that $\gamma$ defined by $\gamma(\omega_k) =(-1)^{k}\omega_k$ 
for all $\omega_k\in\hO^{(0,k)}$ 
satisfies $\gamma^\dagger = \gamma$, $\gamma^2=1$, 
$\gamma b = b\gamma$ and $\gamma D = - D \gamma$, hence 
$(B,\hH, D,\gamma)$ yields an even spectral triple. 
\end{proof}
\begin{rem} \label{scale}
The spectrum of the Dolbeault--Dirac operator $D$ in Theorem \ref{st} 
depends on the positive scaling factors 
$\frac{c_1}{c_0}$, $\frac{c_2}{c_1}$ and $\frac{c_3}{c_2}$, 
where the constants $c_k$ are related to the inner 
product on $\hO^{(0,k)}$, see Lemma \ref{equiv}. On the other hand, 
the Bernstein--Gelfand--Gelfand resolution \eqref{BGG} 
is unique only up to a non-zero scaling factor. 
Corollaries \ref{coreth} and \ref{spec} show that a rescaling of the Dolbeault operators 
$\db_k$ by nonzero constants will have the same effect on the spectrum of $D$, i.\,e., 
the same sets of eigenvalues will be rescaled 
by a positive constant. 
Therefore, the appearance of the scaling factors can be viewed as 
a rescaling of the inner products  
on $(0,k)$-forms, or as a rescaling of the Dolbeault operators $\db_k$ . 
\end{rem}

\subsection*{Acknowledgments} 
This work was partially supported by the CIC-UMSNH 
and by the CONACyT project A1-S-46784. 
R\'OB was supported by GA\v{C}R through the framework of the grant GA19-06357S, 
and supported by the Charles University PRIMUS grant 
\it{Spectral Noncommutative Geometry of Quantum Flag Manifolds}.

\end{document}